\documentclass[12pt,leqno]{article}

\usepackage{amsmath,amssymb,graphicx,theorem,latexsym}
\usepackage[all]{xy}

\newcommand{\op}[1]{\operatorname{#1}}
\newcommand{\opi}[1]{\operatorname{{\it #1}}}

\newcommand{\Ext}{\op{Ext}}
\newcommand{\Hom}{\op{Hom}}
\newcommand{\End}{\op{End}}
\newcommand{\Aut}{\op{Aut}}
\newcommand{\sIso}{\opi{\mycal{I}\!\!som}}
\newcommand{\sDiff}{\opi{\mycal{D}\!iff}}
\newcommand{\sHom}{\opi{\mycal{H}\!\!om}}
\newcommand{\sEnd}{\opi{\mycal{E}\!nd}}
\newcommand{\sExt}{\opi{\mycal{E}\!xt}}
\newcommand{\sAut}{\opi{\mycal{A}\!ut}}
\newcommand{\sTor}{\opi{\mycal{T}\!\!or}}

\newcommand{\bDelta}{\boldsymbol{\Delta}}
\newcommand{\Sh}{\op{\sf Sh}}
\newcommand{\sQ}[1]{\opi{{\sf Q}#1}}

\newcommand{\Omod}[1]{\sMod{\cO_{\sQ{#1}}}} 

\newcommand{\sX}{\st{X}} 
\newcommand{\sY}{\st{Y}} 
\newcommand{\an}{_{\op{an}}} 

\newcommand{\gb}[1]{{\mathbb{#1}}}	
\newcommand{\qu}[1]{{\widetilde{#1}}} 
\newcommand{\st}[1]{\mycal{#1}}     
\newcommand{\cC}{\mathcal{C}} 
\newcommand{\cO}{\mathcal{O}} 
\newcommand{\cD}{\mathcal{D}} 
\newcommand{\cE}{\mathcal{E}} 
\newcommand{\cF}{\mathcal{F}} 
\newcommand{\cG}{\mathcal{G}} 
\newcommand{\cA}{\mathcal{A}} 
\newcommand{\ID}{\op{\sf Id}} 
\newcommand{\Defst}{\op{\sf Def}} 
\newcommand{\RR}{\mathbb{R}}  
\newcommand{\LL}{\mathbb{L}}  
\newcommand{\Coh}{\op{\sf Coh}} 
\newcommand{\supp}{\op{supp}} 
\newcommand{\codim}{\op{codim}} 
\newcommand{\C}{\mathbb{C}}	
\newcommand{\Z}{\mathbb{Z}}	
\newcommand{\mult}{\op{\sf mult}} 
\newcommand{\Mod}[1]{\op{\sf Mod}(#1)}	
\newcommand{\Comod}[1]{\op{\sf Comod}(#1)}	
\newcommand{\sMod}[1]{\op{\underline{\sf Mod}(#1)}}	
\newcommand{\ord}{\op{ord}}
\newcommand{\frP}{{\mathfrak{P}}}

\newcommand{\q}[1]{^{(#1)}} 
\newcommand{\df}{\q1}	

\newcommand{\mycal}[1]{\mathcal{#1}}

\newenvironment{myequation}[1]{
\refstepcounter{subsubsection}\label{#1} 
\begin{equation}
}
{
\end{equation}
}

\def\punkt{\refstepcounter{subsubsection}
           \noindent{\bf \thesubsubsection.\ }}

\begin{document}
\title{$\star$-Quantizations of Fourier-Mukai transforms}
\author{D.~Arinkin \and J.~Block \and T.~Pantev}
\date{}
\maketitle

\tableofcontents

\

\section{Introduction}

Derived equivalences of algebraic varieties are an important and
essential topic in algebraic and complex geometry. The derived
category of a scheme encodes much of its geometry and is also flexible
enough to uncover non-trivial relations between non-isomorphic
schemes.

It is natural to study how deformations of a scheme behave with
respect to derived equivalences \cite{dima.thesis,
donagi.pantev+arinkin,BBP,oren.thesis}.  An interesting phenomenon is
that deformations of a scheme in a particular direction may correspond
to deformations to a rather exotic object on the other side of a
derived equivalence. These exotic objects are typically twisted or
non-commutative spaces, and the deformed derived equivalences often
give an interpretation of the deformations of the moduli spaces as
moduli spaces in their own right. This leads to strong geometric
statements such as the variational derived Torelli theorem
\cite{Huybrechts.Macri.Stellari,hms-formal}.

\subsection{Results}

In this paper we study the deformations of Fourier-Mukai transforms in
a general complex analytic setting; see
section~\ref{ss:other.works} for a discussion of works related to the
algebraic
setting. We
start 
with two complex manifolds $X$ and $Y$ together with a coherent
Fourier-Mukai kernel $P$ on $X\times Y$ implementing an equivalence
between the coherent derived categories of $X$ and $Y$. Then given an
arbitrary formal $*$-quantization of $X$ we construct a unique
$*$-quantization of $Y$ such that the Fourier-Mukai transform deforms
to an equivalence of the derived categories of the quantizations. For
this to hold we have to work with $*$-quantizations in the framework
of stacks of algebroids.  To simplify the exposition we consider
deformation quantizations over the ring $R =
\mathbb{C}[[\hbar]]/(\hbar^{n+1})$. The proofs immediately extend to deformation
quantizations over an arbitrary Artinian local $\C$-algebra $R$. In
the case of deformations over a complete Noetherian local $\C$-algebra
$R$ (such as the familiar case $R=\C[[\hbar]]$), Theorem~A still
holds, because it can be reduced to Artinian quotients of $R$. On the
other hand, Theorem~B requires additional techniques (such as the
results of \cite{KSarxiv}) to deal with the category of coherent
sheaves on deformation quantizations.

\

Our main results are the following two theorems (see
Theorem~\ref{thm:main} and Theorem~\ref{thm:FM})

\

\noindent
{\bfseries Theorem A.} \ {\it Suppose that the support of
$P\in\Coh(X\times Y)$ is proper over both $X$ and $Y$, and that the
integral transform $\Phi:D^b_{coh}(Y)\to D^b_{coh}(X)$ defined by $P$
is fully faithful. Then:
\begin{itemize}
\item for any $*$-quantization $\gb{X}$ of $X$
there exists a $*$-quantization $\gb{Y}$ of $Y$ and a deformation
of $P$ to an $\cO$-module $\qu{P}$ on $\gb{X}\times\gb{Y}^{op}$ (that is,
$\qu{P}$ is $R$-flat, and the reduction $\qu{P}/\hbar\qu{P}$ is
identified with $P$);  
\item the pair $(\gb{Y},\qu{P})$ is unique
up to a $1$-isomorphism, which is unique up to a unique
$2$-isomorphism. 
\end{itemize}
}

\

\medskip

\noindent
{\bfseries Theorem~B.} \ {\it Under the same assumptions we have:
\begin{itemize}
\item  the
  deformation $\qu{P}$ gives a fully faithful
  integral transform $\qu{\Phi} : D^{b}_{coh}(\gb{Y}) \to
  D^{b}_{coh}(\gb{X})$  between the coherent 
  derived categories of $\cO$-modules on $\gb{Y}$ and
  $\gb{X}$; 
\item if in addition $\Phi$ is an equivalence, then $\qu{\Phi}$ is an
  equivalence as well.
\end{itemize}
}

\subsection{Organization of the paper}

The paper is organized as follows. In Section~\ref{sec:results} we
review the setup of quantizations as stacks of algebroids and outline
the main steps in the proof of Theorem~A. (The proof of Theorem~B is
more straightforward.) There are two central results that feed into
the proof of  Theorem~A.

One of them identifies the cohomology
controlling the deformations of the pair  $(Y,P)$ to a pair
$(\gb{Y},\qu{P})$ as in Theorem~A. This is done in Section~\ref{sec:defo}. The
idea of the proof is to study deformations of an induced
$\mathcal{D}$-module associated to $P$. This $\mathcal{D}$-module
carries a comodule structure over a certain coalgebra, the necessary
framework is explained in Section~\ref{sec:defo}. 
and in full detail 

The second result is a vanishing theorem for the requisite cohomology
assuming full faithfulness of $\Phi$. We prove this result in
Section~\ref{sec:coho}, which also contains the proof of Theorem~B. As
a warm-up we start with the case of the usual Fourier-Mukai transform
for complex tori. The general proof is independent of this case.

The last section discusses an alternative approach to sheaves of  
$\cO$-modules on $*$-quantizations and the $2$-category of
$*$-quantizations. This discussion is not used anywhere else in the
paper but we include it as a unifying framework for dealing with
geometric objects over $*$-quantizations.

\subsection{Relation to other works} \label{ss:other.works}

It is instructive to compare our setup to the analogous problem in
algebraic geometry. Quantizations of schemes can be approached in two
ways. We can either consider quantizations geometrically in terms of
stacks of algebroids as above, or we can look at deformations
categorically as deformations of a dg enhancement of the category of
coherent sheaves.  Similarly a Fourier-Mukai transform for schemes can
be studied either geometrically through its kernel or categorically as
a functor between dg categories (the equivalence of the two
approaches follows from derived Morita  theory
\cite{toen-morita} and the uniqueness of enhancements  
\cite{orlov-lunts}).

This dualism fits with the generalized deformation proposal of Bondal
\cite{bondal-mpi}, Bondal-Orlov
\cite{bondal.orlov-icm}, and Kontsevich-Soibelman
\cite{ks-ncgeometry}, according to which the primary object to deform
should be a dg enchancement of the category of coherent sheaves. In
concrete terms the dualism can be expressed as follows.  Algebroid
stack deformations of a scheme are controlled by its Hochschild
cochain complex. A Fourier-Mukai equivalence between two schemes
induces
a $L_{\infty}$
quasi-isomorphism between
their Hochschild cochain
complexes \cite{keller-invariance}. This motivates
Theorem~A in the algebraic setting. This categorical interpretation
of Theorem~A is somewhat deficient since neither its statement nor its
proof give an indication on how the deformed dg categories can be
realized geometrically. A direct geometric approach to proving
Theorem~A for schemes was given in the recent work \cite{A}.

Unfortunately there are serious difficulties with the categorical
approach in the analytic context. The category of coherent sheaves on
a complex manifold $X$ is not a subtle enough invariant of $X$ (see
\cite{verbitsky}) and so its derived category does not carry enough
information. Even if we choose to work with the more flexible coherent
derived category we encounter the problem
that the derived Morita theory is not available in analytic
setting. It is possible that one can overcome this problem by
replacing the coherent derived category with a larger category, or by
simultaneously considering the derived coherent categories of $X\times
Y$ for all test analytic spaces $Y$.

In this paper, we concentrate on the geometric approach. It needs to
be suitably modified to fit the
analytic 
setting:
roughly speaking, the class of deformations must be restricted to
match the local nature of function theory. This restriction (which is
vacuous
in the
algebraic setting) is
encoded in the notion of a $*$-quantization described in
Section~\ref{ssec:star_quantizations}.  The class of $*$-quantizations
is very natural from the geometric point of view; it is equivalent to
the conditions imposed in
\cite{NT-holo,PS,KS,KSarxiv}.

\noindent
{\bfseries Note:} An interesting algebraic
approach to
Theorem~A was suggested by an anonymous referee. The idea is to
utilize the natural analog of  Theorem~A for deformations of dg
algebras \cite{keller-invariance}. Concretely let  $A$ and $B$ be 
$\mathbb{C}$-dg algebras and let $P$ be $A\otimes B^{\op{op}}$-dg
module such that the induced map 
\[
\phi_{P} : B \to \mathbb{R}\op{Hom}_{A}(P,P)
\]
is a quasi-isomorphism of $B$-modules. 
In this situation \cite{keller-invariance}
associates with $\phi_{P}$ a canonical $L_{\infty}$ morphism 
\[
\mathbb{R}\op{Hom}_{A\otimes A^{\op{op}}}(A,A) 
\to 
\mathbb{R}\op{Hom}_{B\otimes B^{\op{op}}}(B,B) 
\]
of the corresponding $L_{\infty}$ algebras of Hochschild
cochains. Since these $L_{\infty}$ algebras control the deformation
theories of $A$ and $B$ the dg algebra analogue of Theorem~A
follows. It is likely that one can combine this construction with the
observation (see e.g. \cite[Theorem~7.1.2 and
Proposition~3.6.1]{bgnt}) that the quantizations of a complex space
$X$ are classified by the Maurer-Cartan elements in the dg Lie algebra
of de Rham-Sullivan forms
$C^{\bullet}(X,CC^{\bullet}(\mathcal{O}_{X}))$. This should lead to an
order by order construction of the quantization $\gb{Y}$ in
Theorem~A. Carrying out this argument rigorously will take a
considerable effort, and we will not pursue the details here.  Instead
we follow a different route that utilizes comodule structures on
induced right $\mathcal{D}$-modules.

\subsection{Extensions and generalizations}
Theorems~A and B are naturally formulated in greater generality.  One
variant involves 
replacing $X$ and
$Y$ with $\cO^\times$-gerbes over complex manifolds, so that $\gb{X}$
and $\gb{Y}$ are $\star$-stacks (see Section~\ref{def:Quantization})
rather than $\star$-quantizations. Although our proof remains valid in
this setting, we leave this case to the reader in order to simplify
the exposition.

It is also natural to
attempt to
relax the hypotheses of Theorems~A and B by assuming that the
Fourier-Mukai kernel $P$ is an object of the derived category
$D^b_{coh}(X\times Y)$ that is not necessarily concentrated in
cohomological dimension $0$. However, such generalization requires
some technical results; we intend to return to this question in the
future.

\subsection{Notation}

\punkt For a complex manifold $M$, ${\sf Coh}(M)$ is the category of
coherent $\cO_M$-modules. Also, $D^b_{coh}(M)$ is the bounded derived
category of coherent sheaves on $M$, and $D^b_{comp}(M)\subset
D^b_{coh}(M)$ is the full subcategory of objects whose support is
compact.

\

\medskip

\punkt \ Let $X$ and $Y$ be complex manifolds and $P\in {\sf
Coh}(X\times Y)$ be a kernel object. Assume that $\supp(P)$ is
proper over $Y$.  To simplify notation we set $Z := X\times Y$, and we
write $p_{X} : Z \to X$ and $p_{Y} : Z \to Y$ for the two projections.

\

\medskip

\punkt \  Let
\[
\Phi=\Phi^P:D^b_{comp}(Y)\to D^b_{coh}(X)
\] 
be the integral transform with respect to $P$.  Explicitly, $\Phi$ is
given by
\[
\Phi(F)=\RR p_{X*}(p_Y^*(F)\otimes^\LL P)\quad(F\in D^b_{comp}(Y)).
\] 
Since $\supp(P)$ is proper over $Y$, the image of $\Phi$ is contained in
$D^b_{comp}(X)$.

\subsection{Acknowledgements}

We are very grateful to V.~Drinfeld for drawing our attention
to $\star$-operations, which (in the form of $\cD$-coalgebras)
play a central role in this paper. We would like to thank
P.~Schapira for valuable discussions about the Fourier-Mukai 
transform, and for bringing the papers
\cite{KS,KSarxiv} to our attention. We are also grateful to
A.~Yekutieli for numerous remarks and suggestions. We would also
like to thank the referees for their valuable comments and suggestions.

D.A. is a Sloan Research Fellow, and he is grateful to the Alfred
P. Sloan Foundation for the support.  J.B. was partially supported by
NSF grant DMS-1007113.  T.P. was partially supported by NSF RTG grant
DMS-0636606 and NSF grants DMS-0700446 and DMS-1001693.

\section{Formulation of main results} \label{sec:results}

Let $M$ be a complex manifold.  We fix the base ring
  of deformation to be $R = \C[\hbar]/\hbar^{n+1}$. To check
  deformations over an arbitrary finite dimensional local
  $\C$-algebra, we need to consider local $\C$-algebras where the
  maximal ideal is not necessarily principal. All the arguments seem
  to extend to this case but we will not discuss such base algebras
  explicitly.

\subsection{$\star$-quantizations}\label{ssec:star_quantizations}

Let us review the notions of $\star$-quantizations and coherent
sheaves on them. Our exposition mostly follows \cite{KSarxiv} (see
also \cite{KS}).

\

\medskip

\punkt {\bf Notation.} Denote by $\cD_M=\sDiff_M(\cO_M;\cO_M)$ the
sheaf of differential operators on $M$. More generally, let
$\sDiff_M(\cO_M,\dots,\cO_M;\cO_M)$ be the sheaf of polydifferential
operators
\[\cO_M\times\dots\times\cO_M\to\cO_M.\]
We discuss (poly)differential operators in more detail in
Section~\ref{ssec:induced_D}.

\

\medskip

\punkt {\bf Definition.} \label{defn:star_quant} {\it A {\bf
$\star$-product} on $\cO_M\otimes_\C R$ is an associative $R$-linear
product that is local and extends the {product in
  $\cO_M$}.  In
other words, the product is given by a bidifferential operator:
\begin{multline*}
f\star g=B(f,g),\quad
B=B_0+B_1\hbar+\dots+B_n\hbar^n\in\sDiff_M(\cO_M,\cO_M;\cO_M)\otimes_\C
R\\ (B_i\in\sDiff_M(\cO_M,\cO_M;\cO_M))
\end{multline*}
with $B_0(f,g)=fg$ being the {product in
  $\cO_M$}.

A {\bf neutralized $\star$-quantization} $\qu{M}$ of $M$ is a ringed
space $(M,\cO_{\qu{M}})$, where $\cO_{\qu{M}}$ is a sheaf of
$R$-algebras on $M$ that is locally isomorphic to $\cO_M\otimes_\C R$
equipped with a $\star$-product.}

\

\medskip

\punkt {\bf Remark.}  Neutralized $\star$-quantizations are usually
called simply $\star$-quantizations in the literature (see
\cite{NT-holo,yekutieli}). For us, the primary objects are
$\star$-quantizations in the sense of stacks of algebroids, as in
\cite{K,PS,L,KSarxiv}; we therefore reserve the name
`$\star$-quantizations' for these objects
{(Definition~\ref{def:*quantization})}. One can
then view neutralized $\star$-quantizations as $\star$-quantizations
with {an} additional structure: neutralization, see
Example~\ref{ex:neutralized}. Roughly speaking, for
$\star$-quantizations in the sense of stacks of algebroids, the sheaf
of algebras $\cO_{\qu{M}}$ is defined only locally; a neutralization
involves a choice of a global sheaf of algebras (which may not exist).

In the terminology of \cite{KSarxiv}, neutralized $\star$-quantizations 
are called DQ-algebras.

\

\medskip

\punkt \ Let  
\begin{myequation}{eq:alphas}
\theta_i:\cO_{U_i}\otimes_\C R\to\cO_{\qu{M}}|_{U_i}\qquad (M=\bigcup U_i)
\end{myequation}
be isomorphisms of Definition~\ref{defn:star_quant}. The isomorphisms
need not agree on double intersections, {so in
  general}
$\theta_i|_{U_i\cap U_j}\ne\theta_j|_{U_i\cap U_j}$.  However, the
discrepancy is `as small as one may hope': the composition
$$\theta_{ij}:\theta_i|^{-1}_{U_i\cap U_j}\circ\theta_j|_{U_i\cap U_j}:
\cO_{U_i\cap U_j}\otimes_\C R\to\cO_{U_i\cap U_j}\otimes_\C R$$
is given by a differential operator:
$$\theta_{ij}=1+T_1\hbar+\dots+T_n\hbar^n\in\cD_{U_i\cap
U_j}\otimes_\C R\quad (T_i\in\cD_{U_i\cap U_j}),$$ see
\cite[Proposition~4.3]{KSarxiv}. (In the sense of
Definition~\ref{defn:sheafstar}, 
$\cO_{\qu{M}}$ acquires a natural $\star$-structure.) Note in
particular that $\cO_{\qu{M}}$ is $R$-flat and that the reduction
$\cO_{\qu{M}}/\hbar\cO_{\qu{M}}$ is identified with $\cO_M$.

Note that in the definition of $\star$-product, one often requires
that $1\in\cO_M\otimes_\C R$ is the unit element. We prefer to avoid
this restriction. The definition of neutralized $\star$-quantization
is not affected by this choice: in \eqref{eq:alphas}, one can always
choose isomorphisms {$\theta_i$} and
$\star$-products on $\cO_{U_i}\otimes_\C R$ in such a way that $1$ is
the unit element. Indeed, it is clear that any $\star$-product has a
unit element (Lemma~\ref{lem:unit}); changing
{$\theta_i$}, we can ensure that $1$ is the
unit.

\

\medskip 

\punkt {\bf Lemma.} \label{lem:unit} {\it Any neutralized
$\star$-quantization $\qu{M}$ has a unit element $1\in\cO_{\qu{M}}$.}

\noindent
{\bf Proof.} Of course, this is a version of the well-known statement
that deformations of unital algebras are unital (see for instance
\cite[Section~20]{GS-units}). Let us sketch the proof.  Actually, we
only need to assume that $\cO_{\qu{M}}$ is a deformation of $\cO_M$:
that is, we assume that $\cO_{\qu{M}}$ is a flat $R$-algebra with
$\cO_{\qu{M}}/\hbar\cO_{\qu{M}}\simeq\cO_M$.

It suffices to check that the unit exists locally on $M$ since a unit
element in an associative algebra is unique. Also, it suffices to
prove existence of a right unit. Locally on $M$, $1\in\cO_M$ can be
lifted to some section $\alpha\in\cO_{\qu{M}}$. The multiplication
maps $(\bullet)\cdot \alpha : \cO_{\qu{M}} \to \cO_{\qu{M}}$ and
$\alpha \cdot (\bullet) : \cO_{\qu{M}} \to \cO_{\qu{M}}$ are bijective
(since they are bijective on the associated graded). Let now $\beta
\in \cO_{\qu{M}}$ be such that $\alpha\cdot \beta = \alpha$. Then
$(y\cdot \alpha)\cdot \beta = y\cdot \alpha$ for all local sections $y
\in \cO_{\qu{M}}$, and since every element in $\cO_{\qu{M}}$ can be
written as $y\cdot \alpha$ for some $y$, it follows that $\beta$ is a
right unit in $\cO_{\qu{M}}$. \ \hfill $\Box$

\

\medskip

\punkt \ We now proceed to define non-neutralized (`gerby')
quantizations. We first define the notion of a $\star$-stack. We view
a $\star$-stack as a quantization of an $\cO^\times$-gerbe on $M$. (In
\cite{KSarxiv}, $\star$-stacks are called DQ-algebroids, and
$\cO^\times$-gerbes are called invertible $\cO$-algebroids). If the
$\cO^\times$-gerbe is neutral, we get a
neutralized $\star$-quantization of $M$.

We use the following notation. Let $\gb{M}$ be an $R$-linear stack of
algebroids (introduced in \cite{K}, see also \cite{PS,L,KSarxiv}) over
a space $M$.  Given an open subset $U \subset M$ and two sections
$\alpha, \beta \in \gb{M}(U)$, we write
$\op{Hom}_{\gb{M}}(\alpha,\beta)$ for the $R$-module of homomorphisms
in the category $\gb{M}(U)$, and $\sHom_{\gb{M}}(\alpha,\beta)$ for
the corresponding sheaf of homomorphisms on $U$. If $\alpha=\beta$, we
write $\sEnd_{\gb{M}}(\alpha)$ instead of
$\sHom_{\gb{M}}(\alpha,\alpha)$ for the corresponding sheaf of
$R$-algebras.

\

\medskip

\punkt {\bf Definition.} \ \label{def:Quantization} {\it Let $\gb{M}$
be an $R$-linear stack of algebroids on $M$. Given an open set
$U\subset M$ and $\alpha\in\gb{M}(U)$, we denote the ringed space
$(U,\sEnd_{\gb{M}}(\alpha))$ by $\qu{U}_\alpha$. We say that $\gb{M}$
is a {\bfseries $\star$-stack} if $\qu{U}_\alpha$ is a neutralized
$\star$-quantization of $U$ for all $U$ and $\alpha$.  }

\

\medskip

\punkt {\bf Remark.} \ Note that for any $\alpha\in\gb{M}(U)$,
$\beta\in\gb{M}(V)$ (where $U,V\subset M$ are open sets), the sheaves
$\sEnd_\gb{M}(\alpha)|_{U\cap V}$ and $\sEnd_\gb{M}(\beta)|_{U\cap V}$
are locally isomorphic. Therefore, in
Definition~\ref{def:Quantization}, it suffices to check that
$\qu{(U_i)}_{\alpha_i}$ is a neutralized $\star$-quantization for some
open cover $M={\bigcup} U_i$ and for some
$\alpha_i\in\gb{M}(U_i)$.

\

\medskip

\punkt {\bf Example: neutralized quantizations.}
\label{ex:neutralized} By definition, a {\bfseries neutralization} of
a $\star$-stack $M$ is a global section $\alpha\in\gb{M}(M)$. Then
$\qu{M}_\alpha$ is a neutralized $\star$-quantization of $M$.

Conversely, given a neutralized $\star$-quantization $\qu{M}$, define
the stack of algebroids $\gb{M}$ by letting $\gb{M}_U$ be the category
of locally free right $\mathcal{O}_{\qu{M}}|_U$-modules of rank $1$.
One can easily see that $\gb{M}$ has a natural structure of a
$\star$-stack; it carries a natural neutralization
$\alpha=\cO_{\qu{M}}\in\gb{M}_M$ (the free module) such that
$\qu{M}=\qu{M}_\alpha$.  In this way, we can view neutralized
$\star$-quantizations as $\star$-stacks equipped with an additional
structure: neutralization.

\

\medskip

\punkt {\bf Example: commutative $\star$-stacks.}\ We say that a
$\star$-stack $\gb{M}$ is {\it commutative} if for every local section
$\alpha$, the sheaf of endomorphisms $\sEnd_{\gb{M}}(\alpha)$ is a
sheaf of commutative algebras.  It is easy to see that for commutative
$\gb{M}$, the sheaf $\sEnd_{\gb{M}}(\alpha)$ on an open set
$U\subset\gb{M}$ does not depend on the choice of
$\alpha\in\gb{M}(U)$. As $U$ varies, these sheaves glue to a sheaf
$\cO_{\underline{M}}$ on $M$. The ringed space
$\underline{M}=(M,\cO_{\underline{M}})$ is an $R$-deformation of $M$
as a complex analytic space, and $\gb{M}$ is simply an
$\mathcal{O}^\times$-gerbe on $\underline{M}$.  \

\medskip

\punkt Let $\gb{M}$ be a $\star$-stack over
$R=\mathbb{C}[\hbar]/\hbar^{n+1}$. For any $n'<n$,
$R'=\mathbb{C}[\hbar]/\hbar^{n'+1}$ is a quotient of $R$; it is easy
to see that $\gb{M}$ induces a $\star$-stack $\gb{M'}$ over $R'$,
which we call the {\it reduction} of $\gb{M}$ to $R'$. Explicitly, for
every open $U$, we let $(\gb{M}_U)'$ be the category whose objects are
the same as $\gb{M}_U$, while the space of morphisms between
$\alpha,\beta\in(\gb{M}_U)'$ equals
\[ 
\Gamma(U,\sHom_{\gb{M}}(\alpha,\beta)\otimes_RR').
\]
We then let $\gb{M}'$ be the stack associated with the pre-stack
$U\mapsto(\gb{M}_U)'$. It is not hard to check that $\gb{M}'$ is
a $\star$-stack over $R'$.

\

\medskip

\punkt {\bf Definition.}\label{def:*quantization} {\it A {\bfseries
$\star$-quantization of $M$} is a $\star$-stack $\gb{M}$ on $M$
together with a neutralization of the reduction of $\gb{M}$ to
$\mathbb{C}$.}

\

\medskip

\punkt {\bf Remark.} \label{rem:neutralized} Equivalently, a
$\star$-quantization can be defined as follows
(cf. \cite[{Example~2.1.2 and
    Proposition~2.1.3}]{KSarxiv}).  Let $\gb{M}$ be a
$\star$-stack. The structure of a $\star$-quantization on $\gb{M}$ is
given by an $\cO$-module $\ell$ on $\gb{M}$ such that for every open
set $U$ and every $\alpha\in\gb{M}(U)$, the action of
$\cO_{\qu{U}_\alpha}=\sEnd_{\gb{M}}(\alpha)$ on $\ell_\alpha$ factors
through $\cO_U=\sEnd_{\gb{M}}(\alpha)/\hbar\sEnd_{\gb{M}}(\alpha)$,
and that $\ell_\alpha$ is invertible as an $\cO_U$-module.

Explicitly, such $\ell$ is given as follows. For every open set $U$
and every $\alpha\in\gb{M}{(U)}$, we specify an
invertible $\cO_U$-module $\ell_\alpha$, which is functorial in
$\alpha$ and respects restriction to open subsets of $U$. Then the
sheaf $\sEnd_{\gb{M}}(\alpha)$ acts on $\ell_\alpha$. We require that
the action of $f\in\sEnd_{\gb{M}}(\alpha)$ on $\ell_\alpha$
{factors through} $\cO_U$.

Note that for any $\alpha,\beta\in\gb{M}_U$, we obtain an isomorphism
\[
\sHom_\gb{M}(\alpha,\beta)/\hbar\sHom_\gb{M}(\alpha,
\beta)\simeq\ell_\alpha^{-1}\otimes\ell_\beta
\]
that agrees with composition and restriction. In a sense,
$\ell_\alpha$ is an anti-derivative of the reduction modulo $\hbar$ of
the cocycle $\sHom_\gb{M}(\alpha,\beta)$.

\medskip

\punkt {\bf Remarks.}\ {\bf (1)} \ Let $\gb{M}$ be a
$\star$-quantization of $M$. One can easily define the {\it opposite}
$\star$-quantization $\gb{M}^{op}$; for any open $U$, the category
$\gb{M}^{op}_U$ is the opposite of $\gb{M}_U$.

{\bf (2)} \ Let $\gb{M}$, $\gb{N}$ be $\star$-quantizations of complex
manifolds $M$ and $N$, respectively. It is not hard to define the {\it
product} $\gb{M}\times\gb{N}$, which is a $\star$-quantization of
$M\times N$.

\

\medskip

\punkt We now turn to $\cO$-modules on $\star$-quantizations. We
define them as representations of stacks of algebroids (as in for
instance \cite{KSarxiv}); equivalent
approaches are discussed in Section~\ref{sec:equivariant}.

Let $\Sh_{R}$ be the stack of sheaves of $R$-modules on $M$: to an
open set $U\subset M$, it assigns the category $\Sh_R(U)$ of sheaves of
$R$-modules on $U$.

\

\medskip

\punkt {\bf Definition.} \label{def:modules_gerby} {\it Let $\gb{M}$
be a $\star$-stack. An $\cO_{\gb{M}}${\bfseries-module} is a
$1$-morphism $\qu{F} : \gb{M} \to \Sh_{R}$ of stacks of $R$-linear
categories over $M$.  }

\

\medskip

\punkt \ Let $\qu{F} : \gb{M} \to \Sh_{R}$ be an
$\cO_{\gb{M}}$-module. For any open $U\subset M$ and any
$\alpha\in\gb{M}(U)$, $\qu{F}(\alpha)$ is a sheaf of $R$-modules on
$U$.  Moreover, the sheaf of $R$-algebras $\sEnd_{\gb{M}}(\alpha)$
acts on $\qu{F}(\alpha)$.  We denote the resulting
$\cO_{\qu{U}_\alpha}$-module by $\qu{F}_{\alpha}$.  In particular, if
$\alpha \in \gb{M}(M)$ is a neutralization, the
{resulting} functor
\[
\xymatrix@R-2pc{
\left(\text{
 $\cO_\gb{M}$-modules}
\right)  \ar[r] &
\left(\text{{$\sEnd_{\gb{M}}(\alpha)$-modules}}\right) \\
\qu{F}  \ar[r] & \qu{F}_{\alpha},
}
\]
is an equivalence. 

\

\medskip

\punkt {\bfseries Definition.} \label{def:modules.coherent} {\it An
  $\cO_{\gb{M}}$-module $\qu{F}$ is {\bfseries coherent} if
  $\qu{F}_{\alpha}$ is coherent for all open $U \subset M$ and all 
$\alpha \in \gb{M}(U)$. The category of coherent
  $\cO_{\gb{M}}$-modules is denoted by $\sf{Coh}(\gb{M})$. The
  {\bfseries bounded coherent derived category} $D^{b}_{coh}(\gb{M})$
  is the derived category of complexes of $\cO_{\gb{M}}$-modules
  that have finitely many non-zero cohomology objects each of which is
  coherent. 
}

\

\medskip

\punkt {\bfseries Definition.} \label{def:modules.support} {\it The
  {\bfseries support} of an $\cO_{\gb{M}}$-module $\qu{F}$ is the
  subset in $M$ which is the 
  union of the supports of $\qu{F}_{\alpha}$ taken over all open $U
  \subset M$ and all  
$\alpha \in \gb{M}(U)$. Denote by 
$D^{b}_{comp}(\gb{M}) \subset D^{b}_{coh}(\gb{M})$
  is the full subcategory of complexes whose cohomology objects are
  compactly supported.
}

\subsection{Theorems~A and~B}

We want to study the $\star$-quantizations of Fourier-Mukai
equivalences between complex manifolds. Guided by the case of
complex tori \cite{BBP}, the general algebraic case \cite{A}, and by
the derived Morita theory \cite{toen-morita} it is  natural to
conjecture that Fourier-Mukai equivalences propagate along any
$\star$-quantization of one of the manifolds involved.

\noindent In fact in the hypotheses of Theorem~A it is not necessary
to assume properness of the support of $P$ over $X$. One should then
consider the functor $\Phi$ only on the subcategory $D^b_{comp}(Y)
\subset D^{b}_{coh}(Y)$ (see  Theorem~\ref{thm:main}
below). Theorem~\ref{thm:main} clearly implies Theorem~A but 
in fact according to Lemma~\ref{lem:compcoh} the two theorems are
equivalent if the support of $P$ is proper over $X$.

\medskip

\punkt {\bf Theorem.} \ \label{thm:main} {\it Suppose that the support
  of $P\in\Coh(X\times Y)$ is proper over $Y$, and that the integral
  transform $\Phi:D^b_{comp}(Y)\to D^b_{comp}(X)$ defined by $P$ is
  fully faithful. Then:
\begin{itemize}
\item  for any $\star$-quantization $\gb{X}$ of $X$
  there exists a $\star$-quantization $\gb{Y}$ of $Y$ and a
  deformation of $P$ to an $\cO$-module $\qu{P}$ on
  $\gb{X}\times\gb{Y}^{op}$ (that is, $\qu{P}$ is $R$-flat, and the
  reduction $\qu{P}/\hbar\qu{P}$ is identified with $P$);  
\item the pair $(\gb{Y},\qu{P})$ is unique up to a $1$-isomorphism, which
  is unique up to a unique $2$-isomorphism. 
\end{itemize}}

\

\medskip

\noindent
We also have slightly more general version of Theorem~B:
\

\smallskip

\punkt {\bf Theorem.} \ \label{thm:FM} {\it Under the
  {assumptions}  of
  Theorem~\ref{thm:main} we have:
\begin{itemize}
\item  the
  deformation $\qu{P}$ gives a fully faithful
  integral transform $\qu{\Phi} : D^{b}_{comp}(\gb{Y}) \to
  D^{b}_{comp}(\gb{X})$; 
\item  if in addition the support of  $P$ is proper over $X$,  then
  $\qu{P}$ gives a fully faithful
  integral transform $\qu{\Phi} : D^{b}_{coh}(\gb{Y}) \to
  D^{b}_{coh}(\gb{X})$; 
\item if the support of  $P$ is proper over $X$ and $\Phi$ is an
  equivalence, then \linebreak $\qu{\Phi}: D^{b}_{coh}(\gb{Y}) \to
  D^{b}_{coh}(\gb{X})$ is an 
  equivalence as well.
\end{itemize}
}

\subsection{Plan of proof of Theorem~A} \label{ssec:plan_of_proof}

\punkt Let $P$ be a coherent sheaf on $Z=X\times Y$ whose support is
proper over $Y$.  Denote by $\sDiff_{Z/X}(\cO_Z;P)$ the sheaf of
differential operators $\cO_Z\to P$ that are $\cO_X$-linear. We
consider $\sDiff_{Z/X}(\cO_Z;P)$ as a (non-coherent) $\cO_Z$-module
under right multiplications by functions. See
Section~\ref{ssec:induced_D} for details.

Note that a section $s\in P$
defines an operator $f\mapsto fs$; this yields a homomorphism
\begin{myequation}{eq:identityP}
P\to\sDiff_{Z/X}(\cO_Z;P).
\end{myequation}

\punkt \ \label{sssec:EP} Consider now in the derived category of $\cO_Z$-modules the
derived sheaf of homomorphisms
\[\RR\sHom(P,\sDiff_{Z/X}(\cO_Z;P)).\]
Its derived push-forward to $Y$ is a complex of $\cO_Y$-modules 
(or, more properly, an object of the derived category)
\[\cE(P):=\RR p_{Y,*}\RR\sHom(P,\sDiff_{Z/X}(\cO_Z;P)).\]
The homomorphism \eqref{eq:identityP} induces a map
$\cO_Z\to\RR\sHom(P,\sDiff_{Z/X}(\cO_Z;P))$, which yields a map
$\iota:\cO_Y\to\cE(P)$.

\

\medskip

\punkt {\bfseries Remark.} The derived sheaf of homomorphisms
$\RR\sHom(P,\sDiff_{Z/X}(\cO_Z;P))$
can be viewed as the derived sheaf of relative differential operators
from $P$ to itself. In fact, if $P$ is locally free, we have a quasi-isomorphism
\[\sDiff_{Z/X}(P;P)\simeq \RR\sHom(P,\sDiff_{Z/X}(\cO_Z;P)),\]
where $\sDiff_{Z/X}(P;P)$ is the sheaf of $\cO_X$-linear differential operators
from $P$ to itself. (This is a relative version of Lemma~\ref{lem:sDiff}.)
Moreover, the statement remains true for $P$ of the form
$P=i_{\Gamma *}V$, where $\Gamma\subset Z$ is a closed submanifold such that the projection
of $\Gamma$ to $X$ is submersive, $i_\Gamma:\Gamma\hookrightarrow Z$ is the embedding,
and $V$ is a locally free coherent sheaf on $\Gamma$. These hypotheses hold
for the Fourier-Mukai transform for complex tori.

\

\medskip

We separate Theorem \ref{thm:main} into two parts: one concerns
deformation theory, while the other is about integral transforms.

\

\smallskip

\punkt {\bf Theorem.} \ \label{thm:obstructions} {\it Suppose that
  $\iota:\cO_Y\to\cE(P)$ is a quasi-isomorphism.  Then the conclusion
  of Theorem~\ref{thm:main} holds: {For} any
  $\star$-quantization $\gb{X}$ of $X$ there exists a
  $\star$-quantization $\gb{Y}$ of $Y$ and a $\cO$-module $\qu{P}$ on
  $\gb{X}\times\gb{Y}^{op}$ deforming $P$.  The pair $(\gb{Y},\qu{P})$
  is unique up to a $1$-isomorphism, which is unique up to a unique
  $2$-isomorphism. }

\

\medskip

\punkt {\bf Theorem.}\ \label{thm:vanishing} {\it Let $P$ satisfy the
  hypotheses of Theorem~\ref{thm:main}. Then $\iota:\cO_Y\to\cE(P)$ is
  a quasi-isomorphism.}

\

\smallskip

\noindent Clearly, Theorems~\ref{thm:obstructions} and
\ref{thm:vanishing} imply Theorem~\ref{thm:main}. Their proofs occupy
Sections~\ref{sec:defo} and \ref{sec:coho}.

\section{Generalities} \label{sec:gen}

We prove Theorem~\ref{thm:obstructions} by replacing the $\cO$-module
$P$ with the induced (right) $\cD$-module $P_\cD$.  In addition to being a
$\cD$-module, $P_\cD$ carries a coaction of a certain coalgebra.  The
idea of the proof is that deformations of $P_\cD$ as a comodule in the
category of $\cD$-modules are easier to study than deformations of
$P$. If $P$ satisfies certain smoothness conditions, deformations of
$P_\cD$ can be interpreted in terms of $\star$-structures, as we
explain in Section~\ref{sec:tangible}.  (However, $\star$-structures
are not used in the proof.)  The approach to $\star$-structures via
$\cD$-modules is based on $\star$-pseudotensor structure of
\cite{chiralalgebras}. We are very grateful to V.~Drinfeld for drawing
our attention to this approach.

In this section, we study differential operators and $\cD$-modules on
$\star$-quantizations. We also explain the correspondence between
$\cO$-modules and comodules in the category of $\cD$-modules.

\subsection{$\cD$-modules on complex manifolds}
\label{ssec:induced_D}
Let us start by reviewing differential operators  (see
\cite{kashiwara-dmodules,bjork-dmodules}) on a (non-quantized) 
complex manifold $M$.

\

\medskip

\punkt {\bf Definition.} \label{def:diffops} {\it Let $Q$ and $P$ be
$\cO_M$-modules, and let $A:Q\to P$ be a $\C$-linear map. Given a
sequence of functions $f_0,f_1,\dots,f_N\in\cO_M$ for some $N\ge 0$,
define a sequence of $\C$-linear maps $A_k:Q\to P$ by $A_{-1}:=A$,
$A_k:=f_kA_{k-1}-A_{k-1}f_k$. We say that $A$ is a {\bfseries
differential operator} if for every point $x\in M$ and every section
$s\in Q$ defined at $x$, there exists a neighborhood $U$ of $x$ and
$N\ge 0$ such that for any open subset $V\subset U$ and any choice of
functions $f_0,\dots,f_N$ on $V$, $A_N(s|_V)$ vanishes.

Clearly, differential operators $Q\to P$ form a sheaf that we denote by
$\sDiff_M(Q;P)$.
}

\

\medskip

\punkt {\bf Remark.} Generally speaking, the bound $N$ (which could be
thought of as the order of the differential
operator $A$) depends on $s$ and on $U$. If we assume that
$Q\in\Coh(M)$, one can choose a uniform bound $N$ that does not depend
on $s$, however, it may still be local (that is, depend on $U$).

For example, consider $A:\cO_M\to P$. Then $A$ is a differential
operator if and only if it can be locally written as
\[
f\mapsto
\sum_{\alpha=(\alpha_1,\dots,\alpha_n)}
\frac{\partial^{\alpha_1}\cdots\partial^{\alpha_n}f}{\partial
  x_1^{\alpha_1}\cdots\partial x_n^{\alpha_n}} a_\alpha,
\]
for sections $a_\alpha\in P$ such that locally only finitely many of
$a_\alpha$'s are non-zero. Here $(x_1,\dots,x_n)$ is a local
coordinate system on $M$.

\medskip

\

\punkt \ It is helpful to restate Definition~\ref{def:diffops} by
splitting it into steps. Let $P$ and $Q$ be
$\cO_M$-modules.

\begin{itemize}
\item A $\C$-linear map $A:\cO_M\to P$ is a {\bfseries differential
operator of order at most $N$} if for any local functions
$f_0,f_1,\dots,f_N\in\cO_M$, we have $A_N=0$.

\item A $\C$-linear map $A:\cO_M\to P$ is a {\bfseries differential
operator} if locally there exists $N$ such that it is a differential
operator of order at most $N$.

\item A $\C$-linear map $A:Q\to P$ is a {\bfseries differential
operator} if for every $s\in Q$, the map
\[\cO_M\to P:f\mapsto A(fs)\]
is a differential operator.
\end{itemize}
It is now easy to see that a composition of differential operators is
a differential operator.

\

\medskip

\punkt Consider the sheaf of differential
operators
\[\cD_M=\sDiff_M(\cO_M;\cO_M).\]
It is a sheaf of algebras on $M$. Note that $\cD_M$ has two structures
of an $\cO_M$-module given by pre- and post-composition with the
multiplication on $\cO_M$. With respect to either structure, $\cD_M$
is a locally free $\cO_M$-module of infinite rank.

Let $P$ be an $\cO_M$-module. Consider the {\bfseries induced} right
$\cD_M$-module
\[P_\cD:=P\otimes_{\cO_M}\cD_M.\]
Here the tensor product is taken with respect to the left action of
$\cO_M$ on $\cD_M$.  We identify $P_\cD$ with the sheaf
$\sDiff_M(\cO_M;P)$ via the isomorphism
\[
P\otimes_{\cO_M}\cD_M\to\sDiff_M(\cO_M;P):p\otimes A\mapsto
A(\bullet)p.
\] 
Under this identification, the action of $\cD_M$ on
$\sDiff_M(\cO_M;P)$ comes from its action on $\cO_M$.

Note that $P$ can be reconstructed from $P_\cD$ as
$P=P_\cD\otimes_{\cD_M}\cO_M$ (here $\mathcal{O}_{M}$
is considered as a left $\mathcal{D}_{M}$-module).  The
induction functor is adjoint to the forgetful functor from
$\cD_M$-modules to $\cO_M$-modules. Both of these statements hold in
the derived sense:

\

\medskip

\punkt {\bf Lemma.} \label{lem:induced_adjunct}
{\it Let $P$ be an $\cO_M$-module.
\begin{itemize}
\item[{\bf (a)}] For any right $\cD_M$-module $F$, there are natural
  isomorphisms 
\begin{align*}
\Ext^i_{\cD_M}(P_\cD,F)=\Ext^i_{\cO_M}(P,F)\qquad&(i\ge 0)\\
\sExt^i_{\cD_M}(P_\cD,F)=\sExt^i_{\cO_M}(P,F)\qquad&(i\ge 0)
\end{align*}
of vector spaces and sheaves of vector spaces, respectively. 
\item[{\bf (b)}] For any left $\cD_M$-module $F$, there is a natural
  isomorphism 
\[\sTor_i^{\cD_M}(P_\cD,F)=\sTor_i^{\cO_M}(P,F)\]
of sheaves of $\cO_M$-modules.
In particular, for $F=\cO_M$, there is an isomorphism
\[P_\cD\otimes^\LL_{\cD_M}\cO_M=P.\]
\end{itemize}
}

\noindent
{\bf Proof.} The isomorphisms are clear for $i=0$; for $i>0$, it is
enough to resolve $F$ by injective $\cD_M$-modules in {\bf (a)} and
flat $\cD_M$-modules in {\bf (b)}. Note that injective $\cD_M$-modules
remain injective over $\cO_M$; the same holds for flat
$\cD_M$-modules.  Alternatively the
  statement follows from the fact that $\mathcal{D}_{M}$ is a locally
  free and hence flat as an $\mathcal{O}_{M}$-module.\ \hfill $\Box$

\

\medskip
  
\noindent
Let us interpret differential operators between $\cO_M$-modules in
term of the induced $\cD_M$-modules.  Directly from the definition, we
derive the following claim.
\

\smallskip

\punkt {\bf Lemma.} \label{lem:diffops_via_ind} {\it 
Let $Q$, $P$ be $\cO_M$-modules. 
\begin{itemize}
\item[{\bf (a)}] Any differential operator $A\in\sDiff_M(Q;P)$ induces
a morphism of $\cO_M$-modules
\[Q\to\sDiff_M(\cO_M;P)=P_\cD:s\mapsto A(\bullet s))\qquad(s\in Q).\]
This provides an identification $\sDiff_M(Q;P)=\sHom_{\cO_M}(Q,P_\cD)$.

\item[{\bf (b)}] Any differential operator $A\in\sDiff_M(Q;P)$ induces
a map
\[\sDiff_M(\cO_M;Q)\to\sDiff_M(\cO_M;P):B\mapsto A\circ B.\] 
This provides an identification
\[
\begin{split}
\sDiff_M(Q;P) & = \sHom_{\cD_M}(\sDiff_M(\cO_M;Q),\sDiff_M(\cO_M;P))
\\
&  = \sHom_{\cD_M}(Q_\cD,P_\cD).
\end{split}
\]
\end{itemize}
\ \hfill $\Box$
} 

\

\medskip

\punkt {\bf Remarks.} \label{rem:diffops_via_ind} \ {\bf (1)} \ Note that
$P_\cD$ has two structures of an $\cO_M$-module. Indeed, besides the
right action of $\cD_M$ (coming from its action on $\cO_M$), it has a
left action of $\cO_M$ coming from its action on $P$. Unless stated
otherwise, we ignore the latter action; in particular, the structure
of $\cO_M$-module on $P_\cD$ comes from the action of $\cO_M$ on
itself.  Equivalently, after we form the tensor product
$P\otimes_{\cO_M}\cD_M$ using the left action of $\cO_M$ on $\cD_M$,
we consider it as an $\cO_M$-module using the right action of $\cO_M$
on $\cD_M$ rather than the left action.  This is the structure used in
Lemma~\ref{lem:diffops_via_ind}(a).

{\bf (2)} \ The two statements of Lemma~\ref{lem:diffops_via_ind} are
related by the isomorphism of Lemma~\ref{lem:induced_adjunct}(b).

{\bf (3)} \ In homological calculations it is often useful to derive the
notion of a differential operator. 
In view of Lemma~\ref{lem:diffops_via_ind}, it is natural to consider 
\[
\RR\sHom_{\cO_M}(Q,P_\cD)=\RR\sHom_{\cD_M}(Q_\cD,P_\cD)
\]
as the higher derived version of the sheaf of differential operators.

\

\bigskip
  
\noindent
We have a similar formalism for polydifferential operators.

\

\smallskip

\punkt {\bf Definition.} \label{def:polydiffops}
{\it Let $P$ be an $\cO_M$-module.

\begin{itemize}
\item A $\C$-polylinear map $A:\cO_M\times\dots\times\cO_M\to P$ (of $n$ arguments)
is a {\bfseries polydifferential operator of (poly)order at most $(N_1,\dots,N_n)$} if 
it is a differential operator of order at most $N_j$ in the $j$-th argument whenever
the remaining $n-1$ arguments are fixed.

\item A $\C$-polylinear map $A:\cO_M\times\dots\times\cO_M\to P$ is a {\bfseries polydifferential operator}
if locally there exist $(N_1,\dots,N_n)$ such that $A$ is a polydifferential
operator of order at most $(N_1,\dots,N_n)$.

\item Let $P_1,\dots,P_n$ be $\cO_M$-modules. A $\C$-polylinear map
\[
A:P_1\times\dots\times P_n\to P
\]
is a {\bfseries polydifferential operator}
if for any local sections $s_1\in P_1,\dots,s_n\in P_n$, the map
\[\cO_M\times\dots\times\cO_M\to P:(f_1,\dots,f_n)\mapsto A(f_1s_1,\dots,f_ns_n)\]
is a polydifferential operator.
\end{itemize}

We denote the sheaf of polydifferential operators
by $\sDiff_M(P_1,\dots,P_n;P)$.}

\

\medskip

\

\punkt {\bf Lemma.} \label{lem:bidiffops_via_ind}
{\it Let $P$ and $P_1$ be $\cO_M$-modules and suppose that $P_1$ is coherent.
\begin{itemize}
\item[{\bf (a)}] Consider on the sheaf $\sDiff_M(P_1;P)$ the two
structures of an $\cO_M$-module coming from the action of $\cO_M$ on
$P$ and on $P_1$. Denote the resulting $\cO_M$-modules $Q$ and
$Q_1$. Then the tautological map $Q\to Q_1$ is a differential
operator, and the same holds for its inverse.
\item[{\bf (b)}] Let $P_2$ be another $\cO_M$-module. Then a map
$A:P_2\times P_1\to P$ is a bidifferential operator if and only if the
map \[P_2\to\sDiff_M(P_1;P):s\mapsto A(s,\bullet)\] is a differential
operator. The statement holds for both $\cO_M$-module structures on
$\sDiff_M(P_1;P)$.
\end{itemize}
} 

\noindent
{\bf Proof.} {\bf (a)} Let $(x_1,\dots,x_n)$ be a local coordinate
system on $M$. To show that $Q\to Q_1$ is a differential operator, we
need to prove that for any $A\in\sDiff_M(P_1;P)$, there is a formula
\[
fA=\sum_{\alpha=(\alpha_1,\dots,\alpha_n)}
a_\alpha\left(\frac{\partial^{\alpha_1}}{\partial
x_1^{\alpha_1}}\cdots \frac{\partial^{\alpha_n}}{\partial
x_n^{\alpha_n}}f\right)\qquad(f\in\cO_M)
\]
for some $a_\alpha\in\sDiff_M(P_1;P)$ such that locally, only finitely
many of $a_\alpha$ are non-zero. Such a formula exists because locally
on $M$, the order of $A$ is finite, since $P_1$ is coherent. The proof
that $Q_1\to Q$ is a differential operator is similar.

{\bf (b)} Note that by part {\bf (a)}, we have an identification
$\sDiff_M(P_2;Q)=\sDiff_M(P_2;Q_1)$. Therefore, it suffices to prove
the claim for one of the two $\cO_M$-module structures. However
the statement is obvious for the
$\cO_M$-module $Q$.  \ \hfill $\Box$

\

\medskip

\punkt {\bf Example.}\label{exm:bidiffops_via_ind} Let $P$ be an
$\cO_M$-module. By Lemma~\ref{lem:bidiffops_via_ind}, 
we get identifications
\[
\sDiff_M(\cO_M,\cO_M;P)=\sDiff_M(\cO_M;P_\cD)
= P_\cD\otimes_{\cO_M}\cD_M=P\otimes_{\cO_M}\cD_M\otimes_{\cO_M}\cD_M. 
\]
In fact, we obtain two such identifications corresponding to two
actions of $\cO_M$ on $P_\cD$. Let us write them explicitly. 

First, consider $P_\cD$ as an $\cO_M$-module using the action of
$\cO_M$ on $P$ (this $\cO_M$-module is denoted by $Q$ in
Lemma~\ref{lem:bidiffops_via_ind}). This corresponds to forming the
tensor product
\[
P\otimes_{\cO_M}\cD_M\otimes_{\cO_M}\cD_M
\] 
using the left action of $\cO_M$ on both copies
of $\cD_M$. We then obtain as isomorphism
\[
\xymatrix@R-1.5pc{
P\otimes_{\cO_M}\cD_M\otimes_{\cO_M}\cD_M
\ar[r] & \sDiff_M(\cO_M,\cO_M;P) & \\
p\otimes A_1\otimes A_2 \ar@{|->}[r] &  
A_2(\bullet)A_1(\bullet)p & (p\in P; A_1,A_2\in\cD_M).
}
\]
Now consider $P_\cD$ as an $\cO_M$-module using the action of $\cO_M$
on itself (this is denoted by $Q_1$ in
Lemma~\ref{lem:bidiffops_via_ind}). This corresponds to forming the
tensor product $P\otimes_{\cO_M}\cD_M\otimes_{\cO_M}\cD_M$ using the
structure of $\cO_M$-bimodule on the first copy of $\cD_M$.  This
provides an isomorphism
\[
\xymatrix@R-1.5pc{
P\otimes_{\cO_M}\cD_M\otimes_{\cO_M}\cD_M\ar[r] &
\sDiff_M(\cO_M,\cO_M;P) & \\
p\otimes A_1\otimes A_2 \ar@{|->}[r] & 
A_1(\bullet A_2(\bullet))p & (p\in P; A_1,A_2\in\cD_M).
}
\]
Of course, similar formulas apply to the identification 
\[
P\otimes_{\cO_M}\cD_M\otimes\dots\otimes_{\cO_M}\cD_M =
\sDiff_M(\cO_M,\dots,\cO_M;P) 
\] 
for differential operators in more
than two variables.

\

\medskip

\punkt {\bf Remark.} Lemma~\ref{lem:bidiffops_via_ind} provides an
interpretation of polydifferential operators using
$\cD$-modules. However, a more canonical interpretation can be given
using $\star$-operations on $\cD$-modules defined in
\cite{chiralalgebras}, see
Remark~\ref{rem:star_operations}
for a short summary.

\

\medskip

\punkt \ Notice that in our study of $\sDiff_M(Q;P)$, we never imposed
any restrictions, such as quasi-coherence, on the $\cO_M$-module $P$
(on the other hand, coherence of $Q$ is necessary for
Lemma~\ref{lem:bidiffops_via_ind}).  For this reason, the above
results extend immediately to the following situation: Let $p_X:Z\to
X$ be a morphism of complex manifolds, let $Q$ be an $\cO_X$-module
(preferably coherent), and let $P$ be an $\cO_Z$-module (or, more
generally, a $p_X^{-1}\cO_X$-module). Then we can consider
differential operators $p_X^{-1}Q\to P$ defined using the action of
$p_X^{-1}\cO_X$. Such differential operators form a sheaf on $Z$ that
we denote by $\sDiff_Z(p_X^{-1}Q;P)$. In particular, if $Q=\cO_X$, we
obtain the sheaf
\[
\sDiff_Z(p_X^{-1}\cO_X;P)=P\otimes_{p_X^{-1}\cO_X}p_X^{-1}\cD_X;
\] 
it is the right $p_X^{-1}\cD_X$-module induced by $P$.

Similarly, suppose we are given two morphisms $p_X:Z\to X$ and
$p_Y:Z\to Y$, and suppose $P_1$ (resp. $P_2$, $P$) is an
$\cO_X$-module (resp. an $\cO_Y$-module, an $\cO_Z$-module). We can
then consider bidifferential operators $p_X^{-1}P_1\times
p_Y^{-1}P_2\to P$; they form a sheaf on $Z$ that we denote
by $\sDiff_Z(p_X^{-1}P_1,p_Y^{-1}P_2;P)$. One can also consider
polydifferential operators in more than two variables in the same
setting. The details are left to the reader.

\

\medskip

\punkt \ All of the above remains valid in the relative setting.
Namely, let $p_X:Z\to X$ be a submersive morphism of complex
manifolds. Given $\cO_Z$-modules $P$, $Q$, we define the sheaf of
{\bfseries relative differential operators} from $P$ to $Q$ to be
\[
\sDiff_{Z/X}(P;Q):=\{A\in\sDiff_Z(P;Q):A(p_X^{-1}f)=(p_X^{-1}f)A 
\text{
  for all } f\in\cO_X\}.
\]
In particular, set $\cD_{Z/X}:=\sDiff_{Z/X}(\cO_Z;\cO_Z)$. 
For any $\cO_Z$-module $P$, the {\bfseries induced} 
right $\cD_{Z/X}$-module 
\[
P_{\cD/X}:=P\otimes_{\cO_Z}\cD_{Z/X}
\]
is identified with $\sDiff_{Z/X}(\cO_Z;P)$. 
Note that $P$ can be reconstructed from $P_{\cD/X}$ 
as 
\[
P=P_{\cD/X}\otimes_{\cD_{Z/X}}\cO_Z.
\]
The induction functor is adjoint to the forgetful functor from
$\cD_{Z/X}$-modules to $\cO_Z$-modules, and the relative versions of
Lemmas~\ref{lem:induced_adjunct} and \ref{lem:diffops_via_ind} hold. 
The proofs are parallel to the absolute case.

\

\medskip

\punkt {\bf Lemma.} \label{lem:rel_induced_adjunct}
{\it Let $P$ be an $\cO_Z$-module.
\begin{itemize}
\item[{\bf (a)}] For any right $\cD_{Z/X}$-module $F$, there are
  natural isomorphisms 
\begin{align*}
\Ext^i_{\cD_{Z/X}}(P_{\cD/X},F)=\Ext^i_{\cO_Z}(P,F)\qquad&(i\ge 0)\\
\sExt^i_{\cD_{Z/X}}(P_{\cD/X},F)=\sExt^i_{\cO_Z}(P,F)\qquad&(i\ge 0)
\end{align*}
of vector spaces and sheaves of vector spaces, respectively. 
\item[{\bf (b)}] For any left $\cD_{Z/X}$-module $F$, there is a 
natural isomorphism
\[
\sTor_i^{\cD_{Z/X}}(P_{\cD/X},F)=\sTor_i^{\cO_Z}(P,F)
\]
of sheaves of $\cO_Z$-modules.
In particular, for $F=\cO_Z$, there is an isomorphism
\[
P_{\cD/X}\otimes^\LL_{\cD_{Z/X}}\cO_Z=P.
\]
\end{itemize}
\ \hfill $\Box$
}

\

\medskip

\punkt {\bf Lemma.} \label{lem:rel_diffops_via_ind} {\it Let $Q$, $P$
be $\cO_Z$-modules.

\begin{itemize}
\item[{\bf (a)}] Any differential operator $A\in\sDiff_{Z/X}(Q;P)$
  induces a map  
\[
Q\to\sDiff_{Z/X}(\cO_Z;P) = P_{\cD/X} : \ s\mapsto A(\bullet s)
\qquad(s\in Q).
\]
This provides an identification $\sDiff_{Z/X}(Q;P) = 
\sHom_{\cO_Z}(Q,P_{\cD/X})$.
\item[{\bf (b)}] Any differential operator $A\in\sDiff_{Z/X}(Q;P)$ 
induces a map 
\[
\sDiff_{Z/X}(\cO_Z;Q)\to\sDiff_{Z/X}(\cO_Z;P) : \ B\mapsto A\circ B.
\] 
This provides an identification 
\[
\begin{split}
\sDiff_{Z/X}(Q;P) & =
\sHom_{\cD_{Z/X}}(\sDiff_{Z/X}(\cO_Z;Q),\sDiff_{Z/X}(\cO_Z;P)) \\
& = \sHom_{\cD_{Z/X}}(Q_{\cD/X},P_{\cD/X}).
\end{split}
\]
\end{itemize}
\ \hfill $\Box$
} 

\

\medskip

\subsection{$\cD$-modules on $\star$-quantizations}
\label{ssec:qinduced_D}

As before, let $M$ be a complex manifold. In this section, we consider differential
operators on $\star$-quantizations of $M$. 
It makes sense
to consider differential operators between $\cO$-modules on different
quantizations of $M$.

\

\medskip

\punkt {\bf Definition.}
\label{def:qudiffops} {\it Let $\qu{M}'$,
$\qu{M}$ be two neutralized $\star$-quantizations of $M$, and let
$\qu{P}$ be a $\cO_\qu{M}$-module. Let $A:\cO_{\qu{M}'}\to\qu{P}$
be an $R$-linear morphism. We say that $A$ is a {\bfseries differential
operator} if for every point $x\in M$, there exists a neighborhood $U$ 
of $x$ and $N\ge 0$ with the following property:

For any open subset $V\subset U$ and any choice of
sections
\[
\qu{f}'_k\in\Gamma(V,\cO_{\qu{M}'}),\quad
\qu{f}_k\in\Gamma(V,\cO_{\qu{M}})\quad\text{such
  that}\quad\qu{f}'_k=\qu{f}_k\mod\hbar\qquad(k=0,\dots,N),
\] 
the sequence of maps $A_k:\cO_{\qu{M}'}|_V\to\qu{P}|_V$ defined recursively by
\[
A_{-1}:=A|_V, \quad  A_k:= \qu{f}_kA_{k-1}-A_{k-1}\qu{f}'_k
\] 
satisfies $A_N=0$.

Let $\qu{P}'$ be an $\cO_{\qu{M}'}$-module.  An $R$-linear map 
$A:\qu{P}'\to\qu{P}$ is a {\bfseries differential operator} if for any local section
$\qu{s}\in\qu{P}'$, the map
\[\cO_{\qu{M}'}\to\qu{P}:\qu{f}\mapsto A(\qu{f}\qu{s})\]
is a differential operator. Differential operators form a sheaf that we denote by 
$\sDiff_M(\qu{P}';\qu{P})$.
}

\

\medskip 

\punkt \ Let us describe explicitly the sheaf
$\sDiff_M(\cO_{\qu{M}'};\qu{P})$. Recall that locally, the structure
sheaf of a neutralized $\star$-quantization of $M$ is isomorphic to
$\cO_M\otimes_\C R$ equipped with a $\star$-product. Choose such
isomorphisms
\[
\theta':\cO_M\otimes_\C
R\to\cO_{\qu{M}'}\qquad\text{and}\qquad\theta:\cO_M\otimes_\C
R\to\cO_\qu{M}
\]
for $\qu{M}$ and for $\qu{M}'$. For simplicity, we assume that
$\theta$ and $\theta'$ exist globally to avoid passing to an open
cover of $M$. Let us also suppose that $M$ is isomorphic to an open subset of $\C^n$;
let us fix global coordinates $(x_1,\dots,x_n)$ on $M$.

\

\medskip
  
\noindent
We claim that Definition~\ref{def:qudiffops} reduces
to the following description.
\

\smallskip

\punkt {\bf Lemma.} \label{lm:quinduced}
{\it A map $A:\cO_{\qu{M}'}\to\qu{P}$ is a differential
operator if and only if it can be written in the form
\begin{myequation}{eq:diffops}
\qu{g}'\mapsto\sum_{\alpha=(\alpha_1,\dots,\alpha_n)} \theta\left(
\frac{\partial^{\alpha_1}\cdots\partial^{\alpha_n}((\theta')^{-1}\qu{g}')}{\partial
x_1^{\alpha_1}\cdots\partial x_n^{\alpha_n}}\right) a_\alpha,
\end{myequation}
for sections $a_\alpha\in\qu{P}$ such that locally only finitely many
of $a_\alpha$'s are non-zero.}

\noindent
{\bf Proof.}  Denote by $\sDiff'_M(\cO_{\qu{M}'};\qu{P})$ the sheaf of
operators $A:\cO_{\qu{M}'}\to\qu{P}$ of the form
\eqref{eq:diffops}. Let us prove that $\sDiff'_M(\cO_{\qu{M}'};\qu{P})
\subset \sDiff_M(\cO_{\qu{M}'};\qu{P})$.

Fix $A\in\sDiff'_M(\cO_{\qu{M}'};\qu{P})$,  
and consider its order given by 
\[
\ord(A)=\max\{\alpha_1+\dots+\alpha_n:a_\alpha\ne 0\}.
\]
Passing to an open cover of $M$, we may assume that
$\ord(A)$ is finite. Let us choose 
a sequence of sections
$\qu{f}'_k\in\cO_{\qu{M}'},\qu{f}_k\in\cO_{\qu{M}}$ ($k\ge 0$) 
such that $\qu{f}'_k=\qu{f}_k\mod\hbar$, and define $A_k$ recursively by
\[
A_{-1}:=A|_V, \quad  A_k:= \qu{f}_kA_{k-1}-A_{k-1}\qu{f}'_k,
\] 
as in Definition~\ref{def:qudiffops}. We now prove that $A_N=0$ for
some $N$ that depends only on $\ord(A)$, the
orders of bidifferential operators giving
the $\star$-products on $\cO_{\qu{M}}$ and $\cO_{\qu{M}'}$, and the
number $r$ such that $\hbar^r\qu{P}=0$.  The proof proceeds by
induction on $r$. For $r=0$, we have $\qu{P}=0$ and we can take $N=0$
(or $N=-1$).

Suppose now $r>0$. Let $\rho(A)$ be the `reduction of $A$ modulo $\hbar$' defined to be the 
composition 
\[
\cO_{\qu{M}'}\to\qu{P}\to\qu{P}/\hbar\qu{P}.
\]
Note that $\qu{P}/\hbar\qu{P}$ is a $\cO_M$-module, and that $\rho(A)$
factors through $\cO_M$. In other words,
$\rho(A)\in\sDiff_M(\cO_M;\qu{P}/\hbar\qu{P})$. It is now easy to see
that
\[
{\ord}(\rho(A))> {\ord}(\rho(A_0))>
{\ord}(\rho(A_1))> \dots,
\]
and therefore $\rho(A_{{\ord}(A)})=0$. In other
words, the image of $A_{{\ord}(A)}$ is contained in
$\hbar P$.

Now notice that $A_0\in\sDiff'_M(\cO_{\qu{M}'};\qu{P})$, and moreover,
${\ord}(A_0)\le{\ord}(A)+C$ for
some constant $C$. We now see that
$A_{{\ord}(A)}$ is a section of
{$\sDiff'_M(\cO_{\qu{M}'};\hbar \qu{P})$} and
that its {order} is bounded by
{$\ord(A)(C+1) + C$}. {By the
  induction hypothesis $A_{\ord(A)}\in \sDiff'_M(\cO_{\qu{M}'};\hbar
  \qu{P})$ so $A \in \sDiff_M(\cO_{\qu{M}};\hbar \qu{P})$.}

To complete the proof of the lemma, we need to check that the
embedding
$\sDiff'_M(\cO_{\qu{M}'};\qu{P})\hookrightarrow\sDiff_M(\cO_{\qu{M}};\qu{P})$
is bijective. Again, we proceed by induction on $r$ such that
$\hbar^r\qu{P}=0$. The statement is clear if $r=1$. For $r>1$, it
follows from the exact sequences
\[
\xymatrix@-1pc{
  0\ar[r]&\sDiff'_M(\cO_{\qu{M}'};\hbar\qu{P})\ar[r]\ar[d]&
  \sDiff'_M(\cO_{\qu{M}'};\qu{P})\ar[r]\ar[d]&
  \sDiff'_M(\cO_{\qu{M}'};\qu{P}/\hbar\qu{P})\ar@{=}[d]\\ 
0\ar[r]&\sDiff_M(\cO_{\qu{M}'};\hbar\qu{P})\ar[r]&
\sDiff_M(\cO_{\qu{M}'};\qu{P})\ar[r]& 
  \sDiff_M(\cO_{\qu{M}'};\qu{P}/\hbar\qu{P}).}
\]
\ \hfill $\Box$

\

\medskip

\punkt {\bfseries Remark.} From the proof of Lemma~\ref{lm:quinduced}, we see that the 
notion of order is not very useful for differential operators on $\star$-quantizations.
The definition of $\ord(A)$ used in the proof is quite artificial and depends on the choice
of trivializations $\theta$ and $\theta'$. Although this `order' is related to the value
of $N$ appearing in Definition~\ref{def:qudiffops} (which is another candidate for 
order of differential operator), the relation is very indirect.

To give another example, consider differential operators $\cO_{\qu{M}}\to\cO_{\qu{M}}$.
Besides viewing $\cO_{\qu{M}}$ as a left module over itself, we can view it as a right module
over itself (or, equivalently, as an $\cO_{\qu{M}}^{op}$-module). It follows from
Lemma~\ref{lm:quinduced} that we obtain the same sheaf of differential operators in these
two cases. However, it seems impossible to define orders in compatible way: for instance,
right multiplication by a function is $\cO_\qu{M}$-linear with respect to the left 
$\cO_{\qu{M}}$-module structure (`has order zero') and non-linear with respect to the right
$\cO_{\qu{M}}$-module structure (`has positive order').

\

\medskip

\punkt \ \label{sssec:coaction} Consider the sheaf
$\cD_{\qu{M}'}^\qu{M}:=\sDiff_M(\cO_{\qu{M}'};\cO_\qu{M})$. It follows
from Lemma~\ref{lm:quinduced} that $\cD_{\qu{M}'}^\qu{M}$ is a flat
$R$-module equipped with an identification
$\cD_{\qu{M}'}^\qu{M}/\hbar\cD_{\qu{M}'}^\qu{M}=\cD_M$. Moreover,
$\cD_{\qu{M}'}^\qu{M}$ has a right action of $\cO_\qu{M}$ coming from
its action on itself that turns it into a locally free right
$\cO_\qu{M}$-module of infinite rank. (We ignore the commuting left
action of $\cO_\qu{M}$ and the two actions of $\cO_\qu{M'}$.) We now
get a natural identification
\[
\sDiff_M(\cO_{\qu{M}'};\qu{P}) =
\cD_{\qu{M}'}^\qu{M}\otimes_{\cO_\qu{M}}\qu{P}.
\]
Moreover, this identification respects the structure of a right module over
\[
\cD_{\qu{M'}} := \cD_{\qu{M}'}^{\qu{M}'} =
\sDiff_M(\cO_{\qu{M}'};\cO_{\qu{M}'}). 
\]
It is easy to see that the functor
\[
\qu{P}\mapsto\sDiff_M(\cO_{\qu{M}'};\qu{P})
\]
has properties similar to those of the induction $P\mapsto P_\cD$
studied in Section~\ref{ssec:induced_D}. Clearly, the functor is
exact, because $\cD_{\qu{M}'}^\qu{M}$ is flat over
$\cO_\qu{M}$. Moreover, given an $\cO_{\qu{M}_1}$-module $\qu{P}_1$
and a $\cO_{\qu{M}_2}$-module $\qu{P}_2$ on two neutralized
$\star$-quantizations $\qu{M}_1$, $\qu{M}_2$ of $M$, there is a
natural identification
\begin{myequation}{eq:qudiffops_via_ind}
\sDiff_M(\qu{P}_1;\qu{P}_2)
= \sHom_{\cD_{\qu{M}'}}(\sDiff_M(\cO_{\qu{M}'};\qu{P}_1),
\sDiff_M(\cO_{\qu{M}'};\qu{P}_2))  
\end{myequation}
defined as in Lemma~\ref{lem:diffops_via_ind}(b).

\

\medskip

\punkt {\bf Remark.} It may seem strange that we assign to an
$\cO_\qu{M}$-module $\qu{P}$ on a quantization $\qu{M}$ a $\cD$-module
on a different quantization $\qu{M}'$. Here is a different approach to
this assignment that should clarify the situation.

Let $\qu{P}$ be an $\cO_{\qu{M}}$-module on a neutralized
$\star$-quantization $\qu{M}$ of $M$. It is most natural to assign to
$\qu{P}$ the right $\cD_\qu{M}$-module $\sDiff_M(\cO_\qu{M};\qu{P})$
(corresponding to the choice $\qu{M}'=\qu{M}$). 

On the other hand, we {have the following} 

\

\noindent
{\punkt {\bf Proposition.}\label{prop:isomonodromic}}
\ {For any two neutralized $\star$-quantizations
  $\qu{M}$ and $\qu{M}'$ of $M$, the categories of right modules over
  $\cD_\qu{M}$ and over $\cD_{\qu{M}'}$ are canonically equivalent
  (the same applies to the categories of left $\cD$-modules, but this
  is irrelevant for our purposes).}

\

\noindent
{{\bfseries Proof.}}\  Indeed, the sheaf 
  $\cD_{\qu{M}'}^\qu{M}=\sDiff_M(\cO_{\qu{M}'};\cO_{\qu{M}})$ is
  naturally a $\cD_\qu{M}$-$\cD_{\qu{M}'}$-bimodule, and
  Lemma~\ref{lm:quinduced} shows that it is locally free of rank one
  as a left $\cD_\qu{M}$-module and as a right
  $\cD_{\qu{M}'}$-module. Therefore, it provides a Morita equivalence
\[
F\mapsto F\otimes_{\cD_\qu{M}}\cD_{\qu{M}'}^\qu{M}\qquad(F\text{ is a
  right $\cD_\qu{M}$-module})
\]
between the two categories of $\cD$-modules. The inverse equivalence
is provided by $\cD_\qu{M}^{\qu{M}'}$.  \ \hfill
{$\Box$}

\

\noindent
Note that the equivalence can
be viewed as a version of `isomonodromic transformation': since
$\qu{M}$ and $\qu{M}'$ are quantizations of the same manifold, they
differ only `infinitesimally', and it is well known that the notion of
a $\cD$-module is `topological' in the sense that it is invariant
under infinitesimal deformations.

We can now use the equivalence to pass between $\cD$-modules on
different neutralized $\star$-quantizations of $M$. It is easy to see
that
\[
\sDiff_M(\cO_\qu{M};{\qu{P}})
\otimes_{\cD_\qu{M}}\cD_{\qu{M}'}^\qu{M} =
\sDiff_M(\cO_{\qu{M}'};{\qu{P}}).
\]
This allows us to interpret the
{identification \eqref{eq:qudiffops_via_ind}}
as follows: the sheaves $\qu{P}_i$ ($i=1,2$) give rise to induced
$\cD_{\qu{M}_i}$-modules. Using the `isomonodromic transformation', we
can put this $\cD$-modules into a single category (of
$\cD_{\qu{M}'}$-modules), so that it makes sense to talk about
homomorphisms between them. One can thus say that the `isomonodromic
transformation' is the reason why differential operators between
$\cO$-modules on different quantizations make sense.

It should be clear from this discussion that, when one works with $\cD$-modules
on a neutralized $\star$-quantization $\qu{M}'$ of $M$, the choice of 
$\qu{M}'$ is largely irrelevant. For this reason, in most
situations we can (and will) take $\qu{M}'$ to be the trivial
deformation of $M$: $\cO_{\qu{M}'}=\cO_M\otimes_\C R$. The advantage
of this choice is that $\cD_{\qu{M}'}$-modules then become simply
modules over $\cD_M\otimes_\C R$. On the other hand, in the relative
situation (discussed below) it is no
longer true that the category of relative $\cD$-modules
on a quantization is independent of the choice of quantization, see Remark~\ref{rem:relative_isomonodromy}. However,
relative $\cD$-modules on quantizations appear only in
Section~\ref{ssec:obstructionsneutr'}, which is independent from the
rest of the paper.

\

\medskip

\punkt {\bf Definition.}
{\it Let  $\qu{M}_1,\dots,\qu{M}_n,\qu{M}$ be neutralized
$\star$-quantizations of $M$ (for
some $n\ge 1$), and let
$\qu{P}$ be a $\cO_\qu{M}$-module. Let $A:\cO_{\qu{M}_1}\times\dots\times\cO_{\qu{M}_n}\to\qu{P}$
be an $R$-polylinear morphism. We say that $A$ is a {\bfseries polydifferential
operator} if locally on $M$, there exists $N\ge 0$ with the following property:

For any $j=1,\dots,n$ and any choice of (local) sections \[
\qu{f}'_k\in\cO_{\qu{M}_j},\quad
\qu{f}_k\in\cO_{\qu{M}}\quad\text{such
  that}\quad\qu{f}'_k=\qu{f}_k\mod\hbar\qquad(k=0,\dots,N),
\] 
the sequence of maps $A_k:\cO_{\qu{M}_1}\times\dots\times\cO_{\qu{M}_n}\to\qu{P}$ defined recursively by
\[
A_{-1}:=A, \quad  A_k:=\qu{f}_kA_{k-1}-A_{k-1}(id\times\dots\times\qu{f}'_k\times\dots\times id)
\] 
satisfies $A_N=0$. Here $(id\times\dots\times\qu{f}'_k\times\dots\times id)$ (with $\qu{f}'_k$ in the
$j$-th position) refers to the map from $\cO_{\qu{M}_1}\times\dots\times\cO_{\qu{M}_N}$ to itself that
multiplies $j$-th argument by $\qu{f}'_k$ and leaves the remaining $n-1$ arguments unchanged.

Let $\qu{P}_1,\dots,\qu{P}_n$ be $\cO$-modules on $\qu{M}_1,\dots,\qu{M}_n,\qu{M}$. An 
$R$-polylinear morphism
\[A:\qu{P}_1\times\dots\qu{P}_n\to\qu{P}\]
is a {\bfseries polydifferential operator} if for any choice of local sections $\qu{s}_1\in\qu{P}_1,\dots,\qu{s}_n\in\qu{P}_n$,
the map
\[\cO_{\qu{M}_1}\times\dots\times\cO_{\qu{M}_n}\to\qu{P}:(\qu{g}_1,\dots,\qu{g}_n)\mapsto A(\qu{g}_1\qu{s}_1,\dots,\qu{g}_n\qu{s}_n)\]
is a polydifferential operator. Polydifferential operators form a sheaf, which we
denote by $\sDiff_M(\qu{P}_1,\dots,\qu{P}_n;\qu{P})$.  }

\

\medskip

\punkt \ Let $p_X:Z\to X$ be a morphism of complex
manifolds. Suppose $\qu{Z}$ is a neutralized $\star$-quantization of
$Z$ and $\qu{X}$ is a neutralized $\star$-quantization of $X$;
$\qu{Z}$ and $\qu{X}$ do not need to be compatible in any way.  Let
$\qu{Q}$ be an $\cO_\qu{X}$-module and let $\qu{P}$ be an
$\cO_\qu{Z}$-module.  It is easy to extend
Definition~\ref{def:qudiffops} to construct the sheaf of differential
operators $p_X^{-1}\qu{Q}\to\qu{P}$, which we denote by
$\sDiff_Z(p_X^{-1}\qu{Q};\qu{P})$.

Explicitly, an $R$-linear map $A:p_X^{-1}\qu{Q}\to\qu{P}$ is a differential operator if
for every point $z\in Z$ and every section $s\in p_X^{-1}\qu{Q}$
defined at $z$, there exists a neighborhood $U$ of $z$ and $N\ge 0$
with the following property:

For any open subset $V\subset U$ and for any choice of sections
\[
\qu{f}'_k\in\Gamma(V,p_X^{-1}\cO_\qu{X}),\quad
\qu{f}_k\in\Gamma(V,\cO_\qu{Z})\quad\text{such
  that}\quad\qu{f}'_k=\qu{f}_k\mod\hbar\qquad(k=0,\dots,N),
\] 
the sequence of maps $A_k:p_X^{-1}\qu{P}|_V\to\qu{Q}|_V$ defined by
\[
A_{-1}:=A|_V, \quad  A_k:= \qu{f}_kA_{k-1}-A_{k-1}\qu{f}'_k
\] 
satisfies $A_N(s|_V)=0$.

\

\medskip

\punkt \ Finally, let us consider relative differential operators on
quantizations. Let \linebreak $p_X:Z\to X$ be a submersive morphism of
complex manifolds, and let $\qu{Z}$ and $\qu{X}$ be neutralized
$\star$-quantizations of $Z$ and $X$, respectively. Suppose that
$p_X:Z\to X$ extends to a morphism of quantizations
$\qu{p}_X:\qu{Z}\to\qu{X}$. In other words, $\qu{p}_X$ is an
$R$-linear morphism of ringed spaces that acts on underlying sets as
$p_X$ and such that the reduction of
$\qu{p}_X^{-1}:p_X^{-1}\cO_\qu{X}\to\cO_\qu{Z}$ modulo $\hbar$ equals
$p_X^{-1}$.  We are mostly interested in the case when
$\qu{Z}=\qu{X}\times\qu{Y}$ for a neutralized $\star$-quantization
$Y$, and $\qu{p}_X$ is the natural projection.

Now suppose $\qu{Z}'$ is another neutralized $\star$-quantization of
$\qu{Z}$ and suppose that $p_X:Z\to X$ also extends to a morphism
$\qu{p}'_X:\qu{Z}'\to\qu{X}$. Now given $\cO$-modules $\qu{P}'$,
$\qu{P}$ on $\qu{Z}'$ and $\qu{Z}$, respectively, we define the sheaf
of {\bfseries relative differential operators} by essentially the same
formula as before:
\[\sDiff_{Z/X}(\qu{P}';\qu{P}):=\{A\in\sDiff_Z(\qu{P}';\qu{P}):
A((\qu{p}'_X)^{-1}\qu{f})=(\qu{p}_X^{-1}\qu{f})A\text{ for all }
\qu{f}\in\cO_\qu{X}\}.\]

\

\medskip

\noindent
The above properties of differential operators remain true in 
the setting of relative quantizations.
\

\smallskip

\punkt {\bf Proposition.}\label{prop:relqudiffops_flat} {\it The sheaf
$\sDiff_{Z/X}(\cO_{\qu{Z}'};\cO_\qu{Z})$ is flat over $R$, and the
obvious `reduction modulo $\hbar$' map
\[
\sDiff_{Z/X}(\cO_{\qu{Z}'};\cO_\qu{Z})/
\hbar\sDiff_{Z/X}(\cO_{\qu{Z}'};\cO_\qu{Z})\to\cD_{Z/X}
\]
is an isomorphism.}

\noindent
{\bf Proof.} Looking at the filtration by powers of $\hbar$, we see 
that it suffices to prove that the `reduction modulo $\hbar$' map is
surjective. Let us fix $A\in\cD_{Z/X}$, and let us lift it to an
operator $\qu{A}\in\sDiff_{Z}(\cO_{\qu{Z}'};\cO_\qu{Z})$ (this can
always be done locally). It then leads to a bidifferential operator
\[
\qu{B}:p_X^{-1}\cO_\qu{X}\times\cO_{\qu{Z}'}\to\cO_\qu{Z}:
(\qu{f},\qu{g})\mapsto\qu{f}\qu{A}(\qu{g})-\qu{A}(\qu{f}\qu{g})
\]
that vanishes modulo $\hbar$. Note that $\qu{B}$ is a cocycle:
\[
\qu{f}_1\qu{B}(\qu{f}_2,\qu{g}) -
\qu{B}(\qu{f}_1\qu{f}_2,\qu{g})+\qu{B}(\qu{f}_1,\qu{f}_2\qu{g}) 
= 
0 \qquad 
(\qu{f}_1,\qu{f}_2\in p_X^{-1}\cO_\qu{X};\qu{g}\in\cO_{\qu{Z}'}).
\] 
We need to show that it is locally a coboundary: 
there exists a differential operator 
\[
\qu{C}:\cO_{\qu{Z}'}\to\cO_\qu{Z}
\]
that vanishes modulo $\hbar$ and satisfies
\[
\qu{B}(\qu{f},\qu{g}) = \qu{f}\qu{C}(\qu{g}) -
\qu{C}(\qu{f}\qu{g})\qquad(\qu{f}\in 
p_X^{-1}\cO_\qu{X};\qu{g}\in\cO_{\qu{Z}'}).
\]
Then $\qu{A}-\qu{C}$ would lift $A$ to a relative differential
operator in $\sDiff_{Z/X}(\cO_{\qu{Z}'};\cO_\qu{Z})$.

Arguing by induction in powers of $\hbar$, we see that it suffices to prove
the following statement:
If $B\in\sDiff_Z(p_X^{-1}\cO_X,\cO_Z;\cO_Z)$ satisfies the cocycle condition
\[
f_1 B(f_2,g) - B(f_1f_2,g) + B(f_1,f_2g) = 
0\qquad(f_1,f_2\in p_X^{-1}\cO_X;g\in
\cO_Z),
\] 
then it is locally a coboundary: locally, there exists $C\in\cD_Z$ such that
\[
B(f,g) = f C(g)-C(fg) \qquad (f\in p_X^{-1}\cO_{X};g\in\cO_Z).
\]
This can be proved directly. 

Indeed, let $(x_1,\dots,x_n)$ be local coordinates on $X$; we denote
their pullbacks to $Z$ in the same way. Then changing $B$ by a
coboundary, we can ensure that $B(x_1,g)=0$. Since $B$ is a cocycle,
this implies that
\[
B(f,x_1g) = B(x_1f,g) = x_1B(f,g) \qquad (f\in
p_X^{-1}\cO_X,g\in\cO_Z). 
\]
In other words, $B$ does not include differentiation with respect to
$x_1$: it is a `bidifferential operator relative to $x_1$'. We can now
treat $x_2$ in the same way, and so on. \ \hfill $\Box$

\

\medskip

\punkt {\bfseries Remark.} Essentially, the proof of
Proposition~\ref{prop:relqudiffops_flat} verifies vanishing of the
first local Hochschild cohomology of $\cO_X$ with coefficients in
$\cD_Z$. Here `local' means that the Hochschild complex is formed
using cochains that are differential operators. (Confusingly, the
complex is also local in a different sense: its terms and cohomology
objects are sheaves rather than vector spaces.)  The proof can be made
more transparent using the framework of comodules introduced in
Section~\ref{ssec:coalgebras}, see Remark~\ref{rem:relqudiffops_flat}.

\

\medskip

\punkt {\bf Corollary.}\label{cor:relqudiffops_ind} {\it Consider on
$\sDiff_{Z/X}(\cO_{\qu{Z}'};\cO_\qu{Z})$ the right action of
$\cO_\qu{Z}$ coming from its right action on itself. Then the map
\[
\sDiff_{Z/X}(\cO_{\qu{Z}'};\cO_\qu{Z})
\otimes_{\cO_\qu{Z}}\qu{P} \to \sDiff_{Z/X}(\cO_{\qu{Z}'};\qu{P}):A\otimes 
p\mapsto A(\bullet)p 
\]
is an isomorphism for any $\cO_\qu{Z}$-module $\qu{P}$.}

\noindent
{\bf Proof.} Consider the filtration by powers of $\hbar$. By
Proposition~\ref{prop:relqudiffops_flat} and
Lemma~\ref{lem:rel_diffops_via_ind}, we see that the map induces an
isomorphism on the associated graded quotients. This implies the
corollary.  \ \hfill $\Box$

\

\medskip

\punkt 
We now obtain a functor
\[
\qu{P}\mapsto\sDiff_{Z/X}(\cO_{\qu{Z}'};\qu{P})
\]
from the category of $\cO_{\qu{Z}}$-modules to the category of right
modules over
\[
\cD_{\qu{Z'}/\qu{X}}:=\sDiff_{Z/X}(\cO_{\qu{Z}'};\cO_{\qu{Z}'}).
\]
By Corollary~\ref{cor:relqudiffops_ind}, this functor can be rewritten as
\[
\qu{P}\mapsto\qu{P}\otimes \sDiff_{Z/X}(\cO_{\qu{Z}'};\cO_\qu{Z})
\otimes_{\cO_\qu{Z}}\qu{P}.\]
Note that $\cD_{Z/X}$ is a flat right $\cO_{Z}$-module, 
and therefore Proposition~\ref{prop:relqudiffops_flat} implies that
$\sDiff_{Z/X}(\cO_{\qu{Z}'};\cO_\qu{Z})$ is a flat right
$\cO_\qu{Z}$-module. 
Hence the functor is
exact. 

\

\medskip

\punkt {\bfseries Remark.} \label{rem:relative_isomonodromy} In the relative
setting, we have a version of the identification \eqref{eq:qudiffops_via_ind},
and of Propositon~\ref{prop:isomonodromic}. However, there is an
important difference between relative and absolute differential operators on quantizations.
Namely, we can consider absolute differential operators between $\cO$-modules on any two neutralized $\star$-quantizations
of a manifold. On the other hand, given a submersive map $Z\to X$ and $\cO$-modules on two neutralized $\star$-quantizations $\qu{Z}$ and
$\qu{Z}'$ of $Z$,
it is necessary to have a neutralized $\star$-quantization $\qu{X}$ of $X$ and quantizations of $p_X$ to morphisms
$\qu{p}_X:\qu{Z}\to\qu{X}$, $\qu{p}'_X:\qu{Z}'\to\qu{X}$ in order for relative differential operators to make sense. 

Thus, a relative version of Proposition~\ref{prop:isomonodromic} provides an equivalence between the categories
of (right) modules over $\cD_{\qu{Z}/\qu{X}}$ and over $\cD_{\qu{Z}'/\qu{X}}$. In this way, the category of
$\cD_{\qu{Z}/\qu{X}}$ does not actually depend on the choice of $\qu{Z}$; however, the category does depend on the quantization 
$\qu{X}$ of $X$. (See \ref{ssec:relativecomodules} for another explanation.)
As an extreme example, suppose $Z=X$ and $\qu{Z}=\qu{X}$. Then $\cD_{\qu{Z}/\qu{X}}=\cO_{\qu{X}}$,
and the category of $\cO_\qu{X}$-modules clearly depends on the choice of quantization $\qu{X}$.

\subsection{Coalgebras in the category of $\cD$-modules}
\label{ssec:coalgebras}

\punkt Let $M$ be a complex manifold. The category of left
$\cD_M$-modules $\Mod{\cD_M}$ is naturally a tensor category: given
$F,G\in\Mod{\cD_M}$, the tensor product $F\otimes_{\cO_M} G$ carries a
left action of $\cD_M$ (the action of vector fields is given by the Leibnitz
rule from the actions on $F$ and on $G$).
Consider coalgebras in this tensor category.

For the rest of this section, let us fix a coalgebra
$\cA\in\Mod{\cD_M}$. Here and in the rest of the paper, all coalgebras
are assumed to be coassociative and counital, but not necessarily
cocommutative.  We also assume that $\cA$ is flat over $\cO_M$. It is
well known that this assumption is necessary to define kernel of
morphisms between comodules.

 Consider now the category of right $\cD_M$-modules
$\Mod{\cD_M^{op}}$. It has an action of the tensor category
$\Mod{\cD_M}$: given $F\in\Mod{\cD_M}$ and $G\in\Mod{\cD_M^{op}}$, the
tensor product $F\otimes_{\cO_M} G$ carries a right action of
$\cD_M$. It makes sense to talk about $\cA$-comodules in the category
of right $\cD_M$-modules.

\

\medskip

\punkt {\bfseries Definition.}\label{def:comodules} {\it Let
$\Comod{\cA}$ be the category of (left) $\cA$-comodules in the
category $\Mod{\cD_M^{op}}$.  We call objects of $\Comod{\cA}$ simply
`$\cA$-comodules'.}

\

\medskip

\punkt {\bfseries Lemma.}\label{lm:comod_inj} {\it The category
  $\Comod{\cA}$ is an abelian category with enough injectives.}

\noindent
{\bf Proof.} Let $f : F\to G$ be a morphism of
  $\cA$-comodules. The fact that $F$ is $\cO_{M}$-flat implies that
  the kernel of $f$ in the category of right $\cD_{M}$-modules is
  automatically an $\cA$-comodule. Thus $\Comod{\cA}$ is abelian.
The forgetful functor $\Comod\cA\to\Mod{\cD_M^{op}}$ admits a right
adjoint: the coinduction functor
\[
\Mod{\cD_M^{op}}\to\Comod\cA \ :  \ F\mapsto\cA\otimes_{\cO_{M}} F.
\]
Being the right adjoint functor of an exact functor, it preserves injectivity. 
Thus we can obtain sufficiently many injective
$\cA$-comodules by coinduction.  \ \hfill $\Box$

\

\medskip

\punkt \ \label{sssec:defo.comodules} Let us now study deformation of
$\cA$-comodules. Let $\Mod{\cD_M\otimes_\C R}$ be the category of left
$\cD_M\otimes_\C R$-modules. Alternatively, $\Mod{\cD_M\otimes_\C R}$
is the category of $R$-modules in the $\C$-linear category
$\Mod{\cD_M}$. The tensor product over $\cO_M\otimes_\C R$ turns
$\Mod{\cD_M\otimes_\C R}$ into a tensor category.

Let $\qu\cA$ be an $R$-deformation of $\cA$. That is, $\qu\cA$ is a
coalgebra in $\Mod{\cD_M\otimes_\C R}$ that is flat over $R$ and
equipped with an isomorphism $\cA\simeq\qu\cA/\hbar\qu\cA$.

We extend the conventions of Definition~\ref{def:comodules} to
$\qu\cA$-comodules. Thus $\qu\cA$-comodules are by definition left
$\qu\cA$-comodules in the category $\Mod{\cD_M^{op}\otimes_\C R}$ of
right $\cD_M\otimes_\C R$-modules; the category of
$\tilde\cA$-comodules is denoted by $\Comod{\qu\cA}$.

Given $Q\in\Comod{\cA}$, we consider the problem of deforming $Q$ to
$\qu{Q}\in\Comod{\qu\cA}$. More precisely, let $\Defst(Q)$ be the
groupoid of comodules $\qu{Q}\in\Comod{\qu\cA}$ that are flat over $R$
together with an isomorphism $Q\simeq\qu{Q}/\hbar\qu{Q}$. As one would
expect, deformations of $Q$ are controlled by $\Ext$'s from $Q$ to
itself.

\

\medskip

\punkt {\bfseries Proposition.}\label{prop:defcomodules} {\it

\begin{itemize} 
\item[\bf(a)] If $\Ext^2_\cA(Q,Q)=0$, the groupoid $\Defst(Q)$ is
non-empty.

\item[\bf(b)] If $\Ext^1_\cA(Q,Q)=\Ext^2_\cA(Q,Q)=0$, the groupoid
$\Defst(Q)$ is connected (that is, all objects are isomorphic).

\item[\bf(c)] If $\Ext^1_\cA(Q,Q)=0$ and $\qu{Q}\in\Defst(Q)$, then
the algebra $\End_{\qu\cA}(\qu{Q})$ is an $R$-flat deformation of
$\End_\cA(Q)$.
\end{itemize} 
} 

\noindent
{\bf Proof.}  The result is essentially classical: for modules over an
algebra, it is sketched in \cite[\S~2]{La_massey} (and mentioned in
\cite[p.~150]{La_moduli}); and the case we are interested in
is similar. Actually, as stated in \cite[p.~150]{La_moduli}, the
deformation theory applies `in all cases where we have a good
cohomology and an obstruction calculus'.

More rigorously, the proposition follows from general theory of
deformation of objects in abelian categories due to Lowen
(\cite{Lowen}).  First note that $\Comod{\qu{\cA}}$ is a flat
deformation of $\Comod{\cA}$ in the sense of
\cite[Proposition~3.4]{Lowen_VanDerBergh}. Since
$\Comod{\qu\cA}$ has enough injectives, we need to check that
injective $\qu\cA$-comodules are
$R$-coflat. It suffices to check this for coinduced injective objects,
that is, for objects of the form $\qu{\cA}\otimes_{\cO_{M}} I$, where $I$ is an injective
  $\cD^{op}_{M}$-module. We need to verify that
  $\Ext^{1}_{R}(-,\qu{\cA}\otimes_{\cO_{M}} I) = 0$. But
  $\Ext^{1}_{R}(-,\qu{\cA}\otimes_{\cO_{M}} I)$ is equal to
  $\qu{\cA}\otimes_{\cO_{M}} \Ext^{1}_{R}(-,I)$ by
  \cite[Proposition~2.10]{Lowen_VanDerBergh} and the desired vanishing
  follows from properties (1) and (2) in \cite[page
    5441]{Lowen_VanDerBergh}.  Now {\bf (a)} and {\bf (b)} follow by
iterated application of the dual version of
\cite[Theorem~A]{Lowen}. To prove {\bf (c)}, we notice that the
filtration of $\qu{Q}$ by powers of $\hbar$ induces a filtration of
$\End_{\qu\cA}(\qu{Q})$, and the associated graded for
 this filtration is $\End_\cA(Q)\otimes_\C R$.  \ \hfill
$\Box$ 

\

\medskip

\punkt {\bfseries Remark.} We consider flat deformations, which could
be viewed as lifts along the functor
\[
\Comod{\qu\cA}\to\Comod\cA:\qu{Q}\mapsto\qu{Q}\otimes_R\C.
\]
On the other hand, \cite[Theorem~A]{Lowen} concerns
  lifting coflat objects, rather than flat ones. That is,
  \cite[Theorem~A]{Lowen}  gives lifts along the functor
\[
\Comod{\qu\cA}\to\Comod\cA:\qu{Q}\mapsto\sHom_R(\C,\qu{Q});
\]
  As long as one considers deformation over
$R=\C[\hbar]/\hbar^{n+1}$ (rather than over an arbitrary Artinian
ring) the two kinds of deformations coincide, because $R$ is
a Gorenstein, zero
  dimensional ring. Indeed, if $P$ is coflat
  $\qu{P}\otimes_{R}\mathbb{C} \cong
  \sHom_{R}(R,\qu{P})\otimes_{R}\mathbb{C} \cong
  \sHom_{R}(\Hom_{R}(\mathbb{C},R),\qu{P})$, see
  e.g. \cite[Proposition~2.9(3)]{Lowen_VanDerBergh}. Since $R$ is a
  local Gorenstein zero dimensional ring, $\Hom_{R}(\mathbb{C},R) =
  R$, and the statement follows. The general case when $R$ is an
  arbitrary Artinian ring would follow from a version of \cite[Theorem~A]{Lowen}
  for flat objects; this `dual version' is left as an exercise in \cite{Lowen}. 

\subsection{Neutralized $\star$-quantizations as
  $\cD$-coalgebras}\label{ssec:quantcoalgebras} 

We now interpret neutralized $\star$-quantizations and $\cO$-modules
on them using the framework of coalgebras and comodules developed in
Section~\ref{ssec:coalgebras}.

\

\medskip

\punkt \ Let $M$ be a complex manifold. Given a neutralized
$\star$-quantization $\qu{M}$ of $M$, consider differential operators
between $\cO_{\qu{M}}$ and the `trivial' deformation $\cO_M\otimes_\C
R$. Set
\[
\cA_\qu{M}:=\sDiff_M(\cO_\qu{M};\cO_M\otimes_\C R).
\]
Composition with differential operators on $\cO_M\otimes_\C R$ equips
$\cA_\qu{M}$ with a structure of a left $\cD_M\otimes_\C R$-module.

\

\medskip

\punkt {\bf Example: $\star$-product.}  Suppose
$\cO_\qu{M}=\cO_M\otimes_\C R$ with multiplication given by some
$\star$-product $f\star g$. By Lemma~\ref{lm:quinduced}, a map
$A:\cO_{\qu{M}}\to\cO_M\otimes_\C R$ is a differential operator if and
only if it is of the form 
\[
A = A_0+ {\hbar A_{1} + } \dots +
{\hbar^{n}} A_n\in\sDiff_M(\cO_M;\cO_M)\otimes_\C R\quad
(A_i\in\sDiff_M(\cO_M;\cO_M)).
\]
(Technically, Lemma~\ref{lm:quinduced} applies only locally on $M$, but
once it is proved that $A$ has the required form locally, it is also true globally.) 
We now see that  $\cA_{\qu{M}}$ is the trivial $\cD_M\otimes_\C R$-module of rank one.

In general, a neutralized $\star$-quantization $\cO_\qu{M}$ is locally
isomorphic to $\cO_M\otimes_\C R$ equipped with a
$\star$-product. Therefore, $\cA_\qu{M}\in\Mod{\cD_M\otimes_\C R}$ is
a rank one locally free $\cD_M\otimes_\C R$-module.

\

\medskip

\punkt \ The $\cD_M$-module $\cA_\qu{M}$ has a natural
coproduct. Indeed, consider the tensor product
\[
\cA_\qu{M}\underset{\cO_M\otimes_\C R}\otimes\cA_\qu{M}
\] 
with the usual $\cD_M$-module structure. It can be identified with the
sheaf of bidifferential operators
$\sDiff(\cO_{\qu{M}},\cO_{\qu{M}};\cO_M\otimes_\C R)$
(cf. Example~\ref{exm:bidiffops_via_ind}).  The identification sends
$A_1\otimes A_2$ ($A_1,A_2\in\cA_\qu{M}$) to the operator
\[
(\qu{f}_1,\qu{f}_2)\mapsto
A_1(\qu{f}_1)A_2(\qu{f}_2)\qquad(\qu{f}_1,\qu{f}_2\in\cO_\qu{M}).
\]
Now define the coproduct
$\bDelta:\cA_{\qu{M}}\to\cA_{\qu{M}}\otimes\cA_{\qu{M}}$ by
$\bDelta(A)=A\circ\mult,$ where
\[
\mult:\cO_{\qu{M}}\times \cO_{\qu{M}}\to\cO_{\qu{M}}
\]
is the multiplication. 

\

\medskip

\punkt {\bf Proposition.}\label{prop:quantization_coalgebras} {\it The
  functor $\qu{M}\mapsto\cA_\qu{M}$ is an equivalence between the
  groupoid of neutralized 
   $\star$-quantizations of $M$ and the groupoid of
  $R$-flat deformations of $\cD_M$ as a
  coassociative and counital coalgebra in
  $\Mod{\cD_M}$.}

\noindent
{\bf Proof.} Let us provide an inverse functor. Let $\qu{\cA}$ be an
$R$-deformation of the left $\cD_M$-module $\cD_M$.  In particular,
$\qu{\cA}$ is a locally free $\cD_M\otimes_\C R$-module of rank
one. Set
\[
\cO_\qu{M} := \sHom_{\cD_M\otimes_\C R}(\qu{\cA},\cO_M\otimes_\C R).
\] 
Then $\cO_\qu{M}$ is $R$-flat and $\cO_\qu{M}/\hbar\cO_\qu{M}=\cO_M$.

Suppose now that $\qu\cA$ carries a coassociative coproduct
$\bDelta:\qu\cA\to\qu\cA\otimes\qu\cA$ extending the canonical
coproduct on $\cD_M$. Then $\bDelta$ induces an associative product on
$\cO_\qu{M}$ via
\[
\qu{f}\cdot\qu{g} := (\qu{f}\otimes\qu{g})\circ \bDelta\qquad(\qu{f},
\qu{g}\in\sHom_{\cD_M\otimes_\C R}(\qu{\cA},\cO_M\otimes_\C R)),
\]
and the identification $\cO_\qu{M}/\hbar\cO_\qu{M}=\cO_M$ respects
this product.

Finally, locally we can lift $1\in\cD_M$ to a section
$\qu{1}\in\qu\cA$. This yields a local
 isomorphism
\[
\cO_\qu{M}\to\cO_M\otimes_\C R:\qu{f}\mapsto\qu{f}(\qu{1})
\] 
of sheaves of $R$-modules.  Under this isomorphism, the product on
$\cO_\qu{M}$ corresponds to the $\star$-product on $\cO_M\otimes_\C R$
given by the bidifferential operator $B$ such that
\[ 
\bDelta(\qu{1})=B(\qu{1}\otimes\qu{1}).
\]
Therefore, $\qu{M}$ is a neutralized $\star$-quantization of $M$.
Verifying that this is indeed the inverse construction to
$\qu{M}\mapsto\cA_\qu{M}$ is straightforward.  \ \hfill $\Box$\

\

\medskip

\punkt {\bf Examples.} \label{exm:oppositeproduct} The coalgebra
$\cA_{\qu{M}^{op}}$ corresponding to the opposite quantization
$\qu{M}^{op}$ is the opposite coalgebra of $\cA_\qu{M}$: they coincide
as $\cD$-modules, but the coproducts are opposite to each other.

Let $\qu{M}$ and $\qu{N}$ be neutralized $\star$-quantizations of
complex manifolds $M$ and $N$, respectively. Consider the quantization
$\qu{M}\times\qu{N}$ of $M\times N$. The corresponding coalgebra
$\cA_{\qu{M}\times\qu{N}}$ can be obtained as the tensor product
$p_M^*\cA_\qu{M}\otimes p_N^*\cA_\qu{N}$. Here $p_M:M\times N\to M$
and $p_N:M\times N\to N$ are the projections, $p_M^*$ and $p_N^*$ are
pullback functors for left $\cD$-modules, and $\otimes$ is the tensor product over $\cO_{M\times
N}\otimes_\C R$ (which gives the tensor structure on
$\Mod{\cD_{M\times N}\otimes_\C R}$).

\

\medskip

\punkt {\bfseries Remark.} Let $p:Z\to X$ be a morphism of complex
manifolds.  In this paper, we denote the corresponding pullback
functor for left $\cD$-modules by $p^*$. Note that $p^*$ coincides
with the $\cO$-module pullback. In the $\cD$-module literature (see
e.g. \cite[{page 232}]{Borel-dmodules}), this
pullback functor is denoted by $p^!$ since under the Riemann-Hilbert
correspondence it corresponds to the right adjoint of the pushforward
with compact supports, and $p^*$ refers to a different functor (its
Verdier dual). In our setup, we always require that $p$ is a
submersive map, which ensures that the two functors coincide up to a
cohomological shift: $p^!=p^*[2(\dim(Z)-\dim(X))]$. For this reason, we
hope that our notation is not unnecessarily confusing.

\

\medskip

\punkt \label{ssec:coaction} Now let us study $\cO$-modules, beginning
with the non-quantized complex manifold $M$. Let $P$ be an
$\cO_M$-module. Recall that the induced right $\cD_M$-module $P_\cD$
is defined by
\[ 
P_\cD = P\otimes_{\cO_M}\cD_M = \sDiff_M(\cO_M;P).
\]
We claim that it carries a natural coaction 
\[ 
\bDelta_P:P_\cD\to\cD_M\otimes_{\cO_M} P_\cD
\] 
of the coalgebra $\cD_M$. (Recall that $\cO_M$ acts on $P_\cD =
{\sDiff_{M}(\mathcal{O}_{M};P)}$  
{through the action of $\cO_{M}$ on itself, so that $\cO_M$ plays the role of $Q_{1}$
  in Lemma~\ref{lem:bidiffops_via_ind}, and
  Example~\ref{exm:bidiffops_via_ind}).} Moreover, $\bDelta_P$ is a
morphism of right $\cD_M$-module. Let us describe it in concrete
terms.

The tensor product $\cD_M\otimes_{{\cO_{M}}} P_\cD$ is identified with the sheaf
of bidifferential operators $\sDiff(\cO_M,\cO_M;P)$, see
Lemma~\ref{lem:bidiffops_via_ind} and
Example~\ref{exm:bidiffops_via_ind}.  The identification sends
$A\otimes B$ for $A\in\cD_M$, $B\in P_\cD$ to the operator
\[
(f,g)\mapsto B(A(f)g)\qquad (f,g\in\cO_M).
\]

Being the tensor product of a left $\cD_M$-module and a right $\cD_M$-module, the sheaf
$\cD_M\otimes_{{\cO_{M}}} P_\cD$ carries a natural structure of a right
$\cD_M$-module.  Explicitly, it is given by the following formulas
\begin{align*}
(A\otimes B)\cdot f&=(fA)\otimes B=A\otimes(Bf),\\
(A\otimes B)\cdot \tau&=-(\tau A)\otimes B+A\otimes(B\tau)
\quad (A\in{\cD_{M}},B\in
{P_\cD},f\in\cO_M,\tau\in{\cal T}_M), 
\end{align*}
where ${\cal T}_M$ is the sheaf of vector fields on $M$.

\

\noindent
{\punkt {\bf Lemma.}}  
{\it Under the identification
\[\cD_M\otimes_{{\cO_{M}}} P_\cD=\sDiff_M(\cO_M,\cO_M;P),\] 
the right action of $\cD_M$ on $\cD_M\otimes_{{\cO_{M}}} P_\cD$ corresponds to its right action on
$\sDiff_M(\cO_M,\cO_M;P)$ coming from the left action of differential operators on the second argument.
}

\

\noindent
{\bf Proof.}  Suppose $g, h \in \cO_{M}$,
  then 
\[
((A\otimes B)\cdot f)(g,h) = (fA\otimes
b)(g,h) = B(A(g)f\cdot h) = B(A(g)(fh))=(A\otimes B)(g,fh). 
\]
In addition, if $\tau \in \mathcal{T}_{M}$, we have
$\tau(A(g)h) =(\tau A(g))h + A(g)\tau(h)$, and therefore 
\begin{multline*}
((A\otimes B)\cdot \tau)(g,h) = (A\otimes(B\tau) -(\tau A)\otimes B)(g,h) =
B\tau(A(g)h) - B((\tau A(g))h) \\
= B(A(g)\tau(h))=(A\otimes B)(g,\tau(h)).
\end{multline*}
\ \hfill $\Box$
\

\noindent

The coaction $\bDelta_P$ is a homomorphism of right $\cD_M$-modules
that sends $A\in P_\cD$ to the operator
$\bDelta_P(A)\in\sDiff_M(\cO_M,\cO_M;P)$ given by
\[
\bDelta_P(A):(f,g)\mapsto fA(g),\qquad f,g\in\cO_M.
\]

\

\medskip

\punkt {\bf Remark.} \label{rem:coaction_on_linear}
Note that the sub-$\cO_M$-module $P\subset P_\cD$ generates the induced right $\cD_M$-module $P_\cD$.
Since the coaction
\[ 
\bDelta_P:P_\cD\to\cD_M\otimes_{\cO_M} P_\cD
\]
is a morphism of right $\cD_M$-modules, it is uniquely determined by its restriction to $P\subset P_\cD$.
Under the identification $P_\cD=\sDiff_M(\cO_M;P)$, a section $p\in P\subset P_\cD$ corresponds to the differential 
operator
\[p(f)=fp,\qquad f\in \cO_M,\]
and therefore the operator $\bDelta_P(A)\in\sDiff_M(\cO_M,\cO_M;P)$ is given by
\[\bDelta_P(p):(f,g)\mapsto fgp,\qquad f,g\in\cO_M.\]
Under the identification $\sDiff_M(\cO_M,\cO_M;P)=\cD_M\otimes P_\cD$, this corresponds to
\[\bDelta_P(p)=1\otimes p\in\cD_M\otimes P_\cD,\qquad (p\in P\subset P_\cD).\]

\

\medskip

\punkt \ \label{sss:inducedD} 
Let now $\qu{M}$ be a neutralized $\star$-quantization. To an
$\cO_\qu{M}$-module $\qu{P}$ we assign the right $\cD_M\otimes_\C R$-module 
$\qu{P}_\cD:=\sDiff_M(\cO_M\otimes_\C R;\qu{P})$.  Note that
the right $\cD_M\otimes_\C R$-module
$(\cO_\qu{M})_\cD=\sDiff_M(\cO_M\otimes_\C R;\cO_\qu{M})$ also has
commuting right and left actions of $\cO_\qu{M}$, and we have a
natural identification
\[
\qu{P}_\cD = (\cO_\qu{M})_\cD\otimes_{\cO_\qu{M}}\qu{P}.
\]
Since $(\cO_\qu{M})_\cD$ is a locally free right $\cO_\qu{M}$-module,
we see that the assignment $\qu{P}\to\qu{P}_\cD$ is an exact functor.
Note also that $(\cO_\qu{M})_\cD$ is the dual of $\cA_\qu{M}$ in the
sense that
\[
(\cO_\qu{M})_\cD = \sHom_{\cD_M\otimes_\C
R}(\cA_\qu{M},\cD_M\otimes_\C R).
\]
The right $\cD_M$-module $\qu{P}_\cD$ is naturally a comodule over the
coalgebra $\cA_\qu{M}$.  The coaction can be written in the same way
as in \ref{ssec:coaction}: the tensor product $\cA_{\qu{M}}\otimes
\qu{P}_\cD$ (taken over $\cO_M\otimes_\C R$) is identified with the
sheaf of bidifferential operators $\sDiff_M(\cO_\qu{M},\cO_M\otimes_\C
R;\qu{P})$ via
\[
A\otimes B \mapsto \left((\qu{f},g)\mapsto B(A(\qu{f})g)\right)\qquad
(A\in\cA_\qu{M},
B\in\qu{P}_\cD,\qu{f}\in\cO_\qu{M},g\in\cO_M\otimes_\C R). 
\]
The coaction 
\[
\bDelta_\qu{P}:\qu{P}_\cD\to\cA_{\qu{M}}\otimes \qu{P}_\cD
\] 
is a homomorphism of right 
$\cD_M\otimes_\C R$-modules that sends $A\in \qu{P}_\cD$
to the operator \linebreak 
$\bDelta_\qu{P}(A)\in\sDiff(\cO_\qu{M},\cO_M\otimes_\C
R;\qu{P})$ given by 
\[
\bDelta_\qu{P}(A):(\qu{f},g)\mapsto \qu{f}A(g),\qquad
\qu{f}\in\cO_\qu{M},g\in\cO_M\otimes_\C R.
\]

\

\medskip

\punkt {\bf Proposition.}\label{prop:absolute_comodules} {\it

\begin{itemize} \item[\bf(a)] The assignment $P\mapsto
\left(P_{\cD},\bDelta_P\right)$ is an equivalence between the category
of $\cO_M$-modules $\Mod{\cO_M}$ and the category of $\cD_M$-comodules
$\Comod{\cD_M}$.  (Recall that `comodule' is short for `left comodule
in the category of right $\cD_M$-modules'.)

\item[\bf(b)] More generally, the functor 
\[\qu{P}\mapsto\left(\qu{P}_{\cD}, \bDelta_\qu{P}\right):
\Mod{\cO_\qu{M}}\to\Comod{\cA_{\qu{M}}}\] 
is an equivalence.  
\end{itemize} 
} 

\noindent
{\bf Proof.} We postpone the proof of {\bf(a)} to avoid repetitions:
it is a particular case of
Proposition~\ref{prop:relative_comodules}. To derive {\bf(b)} from
{\bf(a)} fix an $\cA_\qu{M}$-comodule $F$ and consider
the sheaf of homomorphisms 
\[
\xi(F) := \sHom_{\cA_{{\qu{M}}}}((\cO_\qu{M})_\cD,F).
\]  
Note that the right action of $\cO_\qu{M}$ on itself induces a right
action of $\cO_\qu{M}$ on $(\cO_\qu{M})_\cD$; therefore, $\xi(F)$ is
naturally a $\cO_\qu{M}$-module.  Let us prove that the functors
$F\mapsto\xi(F)$ and $\qu{P}\mapsto\qu{P}_\cD$ are mutually inverse.

Recall that
$\qu{P}_\cD=(\cO_\qu{M})_\cD\otimes_{\cO_\qu{M}}\qu{P}$. We therefore
obtain functorial morphisms $\xi(F)_\cD\to F$ and
$\qu{P}\to\xi(\qu{P}_\cD)$ (essentially because the two functors are
adjoint). We need to prove that the maps are
isomorphisms. By Nakayama's lemma (see
\cite[Proposition~6.3]{Lowen_VanDerBergh}) it is enough to prove this
when $F$ and $\qu{P}$ are annihilated by $\hbar$, in which case it
reduces to statement {\bf(a)}.  \ \hfill $\Box$

\

\medskip

\punkt {\bf Remarks.} \ {\bf (1)} \
Propositions~\ref{prop:absolute_comodules} and
Proposition~\ref{prop:relative_comodules} may be viewed as a form of
Koszul duality, similar {to \cite{positselski.ams}}.

{\bf (2)} \
It follows from {the proof of} 
Proposition~\ref{prop:absolute_comodules} that the inverse of the
functor $\qu{P}\to\qu{P}_\cD$, which we denoted by $\xi$, can be also
given as
\[
\xi(F)=F\otimes_{\cD_M}\cO_M.
\]
However, the description of the $\cO_\qu{M}$-module structure on
$\xi(F)$ is more complicated in this approach. This construction shows
that the $R$-module structure on $\xi(F)$ depends only on the
$\cD_M\otimes_\C R$-module structure on $F$ and not on the
$\cA_\qu{M}$-comodule structure.

\

\medskip

\noindent
We need the following observation about induced
$\cD$-modules. 
\

\smallskip

\punkt {\bfseries Lemma.}\label{lem:inducedflatness}
 {\it Suppose $P\in\Mod{\cO_M}$ and $F\in\Mod{\cD_M}$. The natural
embedding $P\to P_\cD$ induces a morphism
$P\otimes{_{\cO_{M}}} F\to
 P_\cD\otimes{_{\cO_{M}}} F$,
and therefore a morphism of right $\cD_M$-modules 
\[i:(P\otimes{_{\cO_{M}}} F)_\cD\to
 P_\cD\otimes{_{\cO_{M}}} F.\] 
We claim that $i$ is an isomorphism in $\Comod{\cD_M}$.}

\noindent
{\bf Proof.} It is easy to check that $i$ is an isomorphism of
right $\cD_{M}$-modules by looking at the {order}
 filtrations on
$(P\otimes{_{\cO_{M}}} F)_\cD$ 
and on $P_\cD$. It remains to show that $i$
respects the coaction of $\cD_M$; that is, we need to check that the
two compositions
\[
\xymatrix@-1pc{
(P\otimes{_{\cO_{M}}} F)_\cD\ar[r]^i\ar[d] &
  P_\cD\otimes{_{\cO_{M}}} F\ar[d]\\ 
\cD_M\otimes{_{\cO_{M}}}(P\otimes{_{\cO_{M}}}
F)_\cD\ar[r] & \cD_M\otimes{_{\cO_{M}}}
P_\cD\otimes{_{\cO_{M}}} F} 
\] 
coincide. Since the compositions are $\cD_M$-linear, it suffices to
check that they coincide on $P\otimes{_{\cO_{M}}}
F\subset (P\otimes{_{\cO_{M}}} F)_\cD$, because this sub-$\cO_M$-module
generates the induced $\cD_M$-module $(P\otimes{_{\cO_{M}}} F)_\cD$. We claim
that both compositions send $A\otimes B\in
P\otimes_{\cO_{M}} F$ to \[1\otimes A\otimes B\in\cD_M\otimes_{\cO_M} P\otimes_{\cO_M} F\subset \cD_M\otimes_{\cO_M}
P_\cD\otimes_{\cO_M} F.\]

Indeed, the left arrow sends $A\otimes B$ to $1\otimes A\otimes B$ by Remark~\ref{rem:coaction_on_linear}.
The right arrow equals $\bDelta\otimes{id}_F$; therefore, it sends $A\otimes B$ to
$\bDelta(A)\otimes B$, which is also equal to $1\otimes A\otimes B$ by Remark~\ref{rem:coaction_on_linear}.
This completes the
proof.  \ \hfill $\Box$

\

\medskip

\punkt \ Consider now a relative version of the assignment from
Proposition~\ref{prop:absolute_comodules}.  Let $p_X:Z\to X$ be a
submersive morphism of complex manifolds.  Recall that we view the
sheaf of differential operators $\cD_X$ as a coalgebra in
$\Mod{\cD_X}$.  Its inverse image $p_X^*\cD_X$ is a coalgebra in
$\Mod{\cD_Z}$. There is a natural identification
\[
p_X^*\cD_X=\sDiff_Z(p_X^{-1}\cO_X;\cO_Z),
\] 
where the coalgebra structure is induced from the product on $\cO_X$. 
This identification can be used to study $p_X^*\cD_X$-comodules.

For any $\cO_Z$-module $P$, the induced right $\cD_Z$-module $P_\cD$
is naturally equipped with a $p_X^*\cD_X$-comodule structure as
in~\ref{ssec:coaction}. Equivalently, consider the natural map
$\cD_Z\to p_X^*\cD_X$ given by restricting differential operators
$\cO_Z\to\cO_Z$ to $p_X^{-1}\cO_X$. It is a morphism of coalgebras,
and the $p_X^*\cD_X$-comodule structure on $P_\cD$ is induced by the
coaction of $\cD_Z$.

In particular, $\cD_Z=(\cO_Z)_\cD$ is a $p_X^*\cD_X$-comodule. Notice
that it also carries a commuting left action of the algebra of
relative differential operators $\cD_{Z/X}\subset\cD_Z$. We can
therefore factor the functor $P\to P_\cD$ as a composition of two
functors. First, an $\cO_Z$-module $P$ induces a right
$\cD_{Z/X}$-module $P_{\cD/X}=\sDiff_{Z/X}(\cO_X;P)$. Second, a right
$\cD_{Z/X}$-module $F (=P_{\cD/X})$ induces a right $\cD_Z$-module
$F\otimes_{\cD_{Z/X}}\cD_Z$ that carries a coaction
$\bDelta_{F/p_{X}^{*}\cD_{X}}$ of $p_X^*\cD_X$.

\

\medskip 

\punkt {\bf Proposition.} \label{prop:relative_comodules} {\it The
assignment $F\mapsto \left(F\otimes_{\cD_{Z/X}}\cD_Z,
\bDelta_{F/p_{X}^{*}\cD_{X}}\right)$ is an equivalen{
\[
\Mod{\cD_{Z/X}^{op}}\stackrel{\cong}{\to} \Comod{p_X^*\cD_X}
\]
between the category of right $\cD_{Z/X}$-modules and the category
$\Comod{p_X^*\cD_X}$ of comodules over $p_X^*\cD_X$ in the category of
right $\cD_{Z}$-modules.
}
} 

\

\noindent
{{\bfseries Remark:}  Note that
  Proposition~\ref{prop:absolute_comodules} is a special case (an
  absolute version) of Proposition~\ref{prop:relative_comodules} (a
  relative version).}

\noindent
{\bf Proof.} The question is essentially local, because both categories
are global sections of the corresponding stacks of categories, and 
the functor between them is induced by a $1$-morphism of stacks.
We may therefore assume that
$Z$ is an open subset of $\C^{n+m}$ with coordinates
$(x_1,\dots,x_n;y_1,\dots,y_m)$ and that $p_X:Z\to X$ is given by
$p_X(x_1,\dots,x_n;y_1,\dots,y_m)=(x_1,\dots,x_n)\in X\subset\C^n.$
Let us verify that the inverse functor sends $G\in\Comod{p_X^*\cD_X}$
to $\sHom_{p_X^*\cD_X}(\cD_Z,G)$.

For a multiindex $\alpha=(\alpha_1,\dots,\alpha_n)\in(\Z_{\ge0})^n$,
we set
\[
\partial^\alpha := 
\frac{\partial^{\alpha_1}}{\partial x_1^{\alpha_1}}\cdots
\frac{\partial^{\alpha_n}}{\partial x_n^{\alpha_n}},\quad \alpha!
:= \alpha_1!\cdots\alpha_n!,
\quad |\alpha|=\alpha_1+\dots+\alpha_n.
\]
A section $A\in p_X^*\cD_X$ can be written as 
\[
A =
\sum_{\alpha\in(\Z_{\ge0})^n}\frac{1}{\alpha!}a_\alpha{(x,y)}
\partial^\alpha, 
\]
where locally on $Z$, all but finitely many $a_\alpha(x,y)\in\cO_Z$ vanish.
(The factor of $\alpha!$ is included to avoid binomial coefficients in
the coproduct.) The coproduct $\bDelta:p_X^*\cD_X\to p_X^*\cD_X\otimes
p_X^*\cD_X$ is given by
\[
\bDelta\left( 
\sum_{\alpha\in(\Z_{\ge0})^n}
\frac{1}{\alpha!}a_\alpha(x,y)\partial^\alpha\right) =
\sum_{\beta,\gamma\in(\Z_{\ge0})^n}
\frac{1}{\beta!\gamma!}a_{\beta+\gamma}(x,y)
\partial^\beta\otimes\partial^\gamma. 
\]

Let $F$ be a right $\cD_{Z/X}$-module. Then $G=F\otimes_{\cD_{Z/X}}\cD_Z$
can be described as
\[
G=F\otimes_\C\C\left[\frac{\partial}{\partial
    x_1},\dots,\frac{\partial}{\partial x_n}\right].
\]
Note that $F$ is identified with a subsheaf of $G$, and that $F$
generates $G$ as a right $\cD_Z$-module. 
{Note moreover that the coproduct
$\bDelta_{F/p_{X}^{*}\cD_{X}}:G\to p_X^*\cD_X\otimes G$ has the
property $\bDelta_{F/p_{X}^{*}\cD_{X}}(s)=1\otimes s$ for  $s\in
F\subset G$ (even though $\bDelta_{F/p_{X}^{*}\cD_{X}}(A)\neq
1\otimes A$ for general $A\in G$). 

{This property characterizes
  $\bDelta_{F/p_{X}^{*}\cD_{X}}$ uniquely since  $F$ generates $G$ as a right
  $\cD_Z$-module and $\bDelta_{F/p_{X}^{*}\cD_{X}}$ is $\cD_Z$-linear.}

Let us describe explicitly the functor in the
opposite direction. Let $G$ be a right $\cD_Z$-module.  A morphism
$\bDelta_G : G\to p_X^*\cD_X\otimes G$ can be written in the form
\[
\bDelta_G(s)=\sum_{\alpha\in(\Z_{\ge0})^n}\frac{1}{\alpha!}
\partial^\alpha\otimes\bDelta_\alpha(s), 
\] 
where $\bDelta_\alpha : G\to G$ are maps such that for given $s\in G$,
only finitely many $\bDelta_\alpha(s)$ are non-zero locally on $Z$. Then $\bDelta_G$
is a coaction if and only if
$\bDelta_\alpha\circ\bDelta_\beta=\bDelta_{\alpha+\beta}$ for any
$\alpha,\beta\in(\Z_{\ge0})^n$ and $\bDelta_0$ is the identity map. In
other words, a coaction $\bDelta_G$ can be written as
\[
\bDelta_G(s)=\sum_{\alpha\in(\Z_{\ge0})^n}\frac{1}{\alpha!}
\partial^\alpha\otimes
(\bDelta_1^{\alpha_1}\cdots\bDelta_n^{\alpha_n})(s)
\] 
for some
commuting operators $\bDelta_1,\dots,\bDelta_n : G\to G$ that are
locally nilpotent (that is, every $s\in G$ is annihilated by some
power of the operators locally on $Z$).

It remains to impose the restriction that $\bDelta_G$ is a morphism of
right $\cD_Z$-modules.  This results in the following conditions on
$\bDelta_i$:
\begin{align*}
\bDelta_i(sf) & = \bDelta_i(s)f\qquad&(&f\in\cO_Z,1\le i\le n,s\in
G)\\[+1pc] 
\bDelta_i\left(s\frac{\partial}{\partial y_j}\right) & =
\bDelta_i(s)\frac{\partial}{\partial y_j}\qquad&(&1\le i\le n,1\le
j\le m,s\in G) \\[+1pc] 
\bDelta_i\left(s\frac{\partial}{\partial x_j}\right) & =
\bDelta_i(s)\frac{\partial}{\partial x_j}- 
\begin{cases} s,\text{ if }i=j \\[+1pc] 
0,\text{ if }i\ne j\end{cases}
\qquad&(&1\le i\le n,1\le j\le n,s\in G).
\end{align*}
In other words, the operators $\bDelta_i$ commute with $\cD_{Z/X}$,
while the commutativity relations for operators $\bDelta_i$,
$\frac{\partial}{\partial x_j}$ are those of the Weyl algebra
\[
W=\C[t_1,\dots,t_n]\left\langle\frac{\partial}{\partial
  t_1},\dots,\frac{\partial}{\partial t_n}\right\rangle.
\]
Since the operators $\bDelta_i$ are locally nilpotent, Kashiwara's
Lemma {\cite[Theorem~4.30]{kashiwara-dmodules}} implies
that $G$ decomposes 
as a tensor product
\[
G = F\otimes_\C
\C\left[\frac{\partial}{\partial x_1},\dots,\frac{\partial}{\partial
x_n}\right],
\] 
for 
\[
F :=\{s\in G :\bDelta_i(s)=0,\text{ for all
}i=1,\dots,n\}=\{s\in G:\bDelta(s)={1\otimes s}\}.
\] 
Moreover, $F$ is invariant under $\cD_{Z/X}$. Note that 
$F$ is identified with $\sHom_{p_X^*\cD_X}(\cD_Z,G)$. It
is now clear that the two functors are inverse to each other.
\ \hfill $\Box$

\

\medskip

\punkt {\bfseries Remark.}\label{rem:relqudiffops_flat}
Proposition~\ref{prop:relative_comodules} implies
Proposition~\ref{prop:relqudiffops_flat}.  Recall that to prove
Proposition~\ref{prop:relqudiffops_flat}, we had to verify exactness
of a sequence
\[
\sDiff_Z(\cO_Z;\cO_Z)\to\sDiff_Z(p_X^{-1}\cO_X,\cO_Z;\cO_Z) \to
\sDiff_Z(p_X^{-1}\cO_X,p_X^{-1}\cO_X,\cO_Z;\cO_Z).
\] 
It {can be seen} that this sequence is the
beginning of the cobar complex computing the sheaves
$\sExt^i_{p_X^*\cD_X}(\cD_Z,\cD_Z)$ of local $\Ext$'s in the category
$\Comod{p_X^*\cD_X}$.

Under the equivalence of Proposition~\ref{prop:relative_comodules},
$\cD_Z\in\Comod{p_X^*\cD_X}$ 
corresponds to $\cD_{Z/X}\in\Mod{\cD_{Z/X}}$. Therefore,
\[
\sExt^i_{p_X^*\cD_X}(\cD_Z,\cD_Z) =
\sExt^i_{\cD_{Z/X}}(\cD_{Z/X},\cD_{Z/X}) = 
\begin{cases}\cD_{Z/X},&i=0\cr 0,&i\ne 0.\end{cases}
\]
Proposition~\ref{prop:relqudiffops_flat} follows.

\

\medskip

\punkt \label{ssec:relativecomodules} Finally, given a neutralized
$\star$-quantization $\qu{X}$ of $X$, the inverse image
\[
p_X^*\cA_\qu{X} = \sDiff_Z(p_X^{-1}\cO_\qu{X};\cO_Z\otimes_\C R)
\] 
is a coalgebra in the category of left $\cD_Z\otimes_\C R$-modules.
Let $\qu{Z}$ be a neutralized $\star$-quantization of $Z$ such that
$p$ extends to a morphism $\qu{p}:\qu{Z}\to\qu{X}$; in other words,
the morphism of algebras $p_X^{-1}\cO_X\to\cO_Z$ extends to a morphism
$p_X^{-1}\cO_\qu{X}\to\cO_\qu{Z}$. Then for any $\cO_\qu{Z}$-module
$\qu{P}$, the right $\cD_Z\otimes_\C R$-module $\qu{P}_\cD$ is
naturally equipped with a $p_X^*\cA_\qu{X}$-comodule structure.

Repeating the proof of Proposition~\ref{prop:absolute_comodules}(b),
it is easy to show that the category of $p_X^*\cA_\qu{X}$-comodules is
equivalent to the category of right $\cD_{\qu{Z}/\qu{X}}$-modules.
Moreover, the functor $\qu{P}\to\qu{P}_\cD$ decomposes as the
composition of the induction functor
\[
\qu{P} \to \cD_{\qu{Z}/\qu{X}}\otimes_{\cO_{\qu{Z}}}\qu{P} =
\{A\in\sDiff_Z(\cO_\qu{Z};\qu{P}):A(\qu{f}\qu{g}) = \qu{f}A(\qu{g}) \text{
for all }\qu{f}\in p_X^{-1}\cO_\qu{X},\qu{g}\in\cO_\qu{Z}\} 
\]
(the functor appears in the end of Section~\ref{ssec:qinduced_D}) and
this equivalence.  We do not need these statements, so we leave the
details to the reader. The statements are useful because they allow us
to relate two proofs of Theorem~\ref{thm:obstructions}: the proof
using comodules (Section~\ref{ssec:obstructionsneutr}) 
and the proof using relative differential
operators on quantizations
(Sections~\ref{ssec:obstructionsneutr'}).

\

\medskip

\punkt {\bf Remark.} \label{rem:star_operations} As we already
mentioned, the results of Section~\ref{ssec:quantcoalgebras} are
inspired by the $\star$-pseudotensor structure of
\cite{chiralalgebras}. {For instance,
  Proposition~\ref{prop:absolute_comodules} can be extracted from the remark in
  \cite[Section~2.1.8]{chiralalgebras}.}  Let us make the relation more explicit. 

Let $M$ be a complex manifold. Recall (\cite[2.2.3]{chiralalgebras})
that for right $\cD_M$-modules $F, F\q1, \linebreak \dots, F\q{n}$, the
space of {\bfseries $\star$-operations} $F\q1\times\dots\times
F\q{n}\to F$ is defined to be
\[\op{Hom}(\boxtimes_{i=1}^n F\q{i},\Delta_* F).\]
Here the homomorphisms are taken in the category of right
$\cD$-modules on $M^n$, $\Delta:M\to M^n$ is the diagonal embedding,
and $\Delta_*$ is the $\cD$-module direct image functor.  In
particular, let $\qu{P},\qu{P}\q1,\dots,\qu{P}\q{n}$
be $\cO$-modules on neutralized $\star$-quantizations
$\qu{M},\qu{M}\q1,\dots,\qu{M}\q{n}$ of $M$. Consider the induced right
$\cD_M\otimes_\C R$-modules $\qu{P}\q1_\cD\times\dots\times\qu{P}\q{n}_\cD\to\qu{P}_\cD$.
Then the space of $R$-linear $\star$-operations
$\qu{P}\q1_\cD\times\dots\times\qu{P}\q{n}_\cD\to\qu{P}_\cD$ is
identified with the space of polydifferential operators 
$\qu{P}\q{1}\times\dots\times\qu{P}\q{n}\to\qu{P}$
{(see Section~\ref{sss:inducedD})}.

Let $\qu{M}$ be a neutralized $\star$-quantization of $M$. The induced
right $\cD_M\otimes_\C R$-module $(\cO_{\qu{M}})_\cD$ is an associative${}^\star$
algebra (that is, the product is a $\star$-operation). Moreover,
$(\cO_{\qu{M}})_\cD$ is an $R$-deformation of $(\cO_M)_\cD=\cD_M$;
this gives a one-to-one correspondence between neutralized
$\star$-quantizations and deformations of the associative${}^\star$
algebra $\cD_M$.

Note that the associative${}^\star$ algebra $(\cO_{\qu{M}})_\cD$ and
the coalgebra $\cA_\qu{M}$ are related by the duality
(cf. \cite[2.5.7]{chiralalgebras}):
\[
\cA_{\qu{M}}=\sHom_{\cD_M}((\cO_{\qu{M}})_\cD,\cD_M).
\]
The correspondence identifies $(\cO_{\qu{M}})_\cD$-modules and
$\cA_\qu{M}$-comodules.

Similar approach works in the relative situation. Recall that we
consider, for a submersive morphism $p_X:Z\to X$ and a neutralized
$\star$-quantization $\qu{X}$ of $X$, the category of comodules over
the coalgebra $p_X^*\cA_\qu{X}$. Notice however that a pull-back of an
associative${}^\star$ algebra carries a $\star$-operation only after a
cohomological shift: the category of $p_X^*\cA_\qu{X}$-comodules is
equivalent to the category of modules over the associative${}^\star$
algebra $p_X^*(\cO_{\qu{X}})_\cD[\dim(X)-\dim(Z)]$.

\section{Deformation theory} \label{sec:defo}

This section contains the proof of Theorem~\ref{thm:obstructions}.
The statement of the theorem is local in $Y$; we can therefore assume
that $Y$ is a Stein manifold.

\subsection{Proof of Theorem~\ref{thm:obstructions}: neutralized
  case}\label{ssec:obstructionsneutr} 

In this section, we assume that the quantization $\gb{X}$ of $X$ is
neutralized; this simplifies the argument. In
Section~\ref{ssec:gerbydefo}, we explain how to modify the argument
for arbitrary $\star$-quantizations.

\

\medskip

\punkt Let $\gb{Y}$ be a $\star$-quantization of $Y$. By definition, the
reduction of $\gb{Y}$ to $\C$ is neutralized. This neutralization can be 
lifted to a neutralization of $\gb{Y}$: indeed, once such a lift is constructed modulo $\hbar^k$,
the obstructions to lifting modulo $\hbar^{k+1}$ lie in $H^2(Y,\cO_Y)$, which
vanishes because $Y$ is Stein. Similarly, the vanishing of $H^1(Y,\cO_Y)$ implies
that the lifting is unique up to a non-unique automorphism. 
Thus, there exists a global neutralization 
$\alpha\in\gb{Y}(Y)$ that agrees with the neutralization of the reduction of 
$\gb{Y}$ to $\C$. Such
$\alpha$ defines a neutralized $\star$-quantization $\qu{Y}_\alpha$ of
$Y$.

Moreover, $\alpha$ is unique up to a non-unique isomorphism.  Note
that automorphisms of $\alpha$ are required to act trivially modulo
$\hbar$.  The automorphism group of $\alpha$ is identified with the
group
\[
\{\qu{f}\in H^0(Y,\cO_{\qu{Y}_\alpha}):\qu{f}=1\mod \hbar\},
\]
which acts on $\cO_{\qu{Y}_{\alpha}}$ by conjugation.

We can now restate Theorem~\ref{thm:obstructions} using neutralized
$\star$-quantizations.  In other words, we replace pairs
$(\gb{Y},\qu{P})$ with triples $(\gb{Y},\qu{P},\alpha)$, where
$\alpha$ is a neutralization. This leads to the following statement.

\

\medskip

\punkt {\bf Proposition.} \label{prop:obstructionsneutr} {\it Let $X$
and $Y$ be complex manifolds, $P$ a coherent sheaf on $Z=X\times Y$
whose support is proper over $Y$. Assume that $Y$ is Stein and that
the map $\iota:\cO_Y\to\cE(P)$ (defined in
Section~\ref{ssec:plan_of_proof}) is a quasi-isomorphism. Let $\qu{X}$
be a neutralized $\star$-quantization of $X$.
\begin{enumerate}
\item[\bf(a)] There exists a neutralized $\star$-quantization
$\qu{Y}$ of $Y$ and a deformation of $P$ to an $\cO$-module $\qu{P}$
on $\qu{X}\times\qu{Y}^{op}$. 

\item[\bf(b)] The pair $(\qu{Y},\qu{P})$ is unique up to a non-unique
isomorphism.

\item[\bf(c)] The automorphism group 
$\Aut(\qu{Y},\qu{P})$ is identified with the group
\[
\{\qu{f}\in H^0(Y,\cO_\qu{Y}):\qu{f}=1\mod \hbar\},
\] 
which acts on $\cO_{\qu{Y}}$ by conjugation and on $\qu{P}$ using the
structure of an $\cO$-module.
\end{enumerate}}

\

\medskip

\punkt \ Consider the right $\cD_Z$-module $P_\cD$ induced by $P$. It
carries a natural coaction of $\cD_Z$. Alternatively, we view this
coaction as a pair of commuting coactions of $p_X^*\cD_X$ and
$p_Y^*\cD_Y$ using the identification
$\cD_Z=p_X^*\cD_X\otimes p_Y^*\cD_Y$. We claim that the coaction
morphism
\[
\bDelta_Y=\bDelta_{P,Y}:P_\cD\to {p_Y^*\cD_Y\otimes
  P_\cD} 
\]
has the following universal property.

\

\medskip

\punkt {\bf Lemma.} \label{lem:coproduct_universality} {\it In the
assumptions of Proposition~\ref{prop:obstructionsneutr}, let $F$ be
any left $\cD_Y$-module that is flat over $\cO_Y$, and let 
\[
\phi:P_\cD\to {p_Y^* F\otimes P_\cD}
\]
be a morphism of right $\cD_Z$-modules that commutes with the coaction
of $p_X^*\cD_X${, i.e. $\phi$ is also a morphism in the
category $\Comod{p_X^*\cD_X}$}. Then there exists a unique morphism of
left $\cD_Y$-modules $\psi:\cD_Y\to F$ such that $\phi=(id\otimes
p_Y^*(\psi))\circ\bDelta_Y$.

This provides an identification between 
$\Hom_{p_X^*\cD_X}(P_\cD,{p_Y^*F\otimes
  P_\cD})$ (morphisms in the 
category $\Comod{p_X^*\cD_X}$) and
$\Hom_{\cD_Y}(\cD_Y,F)=H^0(Y,F)$. The identification remains valid in
the derived sense:
\[
\mathbb{R}\Hom_{p_X^*\cD_X}(P_\cD,{p_Y^*F\otimes
  P_\cD}) \cong R\Gamma(Y,F) 
\]
and hence  
\[
\Ext^i_{p_X^*\cD_X}(P_\cD,{p_Y^* F\otimes
  P_\cD})=H^i(Y,F). 
\]
}

\noindent
{\bf Proof.} By Proposition~\ref{prop:relative_comodules}, the
category of $p_X^*\cD_X$-comodules is equivalent to the category of
right $\cD_{Z/X}$-modules.  Under this equivalence, $P_\cD$
corresponds to $P_{\cD/X}$. We claim that for any $G\in\Mod{\cD_Z}$,
the comodule  ${G\otimes P_\cD}
\in\Comod{p_X^*\cD_X}$ corresponds to ${G\otimes
P_{\cD/X}}\in\Mod{\cD_{Z/X}^{op}}$. Indeed, 
${G\otimes P_\cD\simeq (G\otimes P)_\cD}$ by
Lemma~\ref{lem:inducedflatness}. The isomorphism respects the coaction
of $\cD_Z$, and therefore also the coaction of
$p_X^*\cD_X$. Similarly, we have an isomorphism of induced right
$\cD_{Z/X}$-modules  ${G\otimes P_{\cD/X}\simeq
(G\otimes P)_{\cD/X}}$.  This implies the claim.

In particular, the $p_X^*\cD_X$-comodule 
{$p_Y^*F\otimes  P_\cD$} 
corresponds to {$p_Y^*F\otimes P_{\cD/X}$}.  Therefore,

\[
\Ext^i_{p_X^*\cD_X}(P_\cD, {p_Y^*
F\otimes  P_\cD})=\Ext^i_{\cD_{Z/X}}(P_{\cD/X},
{p_Y^*F\otimes P_{\cD/X}})=
\Ext^i_{\cO_Z}(P,{p_Y^*F\otimes P_{\cD/X}}),
\]
where the second equality follows from Lemma~\ref{lem:rel_induced_adjunct}.
By the projection formula, 
\[
\begin{split}
\RR \Hom_{\cO_Z}(P,{ p_Y^*F\otimes P_{\cD/X}}) & =
\RR\Gamma(Y,{F\otimes \RR
  p_{Y*}\RR\sHom_{\cO_Z}(P,P_{\cD/X})}) \\  
& = \RR\Gamma(Y,{F\otimes \cE(P)})=\RR\Gamma(Y,F),
\end{split}
\]
as required.
\ \hfill $\Box$

\

\medskip

\punkt \
Recall that $\cA_\qu{X} =
    \sDiff(\cO_{\qu{X}};\cO_{X}\otimes R)$ is a coalgebra in 
  $\Mod{\cD_{X}\otimes R}$. It is a deformation of the coalgebra $p_{X}^{*}\cD_{X}\in\Mod{\cD_X}$; that is,
  it is flat over $R$ and its reduction to $\C$ is identified with $p_{X}^{*}\cD_{X}$, see
  Sections~\ref{sssec:coaction} and \ref{sssec:defo.comodules}. Lemma~\ref{lem:coproduct_universality}
  is a universal property of the pair $(P_\cD,\bDelta_Y)$ consisting of an object and a morphism in the category
  $\Comod{p_X^*\cD_X}$. We claim that this pair admits a deformation in the category $\Comod{\cA_{\qu{X}}}$,
  and that the universal property is preserved. Here is a precise statement.

\

\medskip

\punkt {\bf Lemma} \label{lem:quantized_comodules} {\it In the
assumptions of Proposition~\ref{prop:obstructionsneutr}, we have the
following.
\begin{enumerate}
\item[{\bf (a)}] There is a $R$-flat deformation of
$P_\cD\in\Comod{p_X^*\cD_X}$ to $P'_\cD\in\Comod{p_X^*\cA_\qu{X}}$.
The deformation is unique up to isomorphism.  

\item[{\bf (b)}] The coaction $\bDelta_Y$ extends to a homomorphism of
$p_X^*\cA_\qu{X}$-comodules (and, in particular, right
$\cD_Z$-modules) 
\[
\bDelta'_Y:P'_\cD\to  {p_Y^*\cD_Y\otimes_{\cO_Z}
  P'_\cD} = {
(p_Y^*\cD_Y\otimes_\C R)\underset{\cO_Z\otimes_\C R}{\otimes} P'_\cD}.
\]

\item[{\bf (c)}] The homomorphism $\bDelta'_Y$ is universal in the
following sense: given any $\cD_Y\otimes_\C R$-module $F$ that is flat
over $\cO_Y\otimes_\C R$ and any homomorphism of
$p_X^*\cD_X$-comodules 
\[
\phi:P'_\cD\to  
{p_Y^* F\underset{\cO_Z\otimes_\C R}{\otimes} P'_\cD},
\]
there exists a unique morphism of $\cD_Y\otimes_\C R$-modules
$\psi:\cD_Y\otimes_\C R\to F$ such that $\phi=(id\otimes
p_Y^*(\psi))\circ\bDelta'_Y$.
\end{enumerate}}

\noindent
{\bf Proof.} Since $Y$ is Stein, $H^{1}(\cO_{Y})
  = H^{2}(\cO_{Y}) = 0$. By Lemma~\ref{lem:coproduct_universality}, 
we see that
\[
\Ext^1_{p_X^*\cD_X}(P_\cD,P_\cD)=\Ext^2_{p_X^*\cD_X}(P_\cD,P_\cD)=0.
\]
Now part {\bf (a)} follows from Proposition~\ref{prop:defcomodules}.

To prove part {\bf (b)}, note that obstructions to extending 
{$\bDelta_{Y}$} to  
$\bDelta'_Y$ belong to 
\[
\Ext^1_{{p_X^*\cA_\qu{X}}}(P'_\cD,
{p_Y^*\cD_Y\otimes  \hbar  P'_\cD}). 
\] 
Using the filtration of {$p_Y^*\cD_Y\otimes
  \hbar P'_\cD$} 
by powers of $\hbar$, we see that it suffices to verify the vanishing
  of 
\[
\Ext^1_{{p_X^*\cA_\qu{X}}}(P'_\cD,{
p_Y^*\cD_Y\otimes
(\hbar^k P'_\cD/\hbar^{k+1}
P'_\cD)}) =
\Ext^1_{p_X^*\cD_X}(P_\cD, {p_Y^*\cD_Y\otimes  P_\cD}).
\]
By Lemma~\ref{lem:coproduct_universality}, the vanishing is equivalent to
the vanishing of $H^1(Y,\cD_Y)$, which holds because $Y$ is Stein and $\cD_{Y}$
is a colimit of coherent $\cO_{Y}$-modules.

In part {\bf (c)}, the correspondence $\psi\mapsto
{(p_Y^*(\psi)\otimes \op{id})}\circ\bDelta'_Y$ defines
a map
\[
H^0(Y,F)=\Hom_{\cD_Y\otimes_\C R}(\cD_Y\otimes_\C
R,F)\to\Hom_{{p_X^*\cA_{\qu{X}}}}
(P'_\cD,F\underset{\cO_Z\otimes_\C
R}\otimes P'_\cD),
\]
and we need to verify that it is bijective. Using the filtration of
$F$ by the submodules $\hbar^k F$, we see that it is enough to verify
bijectivity of the map
\[
H^i(Y,F/\hbar F)\to\Ext^i_{{p_X^*\cA_{\qu{X}}}}\left(P'_\cD,
{(F/\hbar F)\underset{\cO_Z\otimes_\C R}{\otimes}
P'_\cD}\right).\] 
However,
\begin{multline*}
\Ext^i_{{p_X^*\cA_{\qu{X}}}}\left(P'_\cD,
{(F/\hbar F)\underset{\cO_Z\otimes_\C R}{\otimes}
P'_\cD}\right)=\Ext^i_{{p_X^*\cA_{\qu{X}}}}\left(P'_\cD,
(F\underset{\cO_Z\otimes_\C R}{\otimes}
P'_\cD)\otimes_R \C\right)=\\
\Ext^i_{p_X^*\cD_X}\left(P'_\cD\otimes_R\C,
(F\underset{\cO_Z\otimes_\C R}{\otimes}
P'_\cD)\otimes_R \C\right)=
\Ext^i_{p_X^*\cD_{X}}\left(P_\cD, (F/\hbar
  F)\underset{\cO_Z}{\otimes} 
P_\cD\right), 
\end{multline*}
and the required identity follows from Lemma~\ref{lem:coproduct_universality}. 
\ \hfill $\Box$

\

\medskip

\punkt \ Let us prove Proposition~\ref{prop:obstructionsneutr},
starting with part {\bf (a)}. Taking into account
Proposition~\ref{prop:quantization_coalgebras} and
Proposition~\ref{prop:absolute_comodules}, we need to construct a
deformation of the coalgebra $\cD_Y$ to a coalgebra $\cA_\qu{Y}$ and a
deformation of $P_\cD$ to a right $\cD_Z\otimes_\C R$-module $P'_\cD$
equipped with commuting coactions of $p_X^*\cA_\qu{X}$ and
$p_Y^*\cA_\qu{Y}^{op}$.

Lemma~\ref{lem:quantized_comodules}(a) provides a
$p_X^*\cA_\qu{X}$-comodule $P'_\cD$. Consider the morphism $\Delta'_Y$
from \linebreak Lemma~\ref{lem:quantized_comodules}(b). It remains to
show that there is a coalgebra structure on $\cD_Y\otimes_\C R$ such
that $\Delta'_Y$ becomes a coaction; we can then let $\cA_\qu{Y}^{op}$
be $\cD_Y\otimes_\C R$ with this coalgebra structure. Consider the
composition
\[
P'_\cD\to P'_\cD\otimes p_Y^*(\cD_Y\otimes_\C R)\to P'_\cD\otimes
p_Y^*(\cD_Y\otimes_\C R)\otimes p_Y^*(\cD_Y\otimes_\C R)
\]
(here and elsewhere in the proof, the tensor products are by default
over $\cO_Z\otimes_\C R$). 
It is a morphism of $p_X^*\cA_X$-comodules, so by
Lemma~\ref{lem:quantized_comodules}(c), it  
corresponds to a morphism
\[
\delta':\cD_Y\otimes_\C R\to(\cD_Y\otimes_\C
R)\otimes_{\cO_Y\otimes_\C R}(\cD_Y\otimes_\C R).
\]
It is easy to check that $\delta'$ is indeed a coassociative coproduct
on $\cD_Y\otimes_\C R$ by using the uniqueness part of
Lemma~\ref{lem:quantized_comodules}(c).

\

\medskip

\punkt \ To prove part {\bf (b)} of
Proposition~\ref{prop:obstructionsneutr}, we investigate uniqueness
properties of $P'_\cD$ and $\cA_\qu{Y}$. Let $P\df_\cD$ be another
deformation of $P_\cD$ to a right $\cD_Z\otimes_\C R$-module equipped
with commuting coactions of $p_X^*\cA_\qu{X}$ and
$p_Y^*\cA_{\qu{Y}\df}^{op}$ for some deformation $\cA_{\qu{Y}\df}$ of
the coalgebra $\cD_Y$. Let us construct an isomorphism
$(P'_\cD,\cA_\qu{Y})\simeq (P\df_\cD,\cA_{\qu{Y}\df})$.

Lemma~\ref{lem:quantized_comodules}(a) provides an isomorphism of
$p_X^*\cA_\qu{X}$-comodules $\phi:P\df_\cD\to P'_\cD$. Using $\phi$ to
identify $P\df_\cD$ and $P'_\cD$, we obtain a coaction
\[
P'_\cD\to P'_\cD\otimes p_Y^*\cA_{\qu{Y}\df}^{op}.
\]
Now Lemma~\ref{lem:quantized_comodules}(c) provides a morphism of left
$\cD_Y\otimes_\C R$-modules
\[
\phi_Y:\cD_Y\otimes_\C R\to p_Y^*\cA_{\qu{Y}\df}^{op},
\]
which is easily seen to be an isomorphism (as it equals identity
modulo $\hbar$). Finally, the uniqueness claim of
Lemma~\ref{lem:quantized_comodules}(c) implies that $\phi_Y$ is a
coalgebra homomorphism if $\cD_Y\otimes_\C R$ is equipped with the
coproduct $\delta'$.

\

\medskip

\punkt \ It remains to prove part {\bf (c)}. This can be done by
analyzing the possible choices for the morphism $\phi$. However, it is also
easy to give a direct proof.

The group 
\[
G=\{\qu{f}\in H^0(Y,\cO_\qu{Y}):\qu{f}=1\mod
\hbar\}
\] 
acts on $(\qu{Y},\qu{P})$; we want to show that the
corresponding morphism of groups $G\to\Aut(\qu{Y},\qu{P})$ is an
isomorphism. Both groups admit filtration by powers of $\hbar$, and
the claim would follow if we verify that the action map
$H^0(Y,\cO_Y)\to\End_{\cO_Z}(P)$ is an isomorphism.

By the hypotheses of Proposition~\ref{prop:obstructionsneutr}, the map
$H^0(Y,\cO_Y)\to H^0(Y,\cE(P))$ is an isomorphism. This map factors 
as 
\begin{myequation}{eq:endomorphisms}
H^0(Y,\cO_Y)\to\End_{\cO_Z}(P)\to\Hom_{\cO_Z}(P,P_{\cD/X})=H^0(Y,\cE(P))
\end{myequation}
(see \ref{sssec:EP}).
Here the map $\End_{\cO_Z}(P)\to\Hom_{\cO_Z}(P,P_{\cD/X})$ is the
composition with \eqref{eq:identityP}. Since the homomorphism
\eqref{eq:identityP} is injective, we see that all morphisms in
\eqref{eq:endomorphisms} are bijective. This completes the proof of
Proposition~\ref{prop:obstructionsneutr}.  We thus proved
Theorem~\ref{thm:obstructions} in the case when $\star$-quantization
$\gb{X}$ is neutralized. \ \hfill $\Box$
 
 \

\medskip

\punkt {\bf Remark.} \ We can also define $\gb{Y}$  
directly (without
assuming that $Y$ is a Stein manifold). For an open set $U\subset Y$,
let $W=X\times U$. Now let $\gb{Y}^{op}(U)$ be the category of triples
$(l,Q'_\cD,\upsilon)$, where $l$ is a line bundle on $U$,
$Q'_\cD\in\Comod{p_X^*\cA_X|_W}$ is flat over $R$, and 
\[
\upsilon:(P|_W\otimes p_U^* l)_\cD\to Q'_\cD/\hbar Q'_\cD
\]
is an isomorphism of comodules over $p_X^*\cD_X|_W$. A morphism of
such triples is given by a homomorphism (not necessarily invertible)
of line bundles $l$ and a compatible map of {the $p_X^*\cA_\qu{X}|_W$}-comodules
$Q'_\cD$. It is not hard to check that this gives a $\star$-quantization of $Y$.
In particular, if $U$ is Stein, we can construct local sections of $\gb{Y}^{op}(U)$
by letting $l=\cO_U$ be trivial and using Lemma~\ref{lem:quantized_comodules} to construct
$Q'_\cD$. On the other hand, there is no
guarantee that $\gb{Y}^{op}$ has any global sections.

In fact $l$ and $\upsilon$ are determined by $Q'_\cD$, so one can
alternatively consider the category of 
$Q'_\cD\in\Comod{{p_X^*\cA_\qu{X}|_W}}$ that admit such $l$ and
$\upsilon$. In other words, for a morphism between two such comodules
$Q'_\cD$, $Q\df_\cD$, there is a unique compatible map between the
corresponding line bundles $l$, $l\df$.  Indeed, the claim is local on
$U$, so we may assume that $l=l\df=\cO_U$. A morphism $Q'_\cD\to
Q\df_\cD$ induces a map between their reductions, which is an element
of
\[
\End_{p_X^*\cD_X|_W}((P|_W)_\cD)=\End_{\cD_{W/X}}((P|_W)_{\cD/X})
=\Hom_{\cO_W}(P|_W,(P|_W)_{\cD/X})=H^0(U,\cO_U),
\]
as required.
Here the last identification uses the isomorphism
$\iota:\cO_Y\to\cE(P)$.

\subsection{Proof of Theorem~\ref{thm:obstructions}: general
  case}\label{ssec:gerbydefo} 
We now drop the assumption that the $\star$-quantization $\gb{X}$ is
neutral.  Let us explain the necessary modifications to the argument
of Section~\ref{ssec:obstructionsneutr}.

We begin by extending the approach of
Section~\ref{ssec:quantcoalgebras} to general (not necessarily
neutralized) quantizations.

\

\medskip

\punkt \ Let $M$ be a complex manifold, and let $\qu{M}$ be a
neutralized $\star$-quantization of $M$. Recall that to $\qu{M}$, we
associate a coalgebra
\[
\cA_\qu{M}=\sDiff_M(\cO_\qu{M};\cO_M\otimes_\C R)
\]
in the category $\Mod{\cD_M}$, and that $\cO_\qu{M}$ is reconstructed from
$\cA_\qu{M}$ as 
\[
\cO_\qu{M} = \sHom_{\cD_M\otimes_\C
  R}({\cA_\qu{M}},\cO_M\otimes_\C R), 
\] 
see Proposition~\ref{prop:quantization_coalgebras}. 

For any $\qu{f}\in\cO_\qu{M}$, we denote the corresponding morphism by
\[
ev(\qu{f}):\cA_\qu{M}\to\cO_M\otimes_\C R.
\] 
It can be viewed as evaluation of differential operators in
$\sDiff_M(\cO_\qu{M};\cO_M\otimes_\C R)$ on $\qu{f}$. 

\

\medskip

\punkt \ Now suppose $\qu{F}\in\Comod{\cA_\qu{M}}$. We then let
$\qu{f}\in\cO_\qu{M}$ act 
on $\qu{F}$ by the composition
\[
{
\xymatrix@1@C+1.75pc@M+0.5pc{
\qu{F} \ar[r]^-{\txt{coaction}} & \cA_{\qu{M}}\otimes \qu{F}
\ar[r]^-{ev(\qu{f})} & (\cO_M\otimes_\C R)\otimes \qu{F} = \qu{F}.}
}
\]
(To simplify notation, in this section all tensor products are assumed
to be over $\cO_M\otimes_\C R$, unless explicitly stated otherwise.)
This provides $\qu{F}$ with a structure of a $\cO_\qu{M}$-module. In
other words, the action of $\cO_\qu{M}$ on $\qu{F}$ is given by the
composition
\[
\alpha_\qu{F}:\cO_\qu{M}\otimes_R\qu{F}\to
\cO_\qu{M}\otimes_R\cA_\qu{M}\otimes\qu{F}\to\qu{F}. 
\]
Note that the action $\alpha_\qu{F}$ commutes with the right action of
$\cD_M$, but not with the coaction of $\cA_\qu{M}$. Thus, 
{$\qu{F}$} 
naturally acquire a structure of an $\cO_\qu{M}\otimes_\C
\cD_M^{op}$-module, so we obtain a functor
\begin{myequation}{eq:induct}
\Comod{\cA_\qu{M}}\to\Mod{\cO_{\qu{M}}\otimes_\C\cD_M^{op}}
\end{myequation}

\

\medskip

\punkt {\bf Remark.} By Proposition~\ref{prop:absolute_comodules},
$\qu{F}$ is induced by $\qu{P}\in\Mod{\cO_\qu{M}}$, so that we can
identify
\[
\qu{F}=\qu{P}_\cD=\sDiff_M(\cO_M\otimes_\C R;\qu{P}).
\]
Under this identification, $\alpha_\qu{F}$ corresponds to the natural
action of $\cO_\qu{M}$ on differential operators $\cO_M\otimes_\C
R\to\qu{P}$. 

Equivalently, we have an isomorphism
\[
\qu{F}=\qu{P}_\cD=(\cO_\qu{M})_\cD\otimes_{\cO_\qu{M}}\qu{P},
\]
where 
\[
(\cO_\qu{M})_\cD=\sDiff_M(\cO_M\otimes_\C R;\cO_\qu{M}).
\]
Note that {$\left(\cO_\qu{M}\right)_{\cD}$} is
a bimodule over $\cO_\qu{M}$ (because 
$\cO_\qu{M}$ is a bimodule over itself); this induces the action
$\alpha_F$.

\

\medskip

\punkt {\bfseries Lemma.}\label{lm:induct}
{\it The functor \eqref{eq:induct} is fully faithful.}

\noindent
{\bf Proof.} Fix $\qu{F},\qu{G}\in\Comod{\cA_\qu{M}}$, and let 
{$\qu{\phi}:\qu{F}\to \qu{G}$} be a morphism between
the corresponding $\cO_\qu{M}\otimes_\C\cD^{op}_M$-modules. We need to
show that $\qu{\phi}$ commutes with the coaction of $\cA_\qu{M}$. In
other words, we need to verify commutativity of the diagram
\[
\xymatrix{
\qu{F}\ar[r]^{\qu{\phi}}\ar[d]&\qu{G}\ar[d]\\
\cA_\qu{M}\otimes\qu{F}\ar[r]^{\op{id}\otimes\qu{\phi}}
&\cA_\qu{M}\otimes\qu{G}.}
\]
Let $\qu{\phi}_2:\qu{F}\to\cA_\qu{M}\otimes\qu{G}$ be the difference
between the two compositions contained in this diagram. Let us check
that $\qu\phi_2=0$.

Any morphism of right $\cD_M\otimes_\C R$-modules 
$\qu\psi:\qu{F}\to\cA_\qu{M}\otimes\qu{G}$ induces a morphism
\[
\cO_\qu{M}\otimes_R\qu{F}\to\cO_\qu{M}\otimes_R
\cA_\qu{M}\otimes\qu{G}\to\qu{G},
\]
which we denote by $a(\qu\psi)$. It is easy to see that
$a(\qu\phi_2)=0$. Therefore, it suffices to check that the map
\[
a:\Hom_{\cD_M\otimes_\C R}(\qu{F},\cA_\qu{M}\otimes\qu{G})\to
\Hom_{\cD_M\otimes_\C R}(\cO_\qu{M}\otimes_R\qu{F},\qu{G})
\]
is injective. This is easy to see directly.

Indeed, $\cA_\qu{M}\otimes\qu{G}=\sDiff_M(\cO_\qu{M};\qu{G})$. Under
this identification, $a$ sends a morphism \linebreak 
$\qu{\psi}:\qu{F}\to\sDiff_M(\cO_\qu{M};\qu{G})$ to the corresponding
bilinear operator  
\[
\cO_\qu{M}\otimes_R\qu{F}\to\qu{G}:
\qu{f}\otimes\qu{s}\mapsto({\qu{\psi}}(\qu{s}))\qu{f}.  
\]
Thus, if $a(\qu\psi)=0$, then $\qu\psi=0$.
\ \hfill $\Box$

\

\medskip

\punkt {\bfseries Remark.} \label{rem-direct} The lemma is also easy to prove 
using Proposition~\ref{prop:absolute_comodules}: indeed, for
$\qu{P}\in\Mod{\cO_\qu{M}}$, an action of $\cO_\qu{M}$ on $\qu{P}_\cD$
allows us to reconstruct its action on $\qu{P}$ using the identification
$\qu{P}=\qu{P}_\cD\otimes_\cD\cO_M$. The advantage of our more
complicated argument is that it is easier to generalize to the
relative situation.

\

\medskip

\punkt {\bfseries Definition.} {\it A module
$\qu{F}\in\Mod{\cO_{\qu{M}}\otimes_\C\cD_M^{op}}$ is {\bfseries of
induced type} if it belongs to the essential image of the functor
\eqref{eq:induct}.}

\

\medskip

\punkt \ Let $\sMod{\cD_M^{op}}$ be the stack of right $\cD_M$-modules
on $M$: to an open $U\subset M$, it assigns the category
$\sMod{\cD_M^{op}}(U)=\Mod{\cD_U^{op}}$.

Let $\gb{M}$ be a $\star$-quantization of $M$ (or, more generally, a
$\star$-stack), and let $\qu{P}$ be an $\cO_\gb{M}$-module. For every
open subset $U\subset M$ and every $\alpha\in\gb{M}(U)$, we have a
neutralized $\star$-quantization $\qu{U}_\alpha$ of $U$ and an
$\cO_{\qu{U}_\alpha}$-module $\qu{P}_\alpha$. The corresponding
induced right $\cD$-module $(\qu{P}_\alpha)_\cD$ carries an action of
$\cO_{\qu{U}_\alpha}$. As $\alpha$ varies, this defines a $\C$-linear
functor
\[
\qu{P}_\cD:\gb{M}\to\sMod{\cD_M^{op}}:\alpha\mapsto(\qu{P}_\alpha)_\cD.
\]

\

\medskip

\punkt {\bfseries Definition.} \label{def:ODbimod}
{\it A $\C$-linear $1$-morphism
$\qu{F}:\gb{M}\to\sMod{\cD_M^{op}}$ is an {\bfseries
$\cO_\gb{M}-\cD_M$-bimodule}. The abelian category of
$\cO_\gb{M}-\cD_M$-bimodules is denoted by
$\Mod{\cO_\gb{M}\otimes_\C\cD^{op}_M}$.

We say that $\qu{F}\in\Mod{\cO_\gb{M}\otimes_\C\cD^{op}_M}$ is
{\bfseries of induced type} if for every open set $U\subset M$ and every
$\alpha\in\gb{M}(U)$, the $\cO_{\qu{U}_\alpha}-\cD_U$-bimodule
$\qu{F}_\alpha$ is of induced type. Let
\[
\Mod{\cO_\gb{M}\otimes_\C\cD^{op}_M}_{ind}\subset
\Mod{\cO_\gb{M}\otimes_\C\cD^{op}_M}
\]
be the full subcategory of bimodules of induced type.}

\

\medskip

\punkt {\bfseries Proposition.} {\it The functor 
\[
\qu{P}\mapsto{\qu{P}_\cD} : 
\Mod{\cO_\gb{M}}\to\Mod{\cO_\gb{M}\otimes_\C\cD^{op}_M}_{ind}
\]
is an equivalence.}

\noindent
{\bfseries Proof.} The statement is local on $M$, so we may assume
that $\gb{M}$ is neutral.  In this case, the claim reduces to
Lemma~\ref{lm:induct}.  \ \hfill $\Box$
 
\

\medskip

\punkt Let us now consider the relative version of the above
notions. Let $p_X:Z\to X$ be a submersive morphisms of complex
manifolds. A neutralized $\star$-quantization $\qu{X}$ of $X$ provides
a coalgebra $p_X^*\cA_\qu{X}\in\Mod{\cD_Z}$. A coaction of $p_X^*\cA_\qu{X}$ on
a right $\cD_Z$-module $\qu{F}$ induces an action of
$p_X^{-1}(\cO_\qu{X})$ on $\qu{F}$.  Lemma~\ref{lm:induct} remains
valid in this situation (with essentially the same proof).

\

\medskip

\punkt {\bfseries Lemma.}\label{lm:relative_induct} {\it The resulting
functor
\[
\Comod{p_X^*\cA_\qu{X}}\to\Mod{p_X^{-1}\cO_\qu{X}\otimes_\C\cD^{op}_Z}
\]
is fully faithful.}

\

\noindent
{ {\bfseries Proof.} As we explained in
  Remark~\ref{rem-direct}, one can prove Lemma~\ref{lm:induct} by
  using Proposition~\ref{prop:absolute_comodules}. The same 
  reasoning yields a proof of Lemma~\ref{lm:relative_induct} from
  Proposition~\ref{prop:relative_comodules}. \ \hfill $\Box$}

\

\medskip

\punkt {\bfseries Definition.} \label{def:ODbimodrel}
{\it Let $\gb{X}$ be a
$\star$-quantization of $X$. An {\bfseries $\cO_\gb{X}-\cD_Z$-bimodule}
is a $\C$-linear 1-morphism of stacks
\[
\qu{F}:\gb{X}\to p_{X*}\sMod{\cD_Z^{op}}.
\]
Here $p_{X*}\sMod{\cD_Z^{op}}$ is the direct image of the stack
$\sMod{\cD_Z^{op}}$ to $X$.  
Thus, for any open subset $U\subset X$,
\[
p_{X*}\sMod{\cD^{op}_Z}(U)=\sMod{\cD^{op}_Z}(p_X^{-1}(U))
=\Mod{\cD^{op}_{p_X^{-1}(U)}}.
\]
For any open set $U\subset X$ and any $\alpha\in\gb{X}$, the right
$\cD_{p_X^{-1}(U)}$-module $\qu{F}_\alpha$ carries a natural left
action of $p_X^{-1}(\cO_{\qu{U}_\alpha})$. We say that $\qu{F}$ is
{\bfseries of induced type} if for all $U$ and $\alpha$, the module
\[
\qu{F}_\alpha\in\Mod{p_X^{-1}(\cO_{\qu{U}_\alpha})\otimes_\C
  \cD^{op}_{p_X^{-1}(U)}}
\]
is induced from a $p_X^*\cA_{\qu{U}_\alpha}$-comodule; that is, it
belongs to the essential 
image of the functor of Lemma~\ref{lm:relative_induct}. 
{If $V
\subset Z$ is an open subset and $W = p_{X}(V) \subset X$ we will write
\[
\cC_{V} =\Mod{p_X^{-1}\cO_\gb{W}\otimes_\C \cD^{op}_V}_{ind}
\]
for the category of $\cO_\gb{W}-\cD_V$-bimodules of induced type.}} 

\

\medskip

\punkt {\bfseries Example.}\label{exm:neutral_induct} Suppose $\gb{X}$
admits a neutralization $\alpha\in\gb{X}(X)$. Then $\qu{F}\in{\mathcal
C}_{Z}$ is uniquely determined by 
\[
\qu{F}_\alpha\in \Mod{p_X^{-1}\cO_{\qu{X}_\alpha}\otimes_\C
  \cD^{op}_Z}.
\]
Moreover, $\qu{F}_\alpha$ must belong to the essential image of the
induction functor of Lemma~\ref{lm:relative_induct}, so that
$\qu{F}_\alpha$ is induced by a
$p_X^*(\cA_{\qu{X}_\alpha})$-comodules. This defines an equivalence
between $\cC$ and $\Comod{\cA_{\qu{X}_\alpha}}$. See Corollary~\ref{cor:trivialize}
for a more general statement.

\

\medskip

\punkt {\bfseries Proposition.} \label{prop-abelianC} 
{\it {$\cC_{Z}$} is an abelian category with enough
  injective objects.}

\noindent
{\bfseries Proof.} The
direct image functor $\cC_V\to{\cC_{Z}}$ preserves
injectivity, because its 
left adjoint is exact. Therefore, the statement of the proposition is
local on $Z$.

Note that the category $\cC_V$ depends only on the restriction of the
$\star$-quantization $\gb{X}$ to the open subset 
{$W = p_X(V)\subset X$}. By choosing $V$ small
enough, we may assume that this restriction is neutral. We are thus in
the {setting} of Example~\ref{exm:neutral_induct}, and
Lemma~\ref{lm:comod_inj} implies the statement. \ \hfill $\Box$

\

\medskip

\punkt Finally, note that {$\cC_{Z}$ is
  equipped with} an action of the tensor category
of $\cD_Z\otimes_\C R$-modules: given
$\qu{F}\in{\cC_{Z}}$ and
$G\in\Mod{\cD_Z\otimes_\C R}$, we define $G\otimes \qu{F}$ by
\[
(G\otimes\qu{F})_\alpha =
G\otimes(\qu{F}_\alpha)\qquad(\alpha\in\gb{X}(U),U\subset X).
\]
The tensor product is over $\cO_Z\otimes_\C R$.

\

\medskip

\punkt \ We can now repeat the argument of
Section~\ref{ssec:obstructionsneutr} using the category 
{$\cC_{Z}$} in place of $\Comod{p_X^*\cA_X}$. The key
step is the following observation.

Let $Z=X\times Y$ be a product of two complex manifolds. Suppose
$\gb{X}$ is a $\star$-quantization of $X$ and $\qu{Y}$ is a
neutralized $\star$-quantization of $Y$. The inverse image
$p_Y^*\cA_\qu{Y}$ is a coalgebra in the category $\Mod{\cD_Z\otimes_\C
  R}$. Since this category acts on {$\cC_{Z}$}, it makes sense to consider
$p_Y^*\cA_\qu{Y}$-modules in this category.  In particular, for any
$\cO_{\gb{X}\times\qu{Y}}$-module $\qu{P}$, the induced bimodule
$\qu{P}_\cD\in{\cC_{Z}}$ carries a coaction
of $p_Y^*\cA_\qu{Y}$.
 
\

\medskip

\punkt {\bfseries Proposition.} {\it The correspondence
$\qu{P}\mapsto\qu{P}_\cD$ is an equivalence between the category of
$\cO_{\gb{X}\times\qu{Y}}$-modules and the category of
$p_Y^*\cA_\qu{Y}$-comodules in {$\cC_{Z}$}.}

\noindent
{\bfseries Proof.} The statement is local on $Z$, so we may assume
without losing generality that $\gb{X}$ corresponds to a neutralized
$\star$-quantization $\qu{X}$ of $X$. Then 
{by Lemma~\ref{lm:induct}  we have that 
$\cC_{Z}=\Comod{p_X^*\cA_\qu{X}}$}, and $p_Y^*\cA_\qu{Y}$-comodules in
{$\cC_{Z}$} are simply comodules over
\[
p_X^*\cA_\qu{X}\otimes p_Y^*\cA_\qu{Y}=\cA_{\qu{X}\times\qu{Y}}.
\]
The proof now follows from
Proposition~\ref{prop:absolute_comodules}. \ \hfill $\Box$

\

\medskip

\punkt \ Let us now outline how the argument of
Section~\ref{ssec:obstructionsneutr} can be adapted to prove
Theorem~\ref{thm:obstructions} in full generality.  We are given a
coherent sheaf $P$ on $Z=X\times Y$ and a 
$\star$-quantization $\gb{X}$ of $X$. As before, we may assume that
$Y$ is Stein.  One can check that Lemma~\ref{lem:quantized_comodules}
remains true if we replace the category $\Comod{p_X^*\cA_\qu{X}}$
with $\cC_Z$; this requires an analogue of Proposition~\ref{prop:defcomodules}
for the category $\cC_Z$, which in turn relies on Proposition~\ref{prop-abelianC}.
Such version of Lemma~\ref{lem:quantized_comodules} provides a deformation of
$P_\cD\in\Comod{p_X^*\cD_X}$ to 
\ $P'_\cD\in{\cC_{Z}}$. Next note that
  $P'_{\cD}$  can
be equipped with a homomorphism
\[
\bDelta'_Y:P'_\cD\to(p_Y^*\cD_Y\otimes_\C R)\otimes P'_\cD,
\]
which satisfies the universal property of
Lemma~\ref{lem:quantized_comodules}(c).  The universal property
then implies that there is a unique coproduct $\delta'$ on
$\cD_Y\otimes_\C R$ such that $\bDelta'_Y$ is a coaction. Define the
neutralized $\star$-quantization $\qu{Y}$ of $Y$ in such a way that
the coalgebra $\cA_\qu{Y}^{op}$ is identified with $\cD_Y\otimes_\C R$
equipped with the coproduct $\delta'$. Then $P'_\cD\in\cC$ is a
comodule over $p_Y^*\cA_\qu{Y}^{op}$, as required.

This proves the existence statement of Theorem~\ref{thm:obstructions};
the proof of uniqueness is completely parallel to the proof of parts
{\bfseries (b)} and {\bfseries (c)} of
Proposition~\ref{prop:obstructionsneutr}. \ \hfill $\Box$

\subsection{Comments on tangible sheaves and
  $\star$-deformations}\label{sec:tangible} 
In our approach to Theorem~\ref{thm:obstructions}, we study
deformations of $\cO$-modules by deforming induced
$\cD$-modules. Under some additional assumptions (which hold, for
instance, in the case of the Fourier-Mukai transform on complex tori),
one can interpret deformations of induced $\cD$-modules in a more
explicit way by looking at `$\star$-deformations of an $\cO$-module'.

Let $M$ be a complex manifold.

\

\medskip

\punkt {\bf Definition.} \label{def:tangible} {\it A coherent
$\cO_M$-module $P$ is {\bfseries tangible} if $P\simeq i_{\Gamma*}V$,
where $i_{\Gamma} : \Gamma\hookrightarrow M$ is a closed analytic
submanifold, and $V$ is a holomorphic vector bundle on $\Gamma$.  }

\

\medskip

\punkt \ Recall that for an $\cO_M$-module $P$, we denote by
$P_\cD=\sDiff_M(\cO_M;P)=P\otimes_{\cO_M}\cD_M$ the induced right
$\cD_M$-module.  From Lemma~\ref{lem:diffops_via_ind}, we obtain an
identification
\[
\sDiff_M(P;P)=\sHom_{\cO_M}(P,P_\cD)=\sHom_{\cD_M}(P_\cD,P_\cD).
\]
If we assume that $P$ is tangible, the identification is valid in the
derived sense; in other words, there are no `higher derived'
differential operators from $P$ to itself.

\medskip

\

\punkt {\bf Lemma.} \label{lem:sDiff} {\it Suppose $P$ is tangible. Then
\[
\sExt^i_{\cO_M}(P,P_\cD)=\sExt^i_{\cD_M}(P_\cD,P_\cD)=0\qquad(i>0).
\]}

\noindent
{\bf Proof.} Suppose $P=i_{\Gamma *}V$, where $i_{\Gamma} : \Gamma
\hookrightarrow M$ is a closed analytic submanifold, and $V$ is a
holomorphic vector bundle on $\Gamma$. By adjunction, we have
\[
\sExt^i_{\cO_M}(P,P_\cD)=i_{\Gamma*}\sHom(V,\RR^i i_{\Gamma}^!P_\cD).
\]
Recall that 
\[
\RR^i i_{\Gamma}^!(\bullet)=\LL_{(\codim\Gamma-i)}\;
i_{\Gamma}^*(\bullet)\otimes\omega_M^{-1}\otimes\omega_\Gamma,
\]
so we need to verify that
\begin{myequation}{eq:kashiwara}
\LL_k\; i_\Gamma^*(P_\cD)=0\quad(k<\codim\Gamma).
\end{myequation}

It is not hard to verify \eqref{eq:kashiwara} directly. However, it
also immediately follows from the Kashiwara Lemma {\cite[Theorem~4.30]{kashiwara-dmodules}},  because
$P_\cD$ is a (right) $\cD_M$-module supported by $\Gamma$.  \ \hfill
$\Box$

\medskip

\

\punkt {\bfseries Definition.} \label{defn:sheafstar}
{\it A {\bfseries $\star$-deformation} of $P$ over $R$ is a pair 
${\frP :=} \left(\qu{P},\mycal{I}_{\qu{P}}\right)$
where
\begin{itemize}
\item $\qu{P}$ is a sheaf of $R$-modules on $M$.
\item $\mycal{I}_{\qu{P}}$ (a {\bfseries $\star$-structure} on
  $\qu{P}$) is a subsheaf in the sheaf of $R$-module isomorphisms
  \linebreak $\sIso_{R}(P\otimes_{\mathbb{C}} R,\widetilde{P})$ that
  is a torsor over the sheaf of groups
\[
\sDiff_{0}\left(P\otimes_{\mathbb{C}} R;P\otimes_{\mathbb{C}} R\right)
:= \left\{\left. \sum_{i=0}^n\displaylimits D_{i} \hbar^{i} \ \right|
\ D_{0} = 1, \ D_{i} \in \sDiff_{M}(P;P) \right\}.
\]
\end{itemize}

\noindent
{The  {\bfseries trivial $\star$-deformation}  
$\frP^{\op{triv}}$ of $P$ is defined as the pair
\[
\frP^{\op{triv}} := \left(P\otimes_{\C} R,
  \sDiff_{0}(P\otimes_{\C} R, P\otimes_{\C} R)\right).
\]}
}
\

\medskip

\punkt {\bf Remarks.} \ \label{rem:1} {\bf (1)} \ Usually for a
  deformation quantization one requires flatness over $R$. In
  Definition~\ref{defn:sheafstar}, $\qu{P}$ is automatically
  $R$-flat. Indeed, $\mycal{I}_{\qu{P}}$ has local sections, so that
  $\qu{P}$ is locally isomorphic to $P\otimes_{\mathbb{C}} R$.

\noindent
{\bf (2)} \ Since $\sDiff_{\!0}(P\otimes_{\mathbb{C}} R;
P\otimes_{\mathbb{C}} R)$ consists of operators that are normalized to
be $1$ modulo $\hbar$, it follows that 
  isomorphisms in $\mycal{I}_{\qu{P}}$ induce a
canonical isomorphism $P \widetilde{\to}
\qu{P}/\hbar{\qu{P}}$.

\

\medskip
  
\noindent
Next we need the notion of $\star$-local maps between
$\star$-deformations.
\

\smallskip

\punkt {\bf Definition.} \label{def:starmap} {\it Let 
  ${\frP_{1} = }(\qu{P}_{1},\mycal{I}_{\qu{P}_{1}})$,  \ldots,
  ${\frP_{k} = }(\qu{P}_{k},\mycal{I}_{\qu{P}_{k}})$,
  ${\frP = }(\qu{P},\mycal{I}_{\qu{P}})$ be
  $\star$-deformations of some 
  coherent
  analytic sheaves $P_{1}$, \ldots,
  $P_{k}$, $P$ on $M$. We say that an $R$-linear sheaf map
\[
\tilde{f} : \qu{P}_{1}\otimes_{R} \qu{P}_{2} \otimes_{R} \ldots
\otimes_{R} \qu{P}_{k} \to \qu{P}.
\]
is {\bfseries $\star$-local} if for every
choice of $\tilde{a}_{i} \in \mycal{I}_{\qu{P}_{i}}$, $i = 1, \ldots,
n$, and $\tilde{a} \in
\mycal{I}_{\qu{P}}$, the map
\[
\tilde{a}^{-1}\circ \tilde{f} \circ \left( \tilde{a}_{1}\otimes \ldots
\otimes \tilde{a}_{k} \right) :
P_{1}\otimes_{\mathbb{C}} \ldots \otimes_{\mathbb{C}}
  P_{n}\otimes_{\mathbb{C}} R \to P\otimes_{\mathbb{C}} R
\]
is a polydifferential operator, that is, it belongs to
$\sDiff_M(P_{1}, \ldots, P_{k}; P)\otimes_{\mathbb{C}} R$. The sheaf
of $\star$-local maps $\tilde{f}$ is denoted by 
{$\sDiff_M(\frP_1,\dots,\frP_k;\frP)$}.}

\

\medskip

\punkt \ Let ${\frP = }
\left(\qu{P},\mycal{I}_{\qu{P}}\right)$ be a 
$\star$-deformation of a tangible $\cO_M$-module $P$. Let
{$\frP_\cD$} 
be the sheaf of $\star$-local maps
\[
{
\sDiff_M(\mathfrak{O}^{\op{triv}}_M;\frP),
}
\]
{where $\mathfrak{O}^{\op{triv}}_M$ denotes the trivial
  $\star$-deformation of $\cO_{M}$.}  Since $\cD_M$ acts on
$\cO_M\otimes_\C R$ on the left, we obtain a right action of
$\cD_M\otimes_\C R$ on {$\frP_\cD$}. It is
clear that {$\frP_\cD$} is $R$-flat and that
{$P_\cD=\frP_\cD/\hbar\frP_\cD$}, so that
{$\frP_\cD$}  is
an $R$-deformation of the right $\cD_M$-module $P_\cD$.

\

\medskip

\punkt {\bf Proposition.} \label{prop:star-induced} 
{\it For any tangible $\cO_M$-module $P$, the correspondence 
\[
{\frP = } 
\left(\qu{P},\mycal{I}_{\qu{P}}\right)\mapsto {\frP}_\cD
\]
provides an equivalence between the category of $\star$-deformations
of $P$ and that of deformations of the right $\cD_M$-module $P_\cD$.
}

\noindent
{\bf Proof.} By definition, any $\star$-deformation 
{$\frP$} of $P$ is locally trivial, that is,
{$\frP$} is locally isomorphic to
{$\frP^{\op{triv}}$} with the obvious
$\star$-structure. Therefore, $\star$-deformations of $P$ are in
one-to-one correspondence with torsors over
$\sDiff_{0}\left(P\otimes_\C R;P\otimes_\C R\right)$: the
correspondence associates $\mycal{I}_\qu{P}$ to
$\frP = \left(\qu{P},\mycal{I}_{\qu{P}}\right)$. On the other hand, for any
Stein open set $U\subset M$,
\[
\Ext^1_{\cD_U}((P_\cD)|_U,(P_\cD)|_U)=0
\]
by Lemma~\ref{lem:sDiff}. Therefore, any deformation of
$P_\cD\in\Mod{\cD_M^{op}}$ is locally trivial.  Again, we see that
deformations of $P_\cD$ are in one-to-one correspondence with torsors
over the same sheaf. This implies the statement.

For the sake of completeness, let us describe the inverse
correspondence. Let $P'_\cD$ be an $R$-deformation of $P_\cD$. Set
$\qu{P}=P'_\cD\otimes_{\cD_M}\cO_M$.  It is easy to see that $\qu{P}$
is $R$-flat and that $\qu{P}/\hbar\qu{P}=P$, because
$\sTor_i^{\cD_M}(P_\cD,\cO_M)=0$ for $i>0$ by
Lemma~\ref{lem:induced_adjunct}.  Note that $\qu{P}$ is equipped with
a natural evaluation map $P'_\cD\to\qu{P}$.

Since $\sExt^i_{\cO_M}(P,P_\cD)=0$ for $i>0$ by Lemma~\ref{lem:sDiff},
the sheaf
\[
\mycal{H}_{\qu{P}}=\sHom_{\cO_M}(P,P'_\cD)=\sHom_{\cD_M}(P_\cD,P'_\cD)
\] 
is $R$-flat and 
\[
\mycal{H}_{\qu{P}}/\hbar\mycal{H}_{\qu{P}} =
\sHom_{\cO_M}(P,P_\cD)=\sDiff_M(P;P).
\]
Let $\mycal{I}_{\qu{P}}\subset\mycal{H}_{\qu{P}}$ be the preimage of
$1\in\sDiff_M(P;P)$.

The evaluation map $P'_\cD\to\qu{P}$ allows us to embed
$\mycal{H}_{\qu{P}}$ in
\[
\sHom_{\mathbb{C}}(P,\qu{P})=\sHom_R(P\otimes_{\mathbb{C}}R,\qu{P}).
\] 
In this way, we view $\mycal{I}_{\qu{P}}$ as a subsheaf of
$\sIso_{R}(P\otimes_{\mathbb{C}}
R,\widetilde{P})\subset\sHom_R(P\otimes_{\mathbb{C}}R,\qu{P})$. It is
automatic that ${\frP :=
}\left(\qu{P},\mycal{I}_{\qu{P}}\right)$ is a $\star$-deformation of
$P$.  

{Suppose that $\frP^{(1)}$
and $\frP$ are $\star$-deformations of two tangible sheaves $P^{(1)}$
and $P$ on $M$. A $\star$-local map $\frP^{(1)}\to\frP$ induces a
homomorphism of $\cD_M$-modules $\frP^{(1)}_\cD\to\frP_\cD$, which
is simply the composition of $\star$-local maps
\[
\cO_M\otimes_{\mathbb{C}}R\to\qu{P}^{(1)}\to\qu{P}.
\]
This provides a one-to-one correspondence between $\star$-local maps
$\frP^{(1)}\to\frP$ and homomorphisms of $\cD_M$-modules
$\frP^{(1)}_\cD\to\frP_\cD$ and completes the proof of the proposition.}
\ \hfill $\Box$\

\

\medskip

\punkt {\bf Remark.} \label{rem:starmaps} More generally, suppose
{$\frP^{(1)},\ldots, \frP^{(k)},\frP$} are
$\star$-deformations of tangible sheaves $P^{(1)},\ldots, P^{(k)},P$
on $M$. {The argument at the end of the proof
  of the previous proposition identifies the sheaf
  $\sDiff_{M}(\frP^{(1)},\ldots,\frP^{(k)};\frP)$ of all
  $\star$-local maps with the sheaf of all $k$-ary $\star$-operations
  from $\{\frP^{(1)}_\cD,\dots,\frP^{(k)}_\cD\}$ to $\frP_\cD$}.  (The
definition of $\star$-operations is recalled in
Remark~\ref{rem:star_operations}.)

\

\medskip

\punkt \ Let $\qu{M}$ be a neutralized $\star$-quantization of $M$,
and let $\qu{P}$ be an $\cO_\qu{M}$-module that is flat over $R$ and
such that $P:=\qu{P}/\hbar\qu{P}$ is a tangible
$\cO_M$-module. 
{Proposition~\ref{prop:star-induced} and
  Proposition~\ref{prop:absolute_comodules} imply that $\qu{P}$
  carries a natural $\star$-structure. Indeed, consider
the sheaf of differential operators $\sDiff_M(P\otimes_\C R;\qu{P})$.
It can be identified with
\[
\sHom_{\cO_M}(P,\qu{P}_\cD),
\]
so Lemma~\ref{lem:sDiff} implies that $\sDiff_M(P\otimes_\C R;\qu{P})$
is an $R$-flat deformation of $\sDiff_M(P;P)$. Therefore,
\[
\mycal{I}_{\qu{P}}:=\{A\in\sDiff(P\otimes_\C R;\qu{P}):A=1\mod\hbar\}
\]
provides a natural $\star$-structure on $\qu{P}$.

Note that if $\frP =
  (\qu{P},\mycal{I}_{\qu{P}})$, then
  Proposition~\ref{prop:star-induced} implies that $\frP_{\cD} =
  \qu{P}_{\cD}$ as right $\cD_{M}$-modules. From now on, we will not
  make a notational distinction between $\frP_{\cD}$ and
  $\qu{P}_{\cD}$.  
  
  Note also that the action map
$\cO_\qu{M}\times\qu{P}\to\qu{P}$ is $\star$-local once we equip
  $\cO_\qu{M}$ and $\qu{P}$ with 
  $\star$-structures $\mycal{I}_{\cO_\qu{M}}$ and
  $\mycal{I}_\qu{P}$.} In addition, we have the
  following

\

\medskip

\punkt {\bf Lemma.} \label{lem:uniquestar} {\it Let $\qu{M}$ be a
neutralized $\star$-quantization of $M$.
\begin{enumerate}
\item[{\bf(a)}] $\mycal{I}_{\cO_\qu{M}}$ is the only $\star$-structure
on $\cO_\qu{M}$ that makes the multiplication map
\[
\cO_\qu{M}\times\cO_\qu{M}\to\cO_\qu{M}
\] 
$\star$-local. This
provides an equivalence between the category of neutralized
$\star$-quantizations of $M$ and the category of $\star$-deformations
$(\cO_\qu{M},\mycal{I}_{\cO_\qu{M}})$ of $\cO_M$ together with a lift
of the algebra structure on $\cO_M$ to an algebra structure on
$\cO_\qu{M}$ that is $\star$-local.

\item[{\bf(b)}] Let $\qu{P}$ be an $\cO_\qu{M}$-module that is flat
over $R$ and such that $P:=\qu{P}/\hbar\qu{P}$ is a tangible
$\cO_M$-module. Then $\mycal{I}_\qu{P}$ is the only $\star$-structure
on $\qu{P}$ that makes the action map
\[
\cO_\qu{M}\times\qu{P}\to\qu{P}
\] 
$\star$-local.
\end{enumerate}
}

\noindent
{\bf Proof.} Both statements of the lemma are easily proved directly
(for {\bf (a)}, see \cite[Proposition~2.2.3]{KSarxiv}).  Alternatively
they can be reduced to our previous results once we properly interpret
them by using Remark~\ref{rem:starmaps}. For part {\bf(a)},
Proposition~\ref {prop:star-induced} identifies the
$\star$-deformations $(\cO_\qu{M},\mycal{I}_{\cO_\qu{M}})$ with
deformations $\widetilde{\cD_M}$ of $\cD_M$ in the category of right
$\cD_M$-modules. By Remark~\ref{rem:starmaps}, a $\star$-local algebra
structure on $\cO_\qu{M}$ is interpreted as a $\star$-operation on
$\widetilde{\cD_M}$. Dualizing (as in
Remark~\ref{rem:star_operations}), we obtain a coalgebra structure on
the right $\cD_M$-module
\[
\widetilde{\cA_{M}}=\sHom_{\cD_M}(\widetilde{\cD_M},\cD_M).
\]
We now see that part {\bf (a)} is equivalent to
Proposition~\ref{prop:quantization_coalgebras}. In the same way, part
{\bf (b)} is equivalent to
Proposition~\ref{prop:absolute_comodules}(a).\ \hfill $\Box$

\

\medskip

\punkt \label{ssec:tangiblerel} Consider now the relative situation. Let $p_X:Z\to X$ be a
submersive morphism of complex manifolds. Let $(\qu{P},\mycal{I}_\qu{P})$ be a
$\star$-deformation of a tangible coherent sheaf $P$ on $Z$, and let $(\qu{Q},\mycal{I}_\qu{Q})$
be a $\star$-deformation of a tangible coherent sheaf $Q$ on $X$. Given an
$R$-linear map $\qu{f}:p_X^{-1}(\qu{Q})\to \qu{P}$, we can state the
property of $\qu{f}$ being $\star$-local completely analogously to
Definition~\ref{def:starmap}.

The notion of $\star$-locality also makes sense for polylinear
maps. In the case that we are particularly interested in we consider a
quantization $\qu{X}$ of $X$. In this case, the sheaf $\cO_\qu{X}$ carries
a natural $\star$-structure that turns it into a $\star$-deformation of $\cO_X$.
Now suppose that we are given an action
\[
\qu{\xi}:p_X^{-1}\cO_\qu{X}\otimes_R\qu{P}\to\qu{P}.
\]
We then say that $\qu{\xi}$ is {\bfseries $\star$-local} if for every
choice of $\tilde{a} \in \mycal{I}_{\qu{P}}$,
$\tilde{b}\in\mycal{I}_{\cO_{\qu{X}}}$, the map
\[
\qu{a}^{-1}\circ\qu{\xi}\circ\left( p_X^{-1}(\qu{b})\otimes
\qu{a}\right): p_X^{-1}(\cO_X)\otimes_{\mathbb{C}} P
\otimes_{\mathbb{C}} R \to P\otimes_{\mathbb{C}} R
\]
is a polydifferential operator, that is, if it belongs to
$\sDiff_Z(p_X^{-1}\cO_X, P; P)\otimes_{\mathbb{C}} R$.  

\

\medskip

\punkt \ Recall that the inverse image $p_X^*\cA_{\qu{X}}$ is a
coalgebra in the category of left $\cD_Z\otimes_\C R$-modules. Also
recall that to $\qu{P}$, we assign a right
$\cD_Z\otimes_\C R$-module $\qu{P}_\cD=\sDiff_Z(\cO_Z\otimes_\C
R;\qu{P})$. It is easy to see that a $\star$-local action $\qu{\xi}$
yields a coaction
\[
\bDelta_{\qu{\xi}}:\qu{P}_\cD\to p_X^*\cA_{\qu{X}}\otimes\qu{P}_\cD,
\]
where the tensor product is over $\cO_Z\otimes_\C R$. Here
$\bDelta_\qu{\xi}$ is given essentially by the same formulas as those
in \ref{ssec:coaction}. Namely, let us identify
$p_X^*\cA_{\qu{X}}\otimes\qu{P}_\cD$ with the sheaf of $\star$-local
maps
\[
p_X^{-1}\cO_{\qu{X}}\otimes_\C \cO_Z\to\qu{P}.
\]
Then $\bDelta_{\qu{\xi}}$ sends a differential operator
\[
\qu{A}\in\sDiff_Z(\cO_Z\otimes_\C R;\qu{P})=\qu{P}_\cD
\] to the
$\star$-local map
\[
\bDelta_{\qu{\xi}}(\qu{A}):\qu{f}\otimes g\mapsto
\qu{\xi}(\qu{f}\otimes\qu{A}(g)).
\]
We now have the following straightforward

\

\medskip

\punkt {\bf Lemma.}\label{lem:star-inducedrel} {\it The correspondence
$\xi\mapsto\bDelta_{\qu\xi}$ is a bijection between $\star$-local
actions $\qu{\xi}$ and coactions of $p_X^*\cA_{\qu{X}}$ on
$\qu{P}_\cD$.} 

\

\noindent
{{\bfseries Proof.}  This follows from
  Lemma~\ref{lem:uniquestar}.  Indeed, let $P \in
  \Mod{\cO_{M}}$ be a tangible sheaf on a manifold $M$. Then by
  Proposition~\ref{prop:star-induced} we have an 
  equivalence of categories}
\[
{
\left\{   
\text{
\begin{minipage}[c]{2.2in} 
$\star$-deformations  $\frP =
    \left(\qu{P},\mycal{I}_{\qu{P}}\right)$   of $P$ over $R$
\end{minipage}}\right\} 
\cong 
\left\{ 
\text{
\begin{minipage}[c]{2.2in} 
deformations $\frP_{\cD}$ of $P_{\cD}$ as a right $\cD_{M}$-module
\end{minipage}
}
\right\}
}
\]
{
Furthermore, when we  consider $P = \cO_{M}$ and take into account the
(co)algebra structure we get equivalences of categories}
\[
\begin{split}
{
\left\{
\text{
\begin{minipage}[c]{2.2in}
$\star$-deformations $\left(\qu{\cO_{M}},\mycal{I}_{\qu{\cO_{M}}}\right)$
  of $\cO_{M}$ over $R$ with a lift of the $\cO_{M}$-algebra structure
  as a \linebreak $\star$-local algebra structure
\end{minipage}
}
\right\} 
} & \quad
{
\stackrel{\text{Lemma~\ref{lem:uniquestar}(a)}}{\cong}  \quad \left\{
\text{
\begin{minipage}[c]{2.1in}
$R$-deformations $\qu{\cA_{M}}$ of $\cD_{M}$ as a coalgebra in
  $\Mod{\cD_{M}}$ 
\end{minipage} 
}
\right\}
} \\[+1pc]
& \quad 
{
\stackrel{Proposition~\ref{prop:quantization_coalgebras}}{\cong} \quad 
\left\{
\text{
\begin{minipage}[c]{1.7in}
$\star$-quantizations $\qu{M}$ of $M$
\end{minipage}
}
\right\}
}
\end{split}
\]
{
These equivalences justify the slight abuse of notation 
$\left(\qu{\cO_{M}},\mycal{I}_{\qu{\cO_{M}}}\right) =
\left(\cO_{\qu{M}},\mycal{I}_{\cO_{\qu{M}}}\right)$  and $\qu{\cA_{M}}
= \cA_{\qu{M}}$. Moreover, the coalgebra structure allows us to consider
comodules over $\cA_{\qu{M}}$. Thus we can extend the previous
equivalences to:
}
\[
\begin{split}
{
\left\{
\text{
\begin{minipage}[c]{2.2in}
$\star$-deformations $\frP = \left(\qu{P},\mycal{I}_{\qu{P}}\right)$ of
  $P$ over $R$ with a lift of the $\cO_{M}$-action to a $\star$-local action of
  $\cO_{\qu{M}}$
\end{minipage}
}
\right\}
} & 
{
\quad \stackrel{\text{Lemma~\ref{lem:uniquestar}(b)}}{\cong} \quad 
\left\{
\text{
\begin{minipage}[c]{1.7in}
$R$-deformations of $P_\cD$ to a comodule
$\frP_{\cD}$ over $\qu{\cA_{M}}$
\end{minipage}
}
\right\}
} \\[+1pc]
& 
{
\quad
\stackrel{\text{Proposition~\ref{prop:absolute_comodules}(b)}}{\cong}
\quad  
\left\{
\text{
\begin{minipage}[c]{1.4in}$R$-deformations of $P$ to a module 
$\qu{P}$ over $\cO_{\qu{M}}$
\end{minipage}
}
\right\}.
}
\end{split}
\]
\ \hfill $\Box$}

\

\medskip

\punkt \ Combining Proposition~\ref{prop:star-induced} and
Lemma~\ref{lem:star-inducedrel}, we can now make the argument of
Section~\ref{ssec:obstructionsneutr} more explicit, at least when $P$
is a tangible coherent sheaf on $Z$.  Let us sketch this reformulation
using the notation of Section~\ref{ssec:obstructionsneutr}.  We assume
that $P$ is tangible.

The first step is to construct a $\star$-deformation 
{$\frP = \left(\qu{P},\mycal{I}_{\qu{P}}\right)$} of
$P$ together with an extension of the action
\[
\xi:p_X^{-1}(\cO_X)\otimes_\C P\to P
\]
to a $\star$-local action
\[
\qu{\xi}:p_X^{-1}(\cO_\qu{X})\otimes_R\qu{P}\to \qu{P}.
\] 
Lemma~\ref{lem:quantized_comodules}(a) claims that this can be done in
a way that is unique up to a non-canonical isomorphism.

The second step is to extend the action 
\[
\eta:p_Y^{-1}(\cO_Y)\otimes_\C P\to P
\]
to a $\star$-local map (not necessarily an action)
\[
\qu{\eta}:p_Y^{-1}(\cO_Y)\otimes_\C \qu{P}\to\qu{P}.
\]
Such an extension is provided by
Lemma~\ref{lem:quantized_comodules}(b).

Finally, by the universal property from
Lemma~\ref{lem:quantized_comodules}(c) it follows that there is a
unique local product on $p_Y^{-1}(\cO_Y)\otimes_\C R$ such that
$\qu{\eta}$ is an action of the resulting algebra. Viewing
$p_Y^{-1}(\cO_Y)\otimes_\C R$ with this product, we obtain a
neutralized $\star$-quantization $\qu{Y}^{op}$ of $Y$ such that
$\qu{P}$ has a structure of a $\cO_{\qu{X}\times\qu{Y}^{op}}$-module.

\subsection{Another proof of Theorem~\ref{thm:obstructions}}
\label{ssec:obstructionsneutr'} 

The proof of Theorem~\ref{thm:obstructions} can be restated using
(relative) differential operators on quantizations, thus avoiding
coalgebras in the category of comodules. The resulting argument is
somewhat more elementary, but less transparent. The approach is also
not completely canonical: essentially, one can work with relative
differential operators for the morphism $\qu{Z}'\to\qu{X}$, where
$\qu{Z}'$ is any neutralized $\star$-quantization of $Z=X\times Y$
that admits a map to the given quantization $\qu{X}$ of $X$. (The
category of $\cD_{\qu{Z}'/\qu{X}}$-modules does not depend on the
choice of $\qu{Z}'$.) We take $\qu{Z}'=\qu{X}\times Y$; more
precisely, $\qu{Z}'$ is the product of $\qu{X}$ and the trivial
neutralized $\star$-quantization of $Y$.

\

\medskip

\punkt \ As in Section~\ref{ssec:obstructionsneutr}, we assume that
the $\star$-quantization of $X$ is neutralized, and our goal is to
prove Proposition~\ref{prop:obstructionsneutr}. One can then extend
the proof to all $\star$-quantizations of $X$ using the approach of
Section~\ref{ssec:gerbydefo}, see
Section~\ref{ssec:sectionsoverquant}. This is not done in this paper.

\punkt \ Set $\qu{Z}'=\qu{X}\times Y$, and let $\cD_{\qu{Z}'/\qu{X}}$
be the sheaf of differential operators $\cO_{\qu{Z}'}\to\cO_{\qu{Z}'}$
that commute with the left action of $p_X^{-1}(\cO_\qu{X})$. Recall
that $\cD_{\qu{Z}'/\qu{X}}$ is a flat sheaf of $R$-algebras by
Proposition~\ref{prop:relqudiffops_flat} and that
\[
\cD_{Z/X}=\cD_{\qu{Z}'/\qu{X}}/\hbar\cD_{\qu{Z}'/\qu{X}}.
\]
Let $P_{\cD/X}$ be the right $\cD_{Z/X}$-module induced from $P$. We
derive Proposition~\ref{prop:obstructionsneutr} from the following
claim (which is essentially a version of
Lemma~\ref{lem:quantized_comodules}).

\

\medskip

\punkt {\bf Lemma.} \label{lem:obstructionsneutr} {\it In the
  assumptions {and notation} of
  Proposition~\ref{prop:obstructionsneutr}, we have the following.
\begin{enumerate}
\item[{\bf (a)}] There is a $R$-flat deformation of $P_{\cD/X}$ to a
$\cD_{\qu{Z}'/\qu{X}}$-module $P'_{\cD/X}$, unique up to a not
necessarily unique isomorphism.
\item[{\bf (b)}] The direct image
$\cE':=p_{Y*}\sEnd_{\cD_{\qu{Z}'/\qu{X}}}(P'_{\cD/X})$ is flat over
$R$. The action of $p_Y^{-1}\cO_Y$ on $P_{\cD/X}$ (coming from its
action on $P$) induces isomorphisms
\[
\cO_Y\to p_{Y*}\sEnd_{\cD_{Z/X}}(P_{\cD/X})\leftarrow \cE'/\hbar\cE'.
\]
\item[{\bf (c)}] Let 
\[m:p^{-1}_Y(\cO_Y)\times P_{\cD/X}\to P_{\cD/X}\]
be the ($\C$-bilinear) morphism induced by the action of
$p^{-1}_Y(\cO_Y)\subset\cO_Z$ on $P$. 
Then $m$ can be lifted to a $\C$-bilinear morphism
\[
m':p^{-1}_Y(\cO_Y)\times P'_{\cD/X}\to P'_{\cD/X}
\]
such that $m'(\qu{f},A)$ is a morphism of
$\cD_{\qu{Z}'/\qu{X}}$-modules for fixed $\qu{f}\in p^{-1}_Y(\cO_Y)$
and a differential operator in $\qu{f}$ for fixed $A\in P'_{\cD/X}$.
(To consider differential operators, we use the $\cO_{\qu{Z}'}$-module
structure on $P'_{\cD/X}$ induced by the embedding
$\cO_{\qu{Z}'}\hookrightarrow\cD_{\qu{Z}'/\qu{X}}$.)
\end{enumerate}}

\noindent
{\bf Proof.} Recall that $\cE(P)=\RR
p_{Y,*}\RR\sHom_{\cO_Z}(P,P_{\cD/X})$, and therefore
\[
{\mathbb
  H}^i(\cE(P))=\Ext^i_{\cO_Z}(P,P_{\cD/X})=\Ext^i_{\cD_{Z/X}}(P_{\cD/X},P_{\cD/X})
\]
by {Lemma~\ref{lem:rel_induced_adjunct}}.
The hypothesis that 
$\iota:\cO_Y\to\cE(P)$ is a quasi-isomorphism implies that
\[
\Ext^i_{\cD_{Z/X}}(P_{\cD/X},P_{\cD/X})
=\begin{cases}\Gamma(Y,\cO_Y)&(i=0)\\0,&(i>0).\end{cases}
\]
This implies part {\bf (a)} by the usual deformation theory argument
(cf. Proposition~\ref{prop:defcomodules}). Indeed, the sheaf of algebras
$\cD_{\qu{Z}'/\qu{X}}$ is an $R$-flat deformation of $\cD_{Z/X}$ by 
Proposition~\ref{prop:relqudiffops_flat}. Therefore, the category
$\Mod{\cD_{\qu{Z}'/\qu{X}}}$ is a flat
deformation of
$\Mod{\cD_{Z/X}}$ in the sense of \cite{Lowen_VanDerBergh}. Now
\cite[Theorem~A]{Lowen} implies that deformations of a
$\cD_{Z/X}$-module are controlled by its derived endomorphisms.

To prove {\bf (b)}, consider on $P'_{\cD/X}$ the filtration by modules
$\hbar^k P'_{\cD/X}$; the corresponding quotients are isomorphic to
$P_{\cD/X}$, so the associated graded module is $P_{\cD/X}\otimes_\C
R$. Since
\[
\RR p_{Y*}\RR\sHom_{\cD_{\qu{Z'}/\qu{X}}}(P'_{\cD/X},P_{\cD/X})=
 \RR p_{Y*}\RR\sHom_{\cD_{Z/X}}(P_{\cD/X},P_{\cD/X})=\cO_Y,
\]
we see that $\cE'$ is filtered by sheaves
\[
p_{Y*}\sHom_{\cD_{\qu{Z'}/\qu{X}}}(P'_{\cD/X},\hbar^k P'_{\cD/X})
\]
with the associated graded being $\cO_Y\otimes_\C R$, as required.

Finally, in {\bf (c)}, let $Q'=\sDiff_Z(p_Y^{-1}\cO_Y;P'_{\cD/X})$ be
the sheaf of differential operators $p_Y^{-1}\cO_Y\to P'_{\cD/X}$.
Note that for $A\in Q'$ and $B\in\cD_{\qu{Z}'/\qu{X}}$, the product
\[
AB:p_Y^{-1}\cO_Y\to P'_{\cD/X}:\qu{f}\mapsto A(\qu{f})B
\]
is again a differential operator. This turns $Q'$ into a right
$\cD_{\qu{Z}'/\qu{X}}$-module.  Note also that $Q'$ is flat over $R$
and that $Q'/\hbar Q'=\sDiff_Z(p_Y^{-1}\cO_Y;P_{\cD/X})$.  We have to
show that $m:P_{\cD/X}\to Q'/\hbar Q'$ lifts to a morphism of
$\cD_{\qu{Z}'/\qu{X}}$-modules $m':P'_{\cD/X}\to Q'$.

The obstruction to lifting $m$ lie in
\[
\Ext^1_{\cD_{\qu{Z}'/\qu{X}}}(P'_{\cD/X},\hbar Q').
\]
Note that $\hbar Q'\subset Q'$ admits a filtration by the modules $\hbar^k
Q'$ with quotients isomorphic to $Q'/\hbar Q'$.  It therefore suffices
to prove the vanishing of
\begin{myequation}{eq:equal.exts}
\Ext^1_{\cD_{\qu{Z}'/\qu{X}}}(P'_{\cD/X},Q'/\hbar Q') 
=\Ext^1_{\cD_{Z/X}}(P_{\cD/X},Q'/\hbar Q').
\end{myequation}
Actually, let us prove that
\[
\RR p_{Y*}\RR\sHom_{\cD_{Z/X}}(P_{\cD/X},Q'/\hbar Q')=\cD_Y.
\]
This implies that \eqref{eq:equal.exts} vanishes, because $Y$ is a
Stein manifold, and $\cD_Y$ is a direct limit of coherent
$\cO_Y$-modules, and therefore $H^{1}(Y,\cD_{Y}) = 0$.

Indeed, $Q'/\hbar Q'$ is equal to the sheaf of bidifferential
operators $p_Y^{-1}(\cO_Y)\times\cO_Z\to P$ that agree with the action
of $p_X^{-1}(\cO_X)$ on $\cO_Z$ and on $P$.  Note that even though
$P_{\cD/X}$ has two structures of an $\cO_Z$-module, both structures
lead to the same class of differential operators $p_Y^{-1}(\cO_Y)\to
P_{\cD/X}$.  We can then identify $Q'/\hbar Q'$ with
$p_Y^{-1}\cD_Y\otimes_{p_Y^{-1}\cO_Y} P_{\cD/X}$, where $\cO_Y$ acts
on both $\cD_Y$ and $P_{\cD/X}$ by left multiplication. We finally get
an isomorphism
\[
\RR
p_{Y*}\RR\sHom_{\cD_{Z/X}}(P_{\cD/X},p_Y^{-1}\cD_Y\otimes_{p_Y^{-1}\cO_Y}
P_{\cD/X})= 
\RR
p_{Y*}\RR\sHom_{\cD_{Z/X}}(P_{\cD/X},P_{\cD/X})\otimes_{\cO_Y}\cD_Y.
\]
\ \hfill $\Box$

\

\medskip

\punkt \ Let us now derive Proposition~\ref{prop:obstructionsneutr},
starting with the existence claim. For $P'_{\cD/X}$ provided by
Lemma~\ref{lem:obstructionsneutr}(a), set
$\qu{P}:=P'_{\cD/X}\otimes_{\cD_{\qu{Z}'/\qu{X}}}\cO_{\qu{Z}'}$.  Then
$\qu{P}$ is a flat sheaf of $R$-modules equipped with an
identification $\qu{P}/\hbar\qu{P}=P$.

Lemma~\ref{lem:obstructionsneutr}(b) now gives a 
 flat $R$-deformation of $\cO_Y$ to a sheaf of $R$-algebras 
 $\cE'$ on $Y$.  
 Note that $m'$ from
 Lemma~\ref{lem:obstructionsneutr}(c) can be viewed as a map
 $p_Y^{-1}(\cO_Y)\to\sEnd_{\cD_{\qu{Z}'/\qu{X}}}(P'_{\cD/X})$.  By
 adjunction, we obtain a map $\cO_Y\to\cE'$.  It follows
 from Lemma~\ref{lem:obstructionsneutr}(b) that the induced map of
 $R$-modules $\mu':\cO_Y\otimes_\C R\to\cE'$ is an isomorphism.

With respect to $\mu'$, the product on $\cE'$ can be written as
\[
\mu'(f)\mu'(g) = \mu'\left(\sum_{k=0}^n
B_k(f,g)\hbar^k\right)\quad(f,g\in\cO_Y\otimes_\C R)
\]
for some bilinear maps $B_k(f,g):\cO_Y\times\cO_Y\to\cO_Y$ with
$B_0(f,g)=fg$. We claim that $B_k(f,g)$ is a bidifferential operator
for any $k$ (that is, the product on $\cE'$ corresponds to a
$\star$-product on $\cO_Y\otimes_\C R$). Indeed, by construction
\[
m'(f,m'(g,A))=\sum_{k=0}^n m'(B_k(f,g),A)\hbar^k\quad(f,g\in\cO_Y,A\in
P'_{\cD/X}),
\]
and the left-hand side is a differential operator in both $f$ and
$g$. Thus by induction in $k$ the $B_{k}$'s are bidifferential
operators.  Define the neutralized $\star$-quantization $\qu{Y}$ of
$Y$ by setting $\cO^{op}_{\qu{Y}}=\cE'$.

\

\medskip

\punkt \ It remains to equip $\qu{P}$ with the structure of an
$\cO$-module on $\qu{X}\times\qu{Y}^{op}$.  Correcting $m'$ by $\mu'$,
we obtain a morphism
\[
m_Y=m'\circ(\mu')^{-1}:p_Y^{-1}(\cO_\qu{Y}^{op})\times P'_{\cD/X}\to
P'_{\cD/X}.
\]
It has the same properties as $m'$: it respects the
$\cD_{\qu{Z}'/\qu{X}}$-module structure on $P'_{\cD/X}$, and it is a
differential operator on $p_Y^{-1}(\cO_Y^{op})$. Moreover, it defines
an action of $p_Y^{-1}(\cO_Y^{op})$ on $P'_{\cD/X}$:
\[
m_Y(\qu{f},m_Y(\qu{g},A)) = m_Y(\qu{f}\qu{g},A)\quad(\qu{f},\qu{g}\in
p_Y^{-1}(\cO_Y^{op}), A\in P'_{\cD/X}).
\]
Recall now that $P'_{\cD/X}$ is a $\cO_{\qu{Z}'}$-module, in
particular, it has a natural action of $p_X^{-1}(\cO_\qu{X})$
\[
m_X:p_X^{-1}(\cO_\qu{X})\times P'_{\cD/X}\to P'_{\cD/X}.
\]
$m_X$ does not preserve the $\cD$-module (or even $\cO$-module)
structure on $P'_{\cD/X}$, but it is still a differential operator in
the $p_X^{-1}(\cO_\qu{X})$-variable.

\

\medskip

\punkt \
Set $\qu{Z}=\qu{X}\times\qu{Y}^{op}$. 
The maps $m_X$ and $m_Y$ define an $R$-bilinear morphism
\[
m_Z:\cO_\qu{Z}\times P'_{\cD/X}\to P'_{\cD/X}
\]
that is a differential operator on $\cO_\qu{Z}$ by the formula 
\[
m_Z(\qu{f}\qu{g},A)=m_X(\qu{f},m_Y(\qu{g},A))\quad(\qu{f}\in
p_X^{-1}(\cO_\qu{X}),\qu{g}\in p_Y^{-1}(\cO_\qu{Y}^{op}), A\in
P'_{\cD/X}).
\]
The map $m_Z$ is an action of $\cO_\qu{Z}$, because the actions $m_X$
and $m_Y$ commute. Finally, if we set
$I:=\{A\in\cD_{\qu{Z'}/\qu{X}}:A(1)=0\}$, we see that the action of
$\cO_\qu{Z}$ on $P'_{\cD/X}$ descends to an action on
\[
\qu{P} =
P'_{\cD/X}\otimes_{\cD_{\qu{Z'}/\qu{X}}}\cO_{\qu{Z}'} =
P'_{\cD/X}/P'_{\cD/X}I. 
\]
Thus $\qu{P}$ is an $\cO_\qu{Z}$-module. This proves part {\bfseries
(a)} of Proposition~\ref{prop:obstructionsneutr}

\

\medskip

\punkt \ Let us now prove part {\bfseries (b)}. Suppose $\qu{Y}\df$ is
a neutralized $\star$-quantization of $Y$, and $\qu{P}\df$ is an
$\cO$-module on $\qu{X}\times(\qu{Y}\df)^{op}$ with
$\qu{P}\df/\hbar\qu{P}\df=P$. Then the sheaf of differential operators
\[\qu{P}\df_{\cD/X}:=\sDiff_{\qu{Z}'/\qu{X}}(\cO_{\qu{Z}'},\qu{P}\df)\] 
is a flat right $\cD_{\qu{Z}'/\qu{X}}$-module that is a flat
$R$-deformation of $P_{\cD/X}$.  By
Lemma~\ref{lem:obstructionsneutr}(a), there is an isomorphism
$\phi:\qu{P}\df_{\cD/X}\to P'_{\cD/X}$.  Choose $\phi$ in such a way
that it reduces to the identity automorphism of $P_{\cD/X}$.  Note
that $\phi$ induces an isomorphism of $R$-modules
\[
\qu{P}\df =P\df_{\cD/X}\otimes_{\cD_{\qu{Z'}/\qu{X}}}\cO_{\qu{Z}'}
\to\qu{P}. 
\]
For any (local) section $\qu{f}\in\cO_{\qu{Y}\df}$, the action of
$p_Y^{-1}(\qu{f})$ on $\qu{P}\df$ commutes with the action of
$p_X^{-1}(\cO_\qu{X})$; therefore, $p_Y^{-1}(\qu{f})$ acts on
$\qu{P}\df_{\cD/X}$, so that we obtain an action
\[
m\df_Y:p_Y^{-1}(\cO_{\qu{Y}\df}^{op}) \times
\qu{P}\df_{\cD/X}\to\qu{P}\df_{\cD/X}.
\]
Under $\phi$, it translates into a morphism of sheaves of algebras
\[
p_Y^{-1}(\cO_{\qu{Y}\df}^{op}) \to
\sEnd_{\cD_{\qu{Z}'/\qu{X}}}(P'_{\cD/X}), 
\]
which by adjunction induces a homomorphism
$\cO_{\qu{Y}}^{op}\to\cE'=\cO_{\qu{Y}}^{op}$ (the sheaf $\cE'$ is
defined in Lemma~\ref{lem:obstructionsneutr}). By
Lemma~\ref{lem:obstructionsneutr}(b), the map is an isomorphism, so it
identifies $\qu{Y}^{op}$ and $\qu{Y}$.

\

\medskip

\punkt \ It remains to check that the identification
$\qu{P}=\qu{P}\df$ is $\cO_{\qu{Z}}$-linear.  Under the isomorphism
$\phi$, the actions $m_X$ and $m_Y$ correspond to the following
actions on differential operators $\qu{P}\df_{\cD/X}$:
\begin{align*}
m_X(\qu{f},A)(\qu{h})&=A(\qu{h}\qu{f})\qquad(A\in\qu{P}\df_{\cD/X},\qu{f}\in 
p_X^{-1}(\cO_{\qu{X}}),\qu{h}\in\cO_{\qu{Z}'}),\\
m_Y(\qu{g},A)(\qu{h})&=\qu{g}A(\qu{h})\qquad(A\in\qu{P}\df_{\cD/X},\qu{g}\in 
p_Y^{-1}(\cO_{\qu{Y}}),\qu{h}\in\cO_{\qu{Z}'}).
\end{align*}
The action $m_Z$ is more exotic, but still given by a differential
operator on $\cO_{\qu{Z}}$. Hence the induced action of $\cO_{\qu{Z}}$
on
\[
\qu{P}\df=\qu{P}\df_{\cD/X}\otimes_{\cD_{\qu{Z'}/\qu{X}}}\cO_{\qu{Z}'}=
\qu{P}\df_{\cD/X}/\qu{P}\df_{\cD/X}I
\]
is again given by a differential operator on $\cO_{\qu{Z}}$. Since
$p_X^{-1}(\cO_{\qu{X}})$ and $p_Y^{-1}(\cO^{op}_\qu{Y})$ act in the
natural way, so does $\cO_{\qu{Z}}$. This proves part {\bfseries (b)}.

\

\medskip

\punkt The proof of part {\bfseries (c)} given in
Section~\ref{ssec:obstructionsneutr} is completely straightforward (it
does not use $\cD$-modules). For this part, no restatement is
required.

This completes the proof of Proposition~\ref{prop:obstructionsneutr}.
\hfill $\Box$

\section{Computation of cohomology} \label{sec:coho}
In this section, we prove Theorem~\ref{thm:vanishing}.

\subsection{An example: The case of complex tori} \label{subsec:vanishing}
As a warmup we first consider the case of families of complex tori,
where the calculation is explicit. To prove the general case we take a
different tack but we hope that seeing a specific calculation first
will be illuminating.

\

\medskip

\punkt \ Suppose that
$f_{X} : X \to B$ is a smooth family of complex tori over a complex manifold $B$.
Here and everywhere else in this paper, a family of tori is supposed to admit a holomorphic
zero section $B\to X$, which is fixed.
Let $f_{Y} : Y \to B$ be the dual family of tori and set
$\Gamma := X\times_{B} Y$. Let $V \to \Gamma$ be
the normalized Poincar\'{e} line bundle, and let $P := i_{\Gamma
  *}V$ denote the corresponding sheaf on $Z = X\times Y$.

Note that $P$ is tangible and that the projection of $\Gamma=\supp(P)$
to $X$ is submersive. This implies
\[ 
\sExt^i(P,P_{\cD/X}) =\sExt^i(P,\sDiff_{Z/X}(\cO_Z;P))=0\qquad(i>0);
\]
this is a relative version of Lemma~\ref{lem:sDiff}. 
Now by Lemma~\ref{lem:rel_diffops_via_ind}, we have
\[
\RR\sHom(P,\sDiff_{Z/X}(\cO_Z;P)) 
=\sHom(P,\sDiff_{Z/X}(\cO_Z;P))=\sDiff_{Z/X}(P;P).
\]
Finally,
\[
\sDiff_{Z/X}(P;P)=i_{\Gamma *}\sDiff_{\Gamma/X}(V;V)
\]
by definition.  Thus Theorem~\ref{thm:vanishing} reduces to showing
that
\[\RR \op{pr}_{Y*} \sDiff_{\Gamma/X}(V;V)=
\mathcal{O}_{Y}.
\]
Here $\op{pr}_Y:=(p_Y)|_\Gamma:\Gamma\to Y$. 

\

\medskip

\punkt \
Consider the filtration of $\sDiff_{\Gamma/X}(V;V)$
\[
\sDiff_{\Gamma/X}^{\leq 0}(V;V) \subset \sDiff_{\Gamma/X}^{\leq
  1}(V;V) \subset \ldots \subset \sDiff_{\Gamma/X}^{\leq i}(V;V)
  \subset \ldots \subset \sDiff_{\Gamma/X}(V;V)
\]
by order of the differential operators on the fibers of
$\op{pr}_{X}:=(p_X)|_\Gamma : \Gamma \to X$. The spectral sequence for
the direct image of the filtered sheaf $\sDiff_{\Gamma/X}(V;V)$ along
the fibers of $\op{pr}_{Y}$ has $E_{1}$-term
\[
E^{ij}_{1} = \RR^{i}\op{pr}_{Y*}S^{j-i}T_{\Gamma/X},
\]
and converges to
  $\RR^{i+j}\op{pr}_{Y*}\sDiff_{\Gamma/X}(V;V)$.

\

\medskip

\punkt {\bf Lemma.}  {\it The differential
\[
d_{1} : \RR^{i}\op{pr}_{Y*}S^{j}T_{\Gamma/X} \to
\RR^{i+1}\op{pr}_{Y*}S^{j-1}T_{\Gamma/X}
\]
is given by the cup product with a class $\alpha_{\Gamma/X} \in
H^{0}(Y,\RR^{1}\op{pr}_{Y*}\Omega^{1}_{\Gamma/X})$. The class
$\alpha_{\Gamma/X}$ \linebreak  is the image of
\[
c_{1}(V) -
\frac{1}{2} c_{1}(\omega_{\Gamma/X}) \in
H^{1}(\Gamma,\Omega^{1}_{\Gamma})
\]
under the natural map
\[
H^{1}(\Gamma,\Omega^{1}_{\Gamma}) \to
H^{0}(Y,\RR^{1}\op{pr}_{Y*}\Omega^{1}_{\Gamma}) \to
H^{0}(Y,\RR^{1}\op{pr}_{Y*}\Omega^{1}_{\Gamma/X}),
\]
which we denote by $\op{pr}_{Y*}$.
}

\noindent
{\bf Proof.} This calculation is done in the absolute case in
\cite[Corollary~2.4.6]{beilinson-bernstein}. The proof of the relative
case presents no difficulties; the key step in the calculation is the
isomorphism of the graded $\mathbb{C}[\hbar]/\hbar^{2}$-Poisson algebra
corresponding to the Rees ring of $\sDiff(V;V)$ and the graded Poisson
algebra corresponding to the twisted cotangent bundle. Since this
isomorphism is written intrinsically in
\cite[Corollary~2.4.5]{beilinson-bernstein} it carries over
immediately to the relative context and so the proof follows. \ \hfill
$\Box$

\

\medskip

\punkt \ Next note that the relative dualizing sheaf
$\omega_{\Gamma/X}$ is trivial locally over $B$, since $\op{pr}_{X} :
\Gamma \to X$ is a family of complex tori. So, in order to understand
the spectral sequence computing the cohomology of
$\sDiff_{\Gamma/X}(V;V)$, we only need to understand the class
$\op{pr}_{Y*}c_{1}(V)$. However the standard Kodaira-Spencer theory
identifies the class
\[
\op{pr}_{Y*}c_{1}(V) \in \RR^{1}\op{pr}_{Y*} \Omega_{\Gamma/Y}^{1} \cong
\left(\RR^{1}\op{pr}_{Y*}\mathcal{O}_{\Gamma}\right) \otimes
\Omega_{Y/B}^{1} = f_{Y}^{*}\RR^{1}f_{X*}\mathcal{O}_{X}\otimes
\Omega_{Y/B}^{1}
\]
with the differential $d\xi_{V} \in
\op{Hom}(T_{Y/B},T_{\op{Pic}^{0}(X/B)/B})$ of the classifying map
\[
\xi_{V} : Y \to \op{Pic}^{0}(X/B)
\]
corresponding to the Poincar\'e line bundle $V$. This immediately implies 
Theorem~\ref{thm:vanishing} in our case:

\

\medskip

\punkt {\bf Corollary.} \label{cor:compute} {\it If $f_{X} : X \to B$
  and $f_{Y} : Y\to B$ are dual families of complex tori, and if $V \to
  \Gamma = X\times_{B} Y$ is the normalized Poincar\'e line bundle, then
  $\RR\op{pr}_{Y*}
  \sDiff_{\Gamma/X}(V;V) = \mathcal{O}_{Y}$. }

\noindent
{\bf Proof.} The differential of the spectral sequence is given by
the cup product with
\[
\alpha_{\Gamma/X} = \op{pr}_{Y*}\left(c_{1}(V) -
\frac{1}{2} c_{1}(\omega_{\Gamma/X})\right) = \op{pr}_{Y*}(c_{1}(V)) =
d\xi_{V} \in \op{Hom}(T_{Y/B},f_{Y}^{*}\RR^{1}f_{X*}\mathcal{O}_{X})
\]
and so
\[
\begin{split}
\left(\RR^{i}\op{pr}_{Y*}S^{j}T_{\Gamma/X},d_{1}\right) & =
\left(\RR^{i}\op{pr}_{Y*}\mathcal{O}_{\Gamma}\otimes
S^{j}T_{Y/B},d_{1}\right) \\
& = \left(f_{Y}^{*}\RR^{i}f_{X*}\mathcal{O}_{X}\otimes
S^{j}T_{Y/B},d_{1}\right) \\
& = \left(f_{Y}^{*}(\wedge^{i}\RR^{1}f_{X*}\mathcal{O}_{X})\otimes
S^{j}T_{Y/B}, \bullet\, \cup  \, d\xi_{V}\right)
\end{split}
\]
is just the Koszul complex. \ \hfill $\Box$

\

\subsection{Properties of the Fourier-Mukai transform}

In this section, we make some observations about the Fourier-Mukai
formalism on complex manifolds. They are used to derive
Theorem~\ref{thm:vanishing} in Section~\ref{ssec:vanishing}.

\

\medskip

\punkt {\bf Lemma.} \label{lem:fibers} {\it Let $M$ be a complex
manifold. Suppose that $F\in D^b_{coh}(M)$ satisfies $\LL i_x^* F=0$
for all embeddings of points $i_x:\{x\}\hookrightarrow M$, $x\in
M$. Then $F=0$.}

\noindent
{\bf Proof.} If $F\ne 0$, let $H^k(F)$ be the top non-vanishing
cohomology of $F$. Then $H^k(\LL i_x^* F)\ne 0$ for any
$x\in\op{supp}(H^k(F))$. \ \hfill $\Box$

\

\smallskip

\noindent Recall that $D^b_{comp}(M)\subset D^b_{coh}(M)$ denotes the
full subcategory of compactly-supported complexes. 

\

\medskip

\punkt {\bf Corollary.} \label{cor:fibers} {\it If $G\in D^b_{coh}(M)$
satisfies $\op{Hom}(F,G)=0$ for all $F\in D^b_{comp}(M)$, then
$G=0$. In other words, the right orthogonal complement to
$D^b_{comp}(M)$ in $D^b_{coh}(M)$ vanishes.}

\noindent
{\bf Proof.} For the embedding of a point $i_x:\{x\}\hookrightarrow
M$, we have
$$\op{Hom}(i_{x,*}{\mathbb{C}}[-k],G)\simeq H^{k-\dim M}(\LL i_x^* G).$$
Now Lemma~\ref{lem:fibers} implies the statement.\ \hfill
$\Box$

\

\medskip

\punkt \ Suppose now that $X$ and $Y$ are two complex manifolds. Fix 
$P\in D^b_{coh}(X\times Y)$ such that $\op{supp}(P)$ is proper over $Y$. Let
\[\Phi=\Phi^P:D^b_{comp}(Y)\to D^b_{comp}(X)\] 
be the integral transform with respect to $P$ given by
\begin{myequation}{eq:FM}
\Phi F =\RR p_{X*}(p_Y^*F\otimes^\LL P).
\end{myequation}

The Serre duality implies that $\Phi$ has a right  adjoint 
\[\Psi=\Psi^P:D^b_{coh}(X)\to D^b_{coh}(Y)\]
given by
\begin{myequation}{eq:adjoint}
\Psi F
=\RR p_{Y*}(\RR\sHom_{X\times Y}(P,p^*_XF))\otimes \omega_Y[\dim Y].
\end{myequation}
The fact that $\Phi$ and $\Psi$ are adjoint is well known if $X$ and $Y$ are both compact,
see e.g. \cite[Proposition~2.2]{FM} for the statement in the framework of algebraic varieties.
The analytic version that we use here relies on the duality theory for coherent sheaves
on complex manifolds, which is developed in \cite{RR}.

\punkt {\bf Remark.} \ Note that $\Psi$ need not preserve the
subcategory of compactly-supported objects. Because of this, the
adjunction relation between $\Phi$ and $\Psi$ is asymmetric. For
instance, a right adjoint of $\Phi$ (in this sense) is not always
unique.

\

\medskip

\punkt \ 
Set
\[
Q:=\RR p_{23*}(\RR \sHom_{X\times Y\times
Y}(p_{12}^*P,p_{13}^*P)\otimes p_3^*\omega_Y[\dim Y])\in
D^b_{coh}(Y\times Y),
\]
where $p_{12}$, for instance, is the projection of $X\times Y\times Y$
onto the product of the first two factors. Base change implies that
$Q$ is the kernel of $\Psi\circ\Phi$:
\[
\Psi\circ\Phi(F) =\RR p_{1*}(Q\otimes^\LL p_2^*F),\quad F\in
D^b_{comp}(Y).
\]
The adjunction between $\Phi$ and $\Psi$ yields a morphism of functors
\begin{myequation}{eq:adjunction}
\ID\to\Psi\circ\Phi.
\end{myequation} 
Clearly, the kernel of the identity functor is $\Delta_*\cO_Y$, where
$\Delta:Y\to Y\times Y$ is the diagonal morphism. We claim that
\eqref{eq:adjunction} is induced by a morphism of kernels
\begin{myequation}{eq:adjunctionkernel}
\Delta_*\cO_Y\to Q
\end{myequation} 
Let us construct \eqref{eq:adjunctionkernel}.

Set
\[
H:=\RR\sHom_{X\times Y\times Y}(p_{12}^*P,p_{13}^*P)\otimes
p_3^*\omega_Y[\dim Y]\in D^b_{coh}(X\times Y\times Y),
\]
so that $Q=\RR p_{23*}H$. The restriction of $H$ to the diagonal
$X\times Y\subset X\times Y\times Y$ is easy to describe:
\[
\LL(\op{id}_X\times\Delta)^*H= \RR\sHom_{X\times Y}(P,P)\otimes
p_Y^*\omega_Y[\dim Y].
\]
We therefore have
\[
\RR(\op{id}_X\times\Delta)^!H=\RR\sHom_{X\times Y}(P,P).
\]
The identity automorphism of $P$ gives a map
\[
\cO_{X\times Y}\to \RR\sHom_{X\times Y}(P,P) =
\RR(\op{id}_X\times\Delta)^!H,
\]
which by adjunction induces a morphism
\[
(\op{id}_X\times\Delta)_*\cO_{X\times Y}\to H.
\]
We finally give \eqref{eq:adjunctionkernel} as the composition
\[
\Delta_*\cO_Y\to \RR p_{23*}(p_{23}^*(\Delta_*\cO_Y)) = \RR
p_{23*}((\op{id}_X\times\Delta)_*\cO_{X\times Y})\to \RR p_{23*}H=Q.
\]
 
\

\medskip

\punkt {\bf Proposition.} \label{prop:adjunction}
{\it
The following conditions are equivalent:
\begin{itemize}
\item[{\bf (a)}] $\Phi:D^b_{comp}(Y)\to D^b_{comp}(X)$ is a fully
faithful functor;
\item[{\bf (b)}] The adjunction homomorphism \eqref{eq:adjunction} is
an isomorphism;
\item[{\bf (c)}] The morphism of kernels \eqref{eq:adjunctionkernel}
is an isomorphism.
\end{itemize}
}

\noindent
{\bf Proof.} $\text{\bf (a)}\Leftrightarrow\text{\bf(b)}$. By
 definition, $\Phi$ is fully faithful if and only if the map
\[
\op{Hom}_Y(F,G)\to\op{Hom}_X(\Phi(F),\Phi(G)) =
\op{Hom}_Y(F,\Psi\circ\Phi(G))
\]
is an isomorphism for all $F,G\in D^b_{comp}(Y)$. In other words, 
\[ 
\op{Hom}_Y(F,\op{Cone}(G\to\Psi\circ\Phi(G))) =0.
\]
By Corollary~\ref{cor:fibers}, this is true if and only if the map
$G\to\Psi\circ\Phi(G)$ is an isomorphism for all $G\in D^b_{comp}(Y)$.

$\text{\bf (c)}\Rightarrow\text{\bf(b)}$ is obvious.

$\text{\bf (b)}\Rightarrow\text{\bf(c)}$. Consider 
\[
Q':=\op{Cone}(\Delta_*\cO_Y\to Q)\in D^b_{coh}(X\times Y).
\]
The integral transform with respect to $Q'$ is
\[
\op{Cone}(\ID\to\Psi\circ\Phi):D^b_{comp}(Y)\to D^b_{coh}(X),
\]
which vanishes by {\bf (b)}. Applying the integral transform to
sky-scraper sheaves, we see that $Q'=0$ by Lemma~\ref{lem:fibers}.  \
\hfill $\Box$ 

\

\medskip

\punkt Let us now assume that $\supp(P)$ is proper over both $X$ and
$Y$. In this case, the situation is more symmetric: the formula
\eqref{eq:FM} defines a functor 
\[\Phi=\Phi^P:D^b_{coh}(Y)\to
D^b_{coh}(X).\] 
By Serre duality, $\Psi:D^b_{coh}(X)\to D^b_{coh}(Y)$
is the right adjoint of $\Phi$. Both $\Phi$ and $\Psi$ preserve the
subcategory of compactly-supported objects.

\

\medskip

\punkt {\bfseries Lemma.} \label{lem:compcoh} {\it The functor
$\Phi=\Phi^P:D^b_{coh}(Y)\to D^b_{coh}(X)$ is fully faithful if and
only if its restriction to $D^b_{comp}(Y)$ is fully faithful (which in
turn is equivalent to the other conditions of
Proposition~\ref{prop:adjunction}).}

\noindent
{\bfseries Proof.} The `only if' direction is obvious. Let us verify
the `if' direction.  It suffices to prove that the adjunction morphism
\eqref{eq:adjunction} is an isomorphism of functors from
$D^b_{coh}(Y)$ to itself. This follows because the morphism between
their kernels \eqref{eq:adjunctionkernel} is an isomorphism by
Proposition~\ref{prop:adjunction}.  \ \hfill $\Box$

\

\medskip

\punkt Finally, let us make some remarks about the bounded below
version of these categories. For any complex
manifold $M$, $D^+(\cO_M)$ is the bounded below derived category of
$\Mod{\cO_M}$. Consider the full subcategories
\begin{align*}
D^+_{coh}(M) &=\{\cF\in D^+(\cO_M):H^i(\cF)\in\Coh(M)\text{ for all }
i\}\\ D^+_{comp}(M) &=\{\cF\in D^+(\cO_M):H^i(\cF)\in\Coh(M)\text{ and
}\op{supp}(H^i(\cF))\text{ is compact for all } i\}.
\end{align*}
As before, suppose the kernel $P\in D^b_{coh}(X\times Y)$ is such that
$\op{supp}(P)$ is proper over $Y$. The same formula \eqref{eq:FM}
defines an extension of $\Phi$ to a functor $D^+_{comp}(Y)\to
D^+_{comp}(X)$, essentially because $P$ has finite
$\op{Tor}$-dimension. We still denote this extension by $\Phi$. If in
addition $\op{supp}(P)$ is proper over $X$, we obtain a functor
$\Phi:D^+_{coh}(Y)\to D^+_{coh}(X)$.

\

\medskip

\noindent
We need the following simple observation.
\

\smallskip

\punkt {\bfseries Lemma.} \label{lem:b+}
{\it Suppose the functor $\Phi:D^b_{comp}(Y)\to D^b_{comp}(X)$
is fully faithful.
\begin{itemize}
\item[{\bf (a)}] For any $\cF\in D^b_{comp}(Y)$, $\cG\in
D^+_{comp}(Y)$, the action of $\Phi$ on homomorphisms
\[
\Hom(\cF,\cG)\to\Hom(\Phi(\cF),\Phi(\cG))
\]
is a bijection.
\item[{\bf (b)}] Assume in addition that $\op{supp}(P)$ is proper over
$X$.  Then for any $\cF\in D^b_{coh}(Y)$, $\cG\in D^+_{coh}(Y)$, the
action of $\Phi$ on homomorphisms
\[
\Hom(\cF,\cG)\to\Hom(\Phi(\cF),\Phi(\cG))
\]
is a bijection.
\end{itemize}}

\noindent
{\bfseries Proof.} Consider the truncation triangle
\[
\tau^{\le N}\cG\to\cG\to\tau^{>N}\cG\to\tau^{\le N}\cG[1].
\]
For $N\gg 0$, we have
\[
\Hom(\cF,\tau^{>N}\cG)=\Hom(\Phi(\cF),\Phi(\tau^{>N}\cG))=0.
\]
Therefore,
\[
\Hom(\cF,\cG) =\Hom(\cF,\tau^{\le N}\cG)=\Hom(\Phi(\cF),\Phi(\tau^{\le
  N}\cG))=\Hom(\Phi(\cF),\Phi(\cG)).
\]
This proves {\bf (a)}. The proof of {\bf (b)} is similar.
\ \hfill $\Box$

\

\medskip

\punkt \ Similarly, one can consider the bounded above 
derived category $D^{-}(\cO_M)$ and its full subcategories $D^-_{coh}(M)$
and $D^-_{comp}(M)$ (where $M$ is a complex manifold). Lemma~\ref{lem:b+}
remains true if we assume that $\cF\in D^-_{comp}(Y)$ and $\cF\in D^-_{coh}(Y)$
in parts {\bf (a)} and {\bf (b)}, respectively.

\subsection{Proof of Theorem~\ref{thm:vanishing}}\label{ssec:vanishing}

\punkt \ We now turn to the proof of Theorem~\ref{thm:vanishing}. It is
convenient to give an alternative description of $\sDiff_{X\times
Y/X}(P;P)$. Fix $P\in\op{Coh}(X\times Y)$.

Recall that 
\[
P_{\cD/X}=\cD_{X\times Y/X}\otimes P.
\]
Applying differential operators in $\cD_{X\times Y/X}$ to functions on
$X\times Y$ that are constant along the fibers of $p_Y$, we can identify
\[
\cD_{X\times Y/X}=\sDiff_{X\times Y/X}(\cO_Y;\cO_{X\times
  Y})=p_Y^*\cD_Y.
\]
Recall that $\cD_Y$ has two actions of $\cO_Y$ (by left and right
multiplication), and that \linebreak $p_Y^*\cD_Y=\cD_{X\times Y/X}$
has two actions of $\cO_{X\times Y}$.  Both actions are used: the
tensor product $p_Y^*\cD_Y\otimes P$ is formed using the left action,
but then we consider it as an $\cO_{X\times Y}$-module using the right
action.

\

\medskip

\punkt \ The two actions of $\cO_Y$ on $\cD_Y$ allow us to view
$\cD_Y$ as a sheaf $D$ of $\cO_{Y\times Y}$-modules.  Let us agree
that the first, respectively second, factor in $Y\times Y$ corresponds
to the action of $\cO_Y$ on $\cD_Y$ by left, respectively right,
multiplication. Note that $D$ is not coherent: it is a union of
coherent sheaves corresponding to differential operators of bounded 
order.

In the same way, we view the inverse image $p_Y^*\cD_Y$ as the sheaf
$p^*_{23}D$ on $X\times Y\times Y$. The tensor product
$p_Y^*\cD_Y\otimes P$ is interpreted as $p^*_{23}D\otimes
p^*_{12}P$. We then consider $p_Y^*\cD_Y\otimes P$ as an $\cO_{X\times
Y}$-module using the right action of $\cO_{Y}$ on $\cD_{Y}$; in other words,
we identify
\[
P_{\cD/X} = p_Y^*\cD_Y\otimes P=p_{13*}(p^*_{23}D\otimes p^*_{12}P).
\]

\

\medskip 

\noindent Of course, the above formulas make sense for arbitrary $P\in
D^b_{coh}(X\times Y)$ that need not be concentrated in a single
cohomological degree. Let us prove Theorem~\ref{thm:vanishing} in this
generality.
\

\smallskip

\punkt {\bf Proposition.} \label{prop:derivedvanishing}
{\it Assume that the support of $P\in D^b_{coh}(X\times Y)$ is
proper over $Y$, and that the integral transform
\[
\Phi:D^b_{comp}(Y)\to D^b_{comp}(X)
\]
with respect to $P$ is fully faithful. Then
\[
\RR p_{Y*}\RR\sHom_{X\times Y}(P, p_{13*}(p^*_{23}D\otimes
p^*_{12}P))=\cO_Y.
\]
}

\noindent
{\bf Proof.}  By Proposition~\ref{prop:adjunction},
\eqref{eq:adjunctionkernel} is an isomorphism. It induces an
identification
\[
\Delta_*(\omega_Y^{-1})[-\dim(Y)] 
=\RR p_{23*}(\RR\sHom(p_{13}^*P,p_{12}^*P)).
\]
Note that compared to \eqref{eq:adjunctionkernel}, we permuted the two
copies of $Y$.

One can check that
\begin{myequation}{eq:Ddiag}
D\otimes^\LL\Delta_*\cO_Y
=\Delta_*\omega_Y[\dim(Y)].
\end{myequation} 
(This can be verified by an explicit local calculation using the Koszul resolution;
see Remark~\ref{rem:delta} for a more conceptual explanation.)
The projection formula gives
\[
\Delta_*\cO_Y =
D\otimes^\LL\Delta_*(\omega_Y^{-1})[-\dim(Y)]=
\RR p_{23*}(\RR\sHom(p_{13}^*P,p_{12}^*P)\otimes^\LL p_{23}^*D).
\]
It remains to take the direct image of both sides with respect to the
second projection $Y\times Y\to Y$. \ \hfill $\Box$

\

\medskip

\punkt {\bf Remark.}\label{rem:delta} Up to a twist, $D$ equals the $\cD_{Y\times
Y}$-module of $\delta$-functions supported on the diagonal.  More precisely,
$\cD_Y\otimes\omega^{-1}_Y$ has two structures of a left
$\cD_Y$-module, so $D\otimes p_2^*\omega^{-1}_Y$ is naturally a left
$\cD_{Y\times Y}$-module. This $\cD$-module is the direct image of the
constant $\cD_Y$-module $\cO_Y$ under $\Delta$. From this point of
view, \eqref{eq:Ddiag} simply states that the restriction of the
direct image to the diagonal is the constant $\cD_Y$-module.

\subsection{Extension of Fourier-Mukai transforms to
  $\star$-quantizations}\label{ssec:qFM} In this section, we study
Fourier-Mukai transforms on $\star$-quantizations, and prove
Theorem~B. A framework for working with Fourier-Mukai transforms on
$\star$-quantizations was built in \cite{KSarxiv}. Our situation is
somewhat simpler, because our $\star$-quantizations are parametrized
by an Artinian ring $R$, rather than $\C[[\hbar]]$, as in
\cite{KSarxiv}. This allows us to give more direct proofs.

\

\medskip

\punkt \ First, let $M$ be a complex manifold, and let
$\gb{M}$ be a 
$\star$-quantization of $M$. The category of $\cO_\gb{M}$-modules is
an abelian category with enough injectives; denote its (bounded)
derived category by $D^b(\cO_\gb{M})$. Recall that
$D^b_{coh}(\gb{M})\subset D^b(\cO_\gb{M})$ is the full subcategory of
complexes with coherent cohomology. Let $D^b_{comp}(\gb{M})\subset
D^b_{coh}(\gb{M})$ be the full subcategory of complexes whose coherent
cohomology {has} compact support.

\

\medskip

\punkt \ The category $\Mod{\cO_M}$ is identified with the full
subcategory of $\Mod{\cO_\gb{M}}$ consisting of modules annihilated by
$\hbar$. Geometrically, this identification is the direct image
functor
\[
i_*:\Mod{\cO_M}\to\Mod{\cO_\gb{M}}
\]
corresponding to the embedding $i:M\hookrightarrow\gb{M}$. The functor
$i_*$ is exact and preserves coherence and support. We therefore
obtain a functor
\[
i_*: D^b(\cO_M)\to D^b(\cO_\gb{M})
\]
such that $i_*(D^b_{coh}(M))\subset D^b_{coh}(\gb{M})$ and
$i_*(D^b_{comp}(M))\subset D^b_{comp}(\gb{M})$.

\

\medskip

\punkt \
For any $\qu{F}\in\Mod{\cO_\gb{M}}$, the tensor product 
\[
\qu{F}\otimes_R \C\in\Mod{\cO_\gb{M}}
\]
is annihilated by $\hbar$. We can therefore consider it as an object
$i^*(\qu{F})\in\Mod{\cO_M}$.

This defines a right exact functor
$i^*:\Mod{\cO_\gb{M}}\to\Mod{\cO_M}$.  Since $\Mod{\cO_\gb{M}}$ has
enough $R$-flat objects, we have a derived functor
\[
\LL i^*:D^b(\cO_\gb{M})\to D^-(\cO_M).
\]

Similarly, the sheaf of homomorphisms
\[
\sHom_R(\C,{\qu{F}})\in\Mod{\cO_\gb{M}}
\]
can be viewed as an object $i^!({\qu{F}})\in\Mod{\cO_M}$, this gives a left
exact functor $i^!:\Mod{\cO_\gb{M}}\to\Mod{\cO_M}$.  Since
$\Mod{\cO_\gb{M}}$ has enough injectives, we have a derived functor
\[
\RR i^!:D^b(\cO_\gb{M})\to D^+(\cO_M).
\]

\

\medskip

\punkt {\bfseries Lemma.} \label{lem:!*} {\it The functors $\LL i^*$
and $\RR i^!$ are left and right adjoint to $i_*$, respectively. Both
functors preserve the subcategories of coherent and compactly supported
coherent sheaves:
\begin{align*}
\LL i^*(D^b_{coh}(\gb{M}))&\subset D^-_{coh}(M) & 
\LL i^*(D^b_{comp}(\gb{M}))&\subset D^-_{comp}(M)\\
\RR i^!(D^b_{coh}(\gb{M}))&\subset D^+_{coh}(M) & 
\RR i^!(D^b_{comp}(\gb{M}))&\subset D^+_{comp}(M).
\end{align*}}

\noindent
{\bf Proof.} Let us check the properties of $\LL i^*$ (the properties
of $\RR i^!$ are proved in the same way).  The adjunction property is
standard; we can construct an identification
\[
\Hom(\cF,i_*\cG)=\Hom(\LL i^*\cF,\cG)\qquad(\cF\in
D^b(\cO_\gb{M}),\cG\in D^b(\cO_M))
\]
by using a $R$-flat resolution of $\cF$ and an injective resolution of
$\cG$.  To show that $\LL i^*\cF$ preserves the subcategories
$D^b_{coh}(\gb{M})$ and $D^b_{comp}(\gb{M})$, we note that $\LL
i^*\cF$ can be computed using a free resolution of the $R$-module
$\C$.  \ \hfill $\Box$

\

\medskip

\punkt \
Now let $X$ and $Y$ be complex manifolds. Let $\gb{X}$ and $\gb{Y}$ be
their $\star$-quantizations. Set $\gb{Z}=\gb{X}\times\gb{Y}^{op}$; it is a
$\star$-quantization of $Z=X\times Y$. 

Fix $\qu{P}\in D^b_{coh}(\gb{Z})$. Suppose that $\op{supp}\qu{P}$ is
proper over $Y$. Also,
let us assume that $\qu{P}$ has finite $\op{Tor}$-dimension as an $R$-module. 
Since $R$ is a local Artinian algebra, this is equivalent to boundedness of the derived
tensor product
\[\qu{P}\otimes^{\LL}_R\C.\] This derived tensor product is naturally an object in the
derived category of $\cO_Z$-modules; we can identify it with the derived restriction
\[P:=\LL i^*\qu{P}.\]
Thus, we see that $P\in D^b_{coh}(Z)$. 

Consider the integral transform functor

\[
\Phi=\Phi_P:D^+(\cO_Y)\to D^+(\cO_X)
\]
with the kernel $P$. We also have a Fourier-Mukai functor 
\[
\qu{\Phi}=\Phi_{\qu{P}}:D^b(\cO_\gb{Y})\to D^b(\cO_\gb{X})
\]
with the kernel $\qu{P}$. The functor $\qu{\Phi}$ is a version of the
`convolutions of kernels' from \cite{KSarxiv}.

\

\medskip

\punkt {\bfseries Remark.} Let us quickly review the definition of
$\qu{\Phi}$.  Let $\Mod{p_X^{-1}\cO_\gb{X}}$ be the category of
$p_X^{-1}\cO_\gb{X}$-modules.  Such modules can be defined as
$R$-linear homomorphisms $\gb{X}\to p_{X*}(\Sh_R)$, where
$p_{X*}(\Sh_R)$ is the stack that assigns to an open subset $U\subset
X$ the category $\Sh_R(p_X^{-1}(U))$ of sheaves of $R$-modules on
$p_X^{-1}(U)\subset Z$.  We have a direct image functor
\[ 
\RR p_{X*}:D^b(p_X^{-1}\cO_\gb{X})\to D^b(\cO_\gb{X}).
\]
Similarly, we can consider the category of
$p_Y^{-1}\cO_\gb{Y}$-modules and the inverse image functor
\[
p_Y^*:D^b(\cO_\gb{Y})\to D^b(p_Y^{-1}\cO_\gb{Y}).
\]

Finally, derived tensor product with $\qu{P}$ makes sense as a functor
\[
D^b(p_Y^{-1}\cO_\gb{Y})\to D^b(p_X^{-1}\cO_\gb{X}).
\]
We can then let $\qu{\Phi}$ be the composition
\[
D^b(\cO_\gb{Y})\to D^b(p_Y^{-1}\cO_\gb{Y})\to D^b(p_X^{-1}\cO_\gb{X})
\to D^b(\cO_\gb{X}).
\]

\

\medskip

\punkt {\bfseries Lemma.} \label{lem:restrictFM} {\it There is a
natural commutative diagram of functors
\[
\xymatrix{D^b(\cO_Y)\ar[r]^{i_*}\ar[d]^\Phi &
D^b(\cO_{\gb{Y}})\ar[r]^{\RR i^!}\ar[d]^{\qu\Phi} & D^+(\cO_Y)\ar[d]^\Phi\\
D^b(\cO_X)\ar[r]^{i_*} & D^b(\cO_{\gb{X}})\ar[r]^{\RR i^!}&
D^+(\cO_X).}\] } 

\noindent
{\bf Proof.} The fact that $i_*$ agrees with the Fourier-Mukai
transform follows directly from definitions. To see that $\RR i^!$ has
the same property, we can compute $\RR i^!$ using a free resolution of
the $R$-module $\C$-module.
\ \hfill $\Box$

\

\medskip

\noindent
Note that although we consider the integral transform $\qu{\Phi}$ on
$\cO$-modules, some finiteness conditions are necessary in order to
ensure that it has reasonable properties, see
\cite[Section~3.2]{KSarxiv}.
\

\smallskip

\punkt {\bfseries Lemma} {\it \begin{itemize}
\item[{\bf (a)}] $\qu{\Phi}(D^b_{comp}(\gb{Y}))\subset
  D^b_{comp}(\gb{X})$. 
\item[{\bf (b)}] Assume that  $\op{supp}({P})$
  is proper over $X$. Then $\qu{\Phi}(D^b_{coh}(\gb{Y}))\subset
  D^b_{coh}(\gb{X})$.
\end{itemize}
}

\noindent
{\bf Proof.} This claim is essentially a version of
\cite[Theorem~3.2.1]{KSarxiv}.  Note that the essential image
$i_*D^b_{comp}(X)$ generates $D^b_{comp}(\gb{X})$, in the sense that
the latter is the smallest triangulated category containing the
former.  This observation, together with Lemma~\ref{lem:restrictFM},
implies part {\bf (a)}.  In the same way, part {\bf (b)} follows from
the observation that $i_*D^b_{coh}(X)$ generates $D^b_{coh}(\gb{X})$.
\ \hfill $\Box$

\

\medskip

\noindent We are now ready to prove Theorem~B. It is implied by the
following statement, which is Theorem~\ref{thm:FM} with relaxed
assumptions on $\qu{P}$.
\

\smallskip

\punkt {\bfseries Proposition.} {\it As above, let $\gb{X}$, $\gb{Y}$
be $\star$-quantizations of complex manifolds $X$ and $Y$, and let
$\qu{P}\in D^b_{coh}(\gb{Z})$, where
$\gb{Z}=\gb{X}\times\gb{Y}^{op}$. Suppose that $\qu{P}$ has finite
$\op{Tor}$-dimension, and set $P=i^*\qu{P}\in D^b(Z)$, $Z=X\times
Y$. Finally, suppose that $\op{supp}(\qu{P})=\op{supp}(P)$ is proper
over $Y$, and that the integral transform functor
$\Phi=\Phi_P:D^b_{comp}(Y)\to D^b_{comp}(X)$ is fully faithful.

\begin{itemize}
\item[{\bf (a)}] The integral transform functor
$\qu{\Phi}=\Phi_{\qu{P}}:D^b_{comp}(\gb{Y})\to D^b_{comp}(\gb{X})$ is
fully faithful.
\item[{\bf (b)}] Let us also assume that $\op{supp}(\qu{P})$ is proper
over $X$. Then $\qu{\Phi}$ is also fully faithful as a functor
$D^b_{coh}(\gb{Y})\to D^b_{coh}(\gb{X})$.
\item[{\bf (c)}] Let us assume that $\op{supp}(\qu{P})$ is proper over
$X$ and that $\Phi$ provides an equivalence $D^b_{coh}(Y)\to
D^b_{coh}(X)$. Then
\[
\qu{\Phi}:D^b_{coh}(\gb{Y})\to D^b_{coh}(\gb{X})
\]
is an equivalence.
\end{itemize}
}

\noindent
{\bfseries Proof.} For part {\bf (a)}, it suffices to show that
$\qu{\Phi}$ is fully faithful when restricted to the essential image
$i_*(D^b_{comp}(Y))$. Indeed, given $\cF,\cG\in D^b_{comp}(Y)$, we
have
\begin{multline*}
\Hom(i_*\cF,i_*\cG)=\Hom(\cF,i^!i_*\cG)=\Hom(\Phi(\cF),\Phi(i^!i_*\cG))=
\Hom(\Phi(\cF),i^!\qu{\Phi}(i_*\cG))\\
=\Hom(i_*\Phi(\cF),\qu{\Phi}(i_*\cG))=\Hom(\qu{\Phi}(i_*\cF),\qu{\Phi}(i_*\cG))
\end{multline*}
by Lemmas~\ref{lem:b+}, \ref{lem:!*} and \ref{lem:restrictFM}.  The
proof of part {\bf(b)} is similar, except that we take $\cF,\cG\in
D^b_{coh}(Y)$.

Let us prove part {\bf(c)}. By {\bf (b)}, the functor $\qu{\Phi}$ is
fully faithful; therefore, its essential image is a triangulated
subcategory of $D^b_{coh}(\gb{X})$.  It suffices to verify that the
image contains $i_* D^b_{coh}(X)$. But this follows from
Lemma~\ref{lem:restrictFM}.  \ \hfill $\Box$

\section{Sheaves on $\star$-quantizations via equivariant stacks} \label{sec:equivariant}

In this paper, we work with many kinds of geometric objects (mostly
sheaves with some additional structure) on $\star$-quantizations. In
this chapter, we suggest a language for defining such objects in a
uniform way. To a large extend, this chapter is independent from the
rest of the paper; it may be considered an appendix.

\subsection{Equivariant stacks over the stack of quantizations}
\label{ssec:equivariance}

\punkt \ Let $M$ be a complex manifold.  Denote by  $M\an$ the
analytic site of $M$. Let $\sQ{M}$ be the stack of local neutralized
$\star$-quantizations of $M$. In other words,  $\sQ{M}$ is the stack of
groupoids on $M\an$ such that for an open subset $U \subset M$, the
fiber $\sQ{M}(U)$ is the groupoid of neutralized $\star$-quantizations
of $U$. There is a natural analytic topology on $\sQ{M}$: a family of
maps
\[
\left\{ \qu{V}_{i} \to \qu{U} \right\}_{i\in I}, \ \qu{V}_{i} \in
\sQ{M}(V_{i}), \ \qu{U} \in \sQ{M}(U)
\]
is a cover if and only if its image in $M\an$ is a cover, 
that is $U = \cup_{i \in I} V_{i}$. We write
$\sQ{M}\an$ for the analytic site of $\sQ{M}$. The site
$\sQ{M}\an$ comes equipped with a natural
sheaf $\cO_{\sQ{M}}$ of $R$-algebras, where
$\cO_{\sQ{M}}(\qu{U}) :=\Gamma(U,\cO_{\qu{U}})$.
Denote by  $\cO_{\sQ{M}}^\times$ the sheaf of
invertible elements in $\cO_{\sQ{M}}$; it is a sheaf of groups on $\sQ{M}\an$.

Let $I_{\sQ{M}}$ be the {\bf inertia sheaf} on $\sQ{M}\an$: for 
$\qu{U}\in\sQ{M}(U)$, 
$I_{\sQ{M}}(\qu{U})=\op{Aut}_{\sQ{M}(U)}(\qu{U})$
is the group of automorphisms of $\qu{U}$ that act trivially on $U$. 
Equivalently, $\op{Aut}_{\sQ{M}(U)}(\qu{U})$ is the group of automorphisms of the sheaf
$\cO_{\qu{U}}$ that act trivially modulo $\hbar$.

\

\medskip

\punkt {\bf Definition.} \ {\it Let $G$ be a sheaf of groups on
$\sQ{M}\an$.  An {\bfseries inertial} action of $G$ on $\sQ{M}$
is a group homomorphism $\rho : G \to I_{\sQ{M}}$ such that
\[
\rho(\alpha)^*(\beta)=\alpha^{-1}\beta\alpha\quad(\alpha,\beta\in
G(\qu{U}),\qu{U}\in\sQ{M}).
\]
}

\

\medskip

\punkt {\bf Example.} \ \label{ex:adjointinertial}
There is a natural inertial action $\op{Ad} :
  \cO_{\sQ{M}}^\times \to I_{\sQ{M}}$, defined as
  follows. Given $\qu{U} \in \sQ{M}\an$ and $g \in
  \cO_{\sQ{M}\an}^\times(\qu{U}) =
  \Gamma(U, \cO_{\qu{U}}^\times)$, we let
  $\op{Ad}(g)\in\op{Aut}_{\sQ{M}(U)}(\qu{U})$ be the automorphism
  of $\qu{U}$ that acts trivially on $U$ and acts on $\cO_{\qu{U}}$
  by conjugation
\[
\left(\op{Ad}(g)\right)^{*}(f) := g^{-1}fg, \quad f \in \cO_{\qu{U}}.
\]

\

\medskip

\punkt \ Let $\sX$ be a stack of categories over the site $\sQ{M}\an$.
Such a stack is the same as a stack over $M$ equipped with a
$1$-morphism $\iota_{\sX}: \sX \to \sQ{M}$.

Since $\sX$ can be considered as a stack on two different sites, we
have two categories of sections. Namely, given an open $U\subset M$
and $\qu{U}\in\sQ{M}(U)$, we consider the categories $\sX(\qu{U})$ and
$\sX(U)$: the essential fibers of the projections of $\sX$ to $\sQ{M}$
and $M\an$, respectively. The natural functor $\sX(\qu{U})\to\sX(U)$
is faithful, but not full.  In particular, for $\alpha\in\sX(\qu{U})$,
we have a natural embedding of automorphism groups:
\[ 
\op{Aut}_{\sX(\qu{U})}(\alpha)\hookrightarrow
\op{Aut}_{\sX(U)}(\alpha)=\op{Aut}_\sX(\alpha).
\]
Explicitly, $\op{Aut}_{\sX}(\alpha)$ is the group of automorphisms of
$\alpha$ in the total category of $\sX$ (note that $\sX(U)\subset\sX$
is a full subcategory), while $\op{Aut}_{\sX(\qu{U})}(\alpha)$
consists of automorphisms that act trivially on $\qu{U}$.

\

\medskip

\punkt {\bf Definition.} \label{def:equivariantstack} {\it Let $G$ be
a sheaf of groups on $\sQ{M}\an$ acting inertially on $\sQ{M}$ and let
$\sX$ be a stack of categories on $\sQ{M}\an$.  A
{\bfseries $G$-equivariant structure on $\sX$} is a collection of
group homomorphisms
\[
a_{\qu{U},\alpha} : G(\qu{U}) \to
\op{Aut}_{\sX}(\alpha)\quad(\qu{U}\in\sQ{M},\alpha \in \sX(\qu{U})) 
\]
satisfying the following conditions:
\begin{enumerate}
\item[(1)]  For any $\alpha \in \sX(\qu{U})$,
$\beta \in \sX(\qu{V})$, any homomorphism $\phi : \alpha \to
  \beta$  in $\sX$, and any $g \in G(\qu{V})$, the diagram
\[
\xymatrix{\alpha \ar[r]^-{\phi} \ar[d]_-{a_{\qu{U},\alpha}(g_{|\qu{U}})} & \beta
\ar[d]^-{a_{\qu{V},\beta}(g)} \\
\alpha \ar[r]_-{\phi} & \beta
}
\]
commutes.
\item[(2)] The composition $G(\qu{U})
  \stackrel{a_{\qu{U},\alpha}}{\to} \op{Aut}_{\sX}(\alpha)
  \to
  \op{Aut}_{\sQ{M}}(\qu{U})$ coincides with the inertial action of $G$.
\end{enumerate}
}

\

\medskip

\punkt {\bf Remark.} \ Suppose $\sX$ is a stack of groupoids. Then 
Definition~\ref{def:equivariantstack} is naturally reformulated in terms of
sheaves on $\sX$. Recall that the inertia sheaf $I_\sX$ assigns to $\alpha\in\sX$
its automorphism group $\op{Aut}_\sX(\alpha)$. It is obviously a presheaf on $\sX$, which
is a sheaf with respect to the natural topology on $\sX$ (a family of maps in $\sX$ is a cover
if and only if its image is a cover in $M\an$). Suppose we are given an
inertial action $\rho : G \to I_{\sQ{M}}$ of a sheaf of groups $G$ on $\sQ{M}$. Then a 
$G$-equivariant structure on $\sX$ is a homomorphism 
$a :\iota_\sX^*G \to I_\sX$ of sheaves of groups on
$\sX$ such that the composition
\[
\iota_\sX^*G \stackrel{a}{\to} I_{\sX} {\to}
\iota_\sX^* I_{\sQ{M}}
\]
coincides with $\iota_\sX^*(\rho) : \iota_\sX^*G \to \iota_\sX^*I_{\sQ{M}}$.

\

\medskip

\punkt \ Let $F:\sX\to\sY$ be a $1$-morphism of stacks on $\sQ{M}$. If both $\sX$ and $\sY$ 
are equipped with $G$-equivariant structures, it is clear what it means for $F$ to be {\bf $G$-equivariant}.

Now suppose $F$ is faithful. If $\sY$ has a $G$-equivariant structure, there is at most one equivariant structure on $\sX$ that makes $F$ an equivariant
$1$-morphism. If such a structure exists, we say that $\sX$ is a {\bf $G$-invariant}. Explicitly, $\sX$ is $G$-invariant if and only if for any 
$\qu{U}\in\sQ{M}$ and any $\alpha\in\sX(\qu{U})$,
the image of $\op{Aut}_\sX(\alpha)$ in $\op{Aut}_\sY(F(\alpha))$ contains the image of $G(\qu{U})$. In particular, 
if $F$ is fully faithful (so $\sX$ is a full substack of $\sY$), $\sX$ is automatically 
$G$-invariant.

\

\medskip

\punkt \ We are primarily interested in $\cO^\times_{\sQ{M}}$-equivariant stacks, where the inertial action of $\cO^\times_{\sQ{M}}$ on $\sQ{M}$ is the adjoint
action from Example~\ref{ex:adjointinertial}. (Sometimes it is natural to consider stacks equivariant under a `congruence subgroup' of $\cO^\times_{\sQ{M}}$,
see below.) Most relevant examples arise when we consider stacks of a kind of geometric objects, such as $\cO$-modules with some additional structure, on neutralized
$\star$-quantizations of open subsets of $M$. We conclude this section with some examples.

\

\medskip

\punkt {\bf Example.} \label{ex:Omod}
Let $\sX = \Omod{M}$ be the stack of left
  $\cO_{\sQ{M}}$-modules. More explicitly, $\sX(\qu{U})$ is the category of
  left $\cO_{\qu{U}}$-modules ($\qu{U}\in\sQ{M}$). The stack $\sX$ is naturally
  $\cO_{\sQ{M}}^\times$-equivariant: for a sheaf
  $\qu{P}$ of left $\cO_{\qu{U}}$-modules, we let $g \in
  \Gamma(U,\cO_{\qu{U}}^\times)$ act on the pair
  $(\qu{U},\qu{P})$ by the automorphism $(\op{Ad(g)}, \qu{m}(g))$,
  where $\qu{m}(g) : \qu{P} \to \qu{P}$ is the left multiplication by
  $g$. Condition (1) in the definition of
  $\cO_{\sQ{M}}^\times$-equivariance is satisfied
  tautologically. To verify condition (2), we must check that for any
  $f \in \Gamma(U, \cO_{\qu{U}})$ and any $s \in
  \Gamma(U,\qu{P})$, we have $\qu{m}(g)(fs) =
  (\op{Ad}(g^{-1})f)(\qu{m}(g)s)$. But by definition,
\[
(\op{Ad}(g^{-1})f)(\qu{m}(g)s) = (gfg^{-1})(gs) = 
gfs = \qu{m}(g)(fs).
\]

\

\medskip

\punkt {\bf Example.} \label{ex:Omodproperty} 
Let $\boldsymbol{\Sigma}$ be a local property of sheaves of $\cO$-modules 
on neutralized $\star$-quantizations. For instance, $\boldsymbol{\Sigma}$ could be one of the following properties of an $\cO$-module $\qu{P}$:
\begin{itemize}
\item $\qu{P}$ is coherent;

\item $\qu{P}$ is $R$-flat;

\item $\qu{P}$ is $R$-flat and its reduction $\qu{P}/\hbar\qu{P}$ is locally free, or tangible.
\end{itemize}

Denote by 
\[
\Omod{M}^{\boldsymbol{\Sigma}} \subset \Omod{M}
\]
the substack of $\cO$-modules having property $\boldsymbol{\Sigma}$. It is a full substack, so it is 
automatically  $\cO_{\sQ{M}}^\times$-invariant. 

\

\medskip

\punkt {\bf Example.} \label{ex:Omodpropertyrel}
There is also a relative version of these examples. Fix a neutralized 
$\star$-quantization $\qu{N}$ of a complex manifold $N$, and let 
$\sX = \Omod{M}/\qu{N}$  be the stack of left
  $\cO_{\sQ{M}\times \qu{N}}$-modules: its objects are pairs $(\qu{U},\qu{P})$, where
$\qu{U}\in\sQ{M}$, and $\qu{P}$ is a left $\cO_{\qu{U}\times\qu{N}}$-module. Similarly,
we can consider $\cO_{\qu{U}\times\qu{N}}$-module having some local property $\boldsymbol{\Sigma}$.

\

\medskip

\punkt {\bf Example.} \label{ex:tangiblerel} Another example is given in Section~\ref{ssec:tangiblerel}. Let
$p_M:Z\to M$ be a submersive morphism of complex manifolds. Given a neutralized $\star$-quantization
$\qu{U}$ of an open subset $U\subset M$, let $\sX(\qu{U})$ be the category of triples
$(P,\qu{P},\qu{\xi})$, where $P$ is a tangible sheaf on $p_M^{-1}(U)$, $\qu{P}$ is its $\star$-deformation,
and 
\[\qu{\xi}:p_M^{-1}\cO_{\qu{U}}\otimes_R\qu{P}\to\qu{P}\] 
is a $\star$-local action. The morphisms between $(P,\qu{P},\qu{\xi})$ and $(P',\qu{P}',\qu{\xi}')$ are
$R$-linear $\star$-local maps $\qu{P}\to\qu{P}'$ that commute with the action. The stack $\sX(\qu{U})$
is naturally $\cO^\times_{\sQ{M}}$-equivariant.

\

\medskip

\punkt {\bf Example.} \label{ex:stackabsolutecomodules}
Examples also arise from looking at modules over coalgebras corresponding to quantizations
(as in Section~\ref{ssec:quantcoalgebras}). Recall that a neutralized $\star$-quantization $\qu{U}$ of an open
subset $U\subset M$ gives rise to a coalgebra $\cA_\qu{U}=\sDiff(\cO_\qu{U};\cO_U\otimes_\C R)$ in the category
of left $\cD_U\otimes_\C R$-modules. We then consider the category $\Comod{\cA_{\qu{U}}}$ of $\cA_{\qu{U}}$-comodules
in the category of right $\cD_U\otimes_\C R$-modules. This provides a stack $\sX$ over $\sQ{M}\an$ defined by 
$\sX(\qu{U})=\Comod{\cA_{\qu{U}}}$. We claim that the stack is naturally $\cO^\times_{\sQ{M}}$-equivariant.

Such equivariant structure is provided by Lemma~\ref{lm:induct}. Indeed, let us also consider the stack $\sY$ given by
\[\sY(\qu{U})=\Mod{\cO_{\qu{U}}\otimes_\C\cD_U^{op}}.\]
This stack has a natural $\cO^\times_{\sQ{M}}$-equivariant structure. (Note that $\sY$ is a stack of modules over a 
sheaf of algebras on $\sQ{M}\an$; the sheaf of algebras contains $\cO_{\sQ{M}}$, which allows us to define the equivariant structure.)
The functor \eqref{eq:induct} provides a $1$-morphism $\sX\to\sY$, which is fully faithful by Lemma~\ref{lm:induct}. Therefore,
$\sX$ is $\cO^\times_{\qu{M}}$-invariant.

Note that by Proposition~\ref{prop:absolute_comodules}, the stack $\sX$ is actually identified with the stack of $\cO_{\sQ{M}}$-modules
$\Omod{M}$ from Example~\ref{ex:Omod}.

\

\medskip 

\punkt {\bf Example.} \label{ex:stackrelativecomodules}
It is more interesting to consider a relative version of Example~\ref{ex:stackabsolutecomodules}. Let $p_M:Z\to M$ be a submersive morphism
of complex manifolds. Given a neutralized $\star$-quantization $\qu{U}$ of an open subset $U\subset M$, the pullback $p_M^*\cA_\qu{U}$
is a coalgebra in the category of left $\cD$-modules on $p_M^{-1}(U)\subset Z$. Set 
$\sX(\qu{U})=\Comod{p_M^*\cA_{\qu{U}}}$. This defines a stack over $\sQ{M}$. 

The stack $\sX$ is naturally $\cO^\times_{\sQ{M}}$-equivariant. Indeed, Lemma~\ref{lm:relative_induct} gives a full embedding $\sX\to\sY$ into
the stack $\sY$ given by 
\[\sY(\qu{U})=\Mod{p_M^{-1}\cO_\qu{U}\otimes_\C\cD^{op}_{p_M^{-1}(U)}}.\]
The stack $\sY$ is naturally $\cO^\times_{\sQ{M}}$-equivariant, and $\sX$ is an invariant substack.

The stack $\sX$ contains the stack of Example~\ref{ex:tangiblerel} as a full substack. The embedding between the two stacks is given by
Lemma~\ref{lem:star-inducedrel}, which associates a $p_M^*\cA_{\qu{U}}$-comodule with a $\star$-deformation of a tangible sheaf equipped
with a local action of $p_M^{-1}(\cO_\qu{U})$. 

Note also that, as we explain in Section~\ref{ssec:relativecomodules}, the category of $p_M^*\cA_{\qu{U}}$-comodules can be identified with
the category of right $\cD_{\qu{V}/\qu{U}}$-modules, where $\qu{V}$ is any neutralized $\star$-quantization of $p_M^{-1}(U)\subset Z$ such that 
$p_M$ lifts to a map $\qu{V}\to\qu{U}$. Note that such $\qu{V}$ does not necessarily exist, and if it does, there may be more than one choice, even though
the category of right $\cD_{\qu{V}/\qu{U}}$-modules does not depend on the choice. At least in the case $Z=M\times N$, it is natural to take 
$\qu{V}=\qu{U}\times N$, and then we get an equivalent description of the stack $\sX$ in the language of relative $\cD$-modules on $\star$-quantizations.

\

\medskip 

\punkt {\bf Remark.} In all of the above examples, the $\cO^\times_\qu{M}$-equivariant structure on a stack $\sX$ actually comes from an action of 
$\cO_\qu{M}$ on $\sX$, in an appropriate sense. We leave the details of the definition to the reader.

\subsection{Sections of equivariant stacks over quantizations}\label{ssec:sectionsoverquant}

\punkt \label{ex:quant-equiv} Another source of
$\cO_{\sQ{M}}^\times$-equivariant stacks are $\star$-quantizations of $M$.
Let $\gb{M}$ be a $\star$-quantization of $M$, or more generally a $\star$-stack on $M$
(see Definition~\ref{def:Quantization}). By definition,
$\gb{M}$ is a stack of algebroids over $M$ and for any
open subset $U \subset M$ and any $\alpha \in \gb{M}(U)$, we get a natural
  neutralized $\star$-quantization $\qu{U}_{\alpha}\in\qu{M}(U)$: 
\begin{myequation}{eq:neutralize}
\qu{U}_{\alpha} =
\left(U,\sHom_{\gb{M}}(\alpha,\alpha)\right).
\end{myequation}
The assignment
\[
\gb{M}(U)\ni \alpha \mapsto
\qu{U}_{\alpha} \in \sQ{M}(U)
\]
is a morphism $\gb{M} \to \sQ{M}$
of stacks over $M$, so we can view $\gb{M}$ as a stack on $\sQ{M}\an$.
It is easy to see that $\gb{M}$ is naturally $\cO_{\sQ{M}}^\times$-equivariant. 
Indeed, if $U \subset M$ is open and $\alpha \in
\gb{M}(U)$, then
\[
\qu{U}_{\alpha} = \left( U,
\sEnd_{\gb{M}}(\alpha)\right) \qquad  \text{and}
\qquad
\cO_{\sQ{M}}^\times(\qu{U}_{\alpha}) =
\Gamma\left(U,\cO^\times_{\qu{U}_{\alpha}}\right)  =
\op{Aut}_{\gb{M}(U)}(\alpha)
\]
acts on $\alpha$ tautologically. It is clear that the action turns
$\gb{M}$ into an $\cO_{\sQ{M}}^\times$-equivariant stack.

Similarly we can consider the gerbe $\gb{M}^\times$ over $M$.
Thus, for every open subset $U\subset M$, we let $\gb{M}^\times(U)$
be the groupoid corresponding to the category $\gb{M}(U)$: it has the same objects,
and its morphisms are isomorphisms of $\gb{M}(U)$. We can view $\gb{M}^\times$
as a stack of groupoids on $\sQ{M}\an$. (Note that $\gb{M}^\times$ is {\it not} a gerbe over
$\sQ{M}\an$ if $\dim(M)>1$.) The same construction provides a natural 
$\cO_{\sQ{M}}^\times$-equivariant structure on $\gb{M}^\times$. In other words,
$\gb{M}^\times$ is an $\cO_{\sQ{M}}^\times$-invariant substack of $\gb{M}$.

\

\medskip

\punkt \ We now turn to definition of natural objects (such as sheaves of modules with additional
structure)
on $\star$-quantizations, or general $\star$-stacks. Formally, the construction applies
to objects that form an $\cO_{\sQ{M}}^\times$-equivariant stack over $\sQ{M}\an$, such as those considered
in Section~\ref{ssec:equivariance}. Given an $\cO_{\sQ{M}}^\times$-equivariant stack, we define
its sections over arbitrary $\star$-stacks using the following as follows. 

\

\medskip

\punkt {\bf Definition.} \label{def:equiv.mor} 
{\it Let $\sX$ be an
$\cO_{\sQ{M}}^\times$-equivariant stack, and let $\gb{M}$ be a $\star$-stack on $M$. Set
$\sX(\gb{M})$ (the category of {\bf sections} of $\sX$ over $\gb{M}$) to be the category
of $\cO_{\sQ{M}}^\times$-equivariant $1$-morphisms $\gb{M}^\times\to\sX$ of stacks. Here
we view $\gb{M}^\times$ with the equivariant structure described in \ref{ex:quant-equiv}.}

\

\medskip

\punkt \
Applying Definition~\ref{def:equiv.mor} to equivariant stacks from examples in Section~\ref{ssec:equivariance},
we define $\cO$-modules on $\gb{M}$, coherent $\cO_\gb{M}$-modules, $\cO_M-\cD_\gb{M}$-bimodules, and so on. 
This provides a `uniform' reformulation of Definitions~\ref{def:modules_gerby}, \ref{def:modules.coherent}, \ref{def:ODbimod},
and \ref{def:ODbimodrel}. The equivalence of the two styles of definitions is the subject of Section~\ref{ssec:linearity}.
 
Note that if $\gb{M}$ is neutral, Definition~\ref{def:equiv.mor} is consistent. This follows from simple abstract observations about stacks 
(Proposition~\ref{prop:mapsfromgerbe}).

\

\medskip

\noindent
 Let $\sX$ be a gerbe over $M$ and let $\sY$ be an arbitrary stack. Suppose
$\sX$ admits a global neutralization $\alpha\in\sX(M)$. Any $1$-morphism of stacks
$F:\sX\to\sY$ yields an object $\beta=F(\alpha)\in\sY(M)$ and a morphism of sheaves
$\phi:\sAut_\sX(\alpha)\to\sAut_\sY(\beta).$
\

\smallskip

\punkt {\bf Proposition.}\label{prop:mapsfromgerbe} {\it
\begin{itemize}
\item[{\bf (a)}] The above correspondence provides an equivalence between the
category of $1$-morphisms $F:\sX\to\sY$ and the category of pairs
$$(\beta\in\sY(M),\phi:\sAut_\sX(\alpha)\to\sAut_\sY(\beta)).$$

\item[{\bf (b)}] Suppose in addition that $\sX$ and $\sY$ are $G$-equivariant stacks on
$\sQ{M}$ for a sheaf of groups $G$ that acts inertially on $\sQ{M}$. Assume that the action
of $G$ on $\sX$ is {\bf simple transitive}: for any $\qu{U} \in\sQ{M}$ and any $\gamma\in\sX(\qu{U})$, the map $G(\qu{U})\to\op{Aut}_{\sX}(\gamma)$ is an isomorphism.
Then the correspondence
$$F\mapsto F(\alpha)$$
yields an equivalence between the category of $G$-equivariant $1$-morphisms $\sX\to\sY$ and 
$\sY(\qu{M}_\alpha)$. Here $\qu{M}_\alpha\in\sQ{M}(M)$ is the image of $\alpha\in\sX(M)$.
\end{itemize}}

\noindent
{\bf Proof.} {\bf (a)} It is easy to see that the correspondence is faithful. Indeed, we need to show that given two $1$-morphisms $F,F':\sX\to\sY$
and two $2$-morphisms $f,f':F\to F'$ such that $f(\alpha)=f'(\alpha):F(\alpha)\to F'(\alpha)$, we have $f(\gamma)=f'(\gamma):F(\gamma)\to F'(\gamma)$ for
any $\gamma\in\sX$. This follows because $\gamma$ is locally isomorphic to $\alpha$.

To see that the correspondence is full, we need to show that for two such $1$-morphisms $F$ and $F'$,
any morphism $f:F(\alpha)\to F'(\alpha)$ that agrees with the action $\phi$ extends to a morphism
of functors. To do this, we need to construct $f(\gamma):F(\gamma)\to F'(\gamma)$ for all $\gamma\in\sX$. Again, this can be easily done locally: an isomorphism 
$\alpha|_U\simeq\gamma|_U$ for an open set $U$ allows us to view $f|_U$ as a map
$$f_{\gamma,U}:F(\gamma)|_U\to F'(\gamma)|_U.$$
Since $f$ agrees with $\phi$, the map $f_{\gamma,U}$ does not depend on the choice of isomorphism
$\gamma|_U\simeq\alpha|_U$; this allows us to glue the local maps $f_{\gamma,U}$ into a map $f(\gamma)$.

Finally, let us show the correspondence is essentially surjective. Given $\beta\in\sY(M)$ and
$\phi:\sAut_\sX(\alpha)\to\sAut_\sY(\beta)$, let us construct the corresponding functor $F$. 
Given $\gamma\in\sX(U)$ for open $U$, we choose an open cover $U=\bigcup U_i$ together with isomorphisms $$s_i:\alpha|_{U_i}\to\gamma|_{U_i}.$$ We then obtain a descent datum for $\sY$:
on sets $U_i$, we have objects $\beta|_{U_i}\in\sY(U_i)$, and on intersections $U_i\cap U_j$, we have isomorphisms $\phi(s_j^{-1}s_i)$ between the restrictions of these objects. Since $\sY$ is
a stack, the descent datum gives rise to an object $\delta\in\sY(U)$, unique up to a unique isomorphism.
We set $F(\gamma)=\delta$. It is easy to see that $F$ defined in this way is indeed a $1$-morphism
of stacks. 

Part {\bf (b)} follows from part {\bf (a)}. Indeed, consider a sheaf $G_\alpha$ on $M$ that
assigns to an open set $U\subset M$ the group
$G(\qu{U}_\alpha)$, where $\qu{U}_\alpha\in\sQ{M}(U)$ is the image of $\alpha|_U$ (which is an open subset of
$\qu{M}_\alpha$). The $G$-equivariant
structure on $\sX$ defines a morphism $a_\alpha:G_\alpha\to\sAut_\sX(\alpha)$; similarly, the $G$-equivariant structure on $\sY$ defines a morphism $a_\beta:G_\alpha\to\sAut_\sY(\beta)$ for any $\beta\in\sY(\qu{M}_\alpha)$. A pair $(\beta,\phi)$ as in part {\bf (a)} corresponds to
a $G$-equivariant $1$-morphism $\sX\to\sY$ if and only if $\phi a_\alpha=a_\beta$. Since
the action of $G$ on $X$ is simple transitive,  $a_\alpha$ is an isomorphism, so if $\beta$ is given, $\phi$ is uniquely determined.
\ \hfill $\Box$

\

\medskip

\punkt {\bf Corollary.} \ \label{cor:trivialize}
{\it Suppose $\gb{M}$ admits a global neutralization $\alpha \in
\gb{M}(M)$. Let $\qu{M}_{\alpha}\in\sQ{M}$ be the corresponding global
neutralized quantization, as in \eqref{eq:neutralize}. Then for any
$\cO_{\sQ{M}}^\times$-equivariant stack $\sX$, the the
  categories $\sX(\gb{M})$ and
  $\sX(\qu{M}_{\alpha})$ are naturally equivalent.
The equivalence is given by
$$F\mapsto F(\alpha)\in\sX(\qu{M}_\alpha)\qquad(F\in\sX(\gb{M}));$$
recall that $F\in\sX(\gb{M})$ is an
$\cO_{\sQ{M}}^\times$-equivariant 1-morphism $\gb{M}^\times\to\sX$.
}

\noindent
{\bf Proof.}  Follows from  Proposition~\ref{prop:mapsfromgerbe}(b), because $\gb{M}^\times$
is a gerbe on which $\cO_{\sQ{M}}^\times$ acts simply transitively. \ \hfill $\Box$

\

\medskip

\punkt \ Let us also mention a version of Definition~\ref{def:equiv.mor} that is specifically
adapted to $\star$-quantizations of $M$, rather than $\star$-stacks.
Recall (Definition~\ref{def:*quantization}) that a $\star$-quantization of $M$ is a 
$\star$-stack $\gb{M}$ whose reduction to $\mathbb{C}$ is neutralized. This additional structure
allows us to define a wider class of objects on $\gb{M}$. From the point of view of equivariant
stacks, we can define $\sX(\gb{M})$ for a stack $\sX$ that is equivariant with respect to a group smaller than $\cO_{\sQ{M}}^\times$. 

Let $\gb{M}$ be a $\star$-quantization of $M$. For any open set $U\subset M$ and any $\alpha\in\gb{M}(U)$, let $\alpha'$ be the reduction of $\alpha$ to $\mathbb{C}$; then
$\alpha'\in\gb{M}'(U)$, where $\gb{M}'$ is the reduction of $\gb{M}$ to $\mathbb{C}$. Recall
that $\gb{M}'$ is equipped with a global neutralization $\beta\in\gb{M}'(M)$.

We define $\gb{M}^{(1)}$ to be the stack of collections $(\alpha,s')$, where $\alpha\in\gb{M}$,
and $s'$ is the trivialization of $\alpha'$. That is, $s'$ is an isomorphism between 
$\alpha'$ and the restriction $\beta|_U$. Clearly, $\gb{M}^{(1)}$ is a stack of groupoids on $M$.

We have a natural faithful functor $\gb{M}^{(1)}\to\gb{M}^\times$; in particular, it allows us to
view $\gb{M}^{(1)}$ as a stack over $\sQ{M}$. Note that $\gb{M}^{(1)}$ is not invariant with respect
to the canonical action of $\cO^\times_{\sQ{M}}$ on $\gb{M}^\times$ (described in \ref{ex:quant-equiv}). However, $\gb{M}^{(1)}$ is invariant for the 
`congruence' subgroup $\cO_{\sQ{M}}^{(1)}\subset\cO_{\sQ{M}}^\times$ defined by
\[\cO_{\sQ{M}}^{(1)}(\qu{U}):=\{\qu{f}\in\Gamma(U,\cO^\times_{\qu{U}})|\qu{f}=1\mod\hbar\}\qquad(\qu{U}\in\sQ{M}).\]
Thus $\gb{M}^{(1)}$ is naturally an $\cO_{\sQ{M}}^{(1)}$-equivariant stack.
It is easy to see that $\gb{M}^{(1)}$ is a gerbe over $M$ on which $\cO_{\sQ{M}}^{(1)}$ acts simply and transitively.  

\

\medskip

\punkt {\bf Remark.} Let $\gb{M}$ be a $\star$-stack. To turn $\gb{M}$ into a $\star$-quantization, we need to choose a trivialization of the reduction of $\gb{M}$ to $\mathbb{C}$.
Alternatively, one can view this choice as a reduction of the `structure group' of 
$\gb{M}^\times$ from $\cO^\times_{\sQ{M}}$ to $\cO^{(1)}_{\sQ{M}}$. Here we essentially identify
$\star$-stacks with $\cO^\times_{\sQ{M}}$-gerbes and $\star$-quantizations with 
$\cO^{(1)}_{\sQ{M}}$-gerbes over $\sQ{M}$. We come back to this point of view in 
Section~\ref{ssec:2stacks}.

\

\medskip

\noindent
We can now modify Definition~\ref{def:equiv.mor}.
\

\smallskip

\punkt {\bf Definition.} \label{def:equiv.mor1}
{\it Let $\sX$ be an
$\cO_{\sQ{M}}^{(1)}$-equivariant stack, and let $\gb{M}$ be a $\star$-quantization of $M$. Set
$\sX(\gb{M})$ (the category of {\bf sections} of $\sX$ over $\gb{M}$) to be the category
of $\cO_{\sQ{M}}^{(1)}$-equivariant $1$-morphisms $\gb{M}^{(1)}\to\sX$ of stacks.}

\

\medskip

\punkt \ Definition~\ref{def:equiv.mor1} agrees with other situations when the category of 
sections is defined. First of all, suppose the $\star$-quantization $\gb{M}$ is neutralized. In other
words, we are given $\alpha\in\gb{M}(M)$ and an isomorphism between the reduction $\alpha'$ of 
$\alpha$ to $\mathbb{C}$ and the standard neutralization of $\gb{M}'$ (which is part of the 
$\star$-quantization structure). Then $\sX(\gb{M})$ is naturally equivalent to $\sX(M_\alpha)$; the
statement and its proof are completely parallel to Corollary~\ref{cor:trivialize}.

Also, Definitions~\ref{def:equiv.mor1} and Definitions~\ref{def:equiv.mor} agree whenever they both make sense. (We give the statement in concrete terms, but the proof is an abstract
stack argument. The abstract formulation of the statement is left to the reader.)

\

\medskip

\punkt {\bf Lemma.} {\it Let $\gb{M}$ be a $\star$-quantization of $M$, and let $\sX$ be
an $\cO^\times_{\sQ{M}}$-equivariant stack. Given an $\cO^\times_{\sQ{M}}$-equivariant $1$-morphism
$\gb{M}^\times\to\sX$, the composition $\gb{M}^{(1)}\to\gb{M}^\times\to\sX$ is clearly
$\cO^{(1)}_{\sQ{M}}$-equivariant. We claim that this provides an equivalence between the category
of $\cO^\times_{\sQ{M}}$-equivariant $1$-morphisms
$\gb{M}^\times\to\sX$ and the category of $\cO^{(1)}_{\sQ{M}}$-equivariant $1$-morphisms
$\gb{M}^{(1)}\to\sX$.}

\noindent
{\bf Proof.} Given two stacks $\sY_1,\sY_2$ over $M$, we can consider the stack of $1$-morphisms $\sY_1\to\sY_2$, whose category of global sections over open $U\subset M$ is the category of
$1$-morphisms $\sY_1|_U\to\sY_2|_U$. Similarly, we can consider stacks of equivariant 
$1$-morphisms between equivariant stacks. We can thus restate the lemma using stacks: we need to show the natural $1$-morphism between the stack of $\cO^\times_{\sQ{M}}$-equivariant $1$-morphisms
$\gb{M}^\times\to\sX$ and the stack of $\cO^{(1)}_{\sQ{M}}$-equivariant $1$-morphisms
$\gb{M}^{(1)}\to\sX$ is an equivalence. 

In this formulation, the claim becomes local on $M$, so we can assume without losing generality that
$\gb{M}$ is a neutralized $\star$-quantization. Then the claim follows from 
Proposition~\ref{prop:mapsfromgerbe}(b) (cf. Corollary~\ref{cor:trivialize}.) \ \hfill $\Box$

\

\medskip

\punkt {\bf Examples.} Obviously, any $\cO^\times_{\sQ{M}}$-equivariant stack $\sX$, such as those considered in
Section~\ref{ssec:equivariance}, is also equivariant for $\cO^{(1)}_{\sQ{M}}$ (and any other subgroup of 
$\cO^\times_{\sQ{M}}$). We obtain other examples of $\cO^{(1)}_{\sQ{M}}$-equivariant stacks by looking at objects
whose reduction modulo $\hbar$ is fixed. For instance, let us fix an $\cO_M$-module $P$. Consider a stack $\sX$
over $\sQ{M}$ such that for any open $U\subset M$ and any $\qu{U}\in\sQ{M}(U)$, the category $\sX(U)$ consists of pairs
$(\qu{P},i)$, where $\qu{P}$ is a $R$-flat $\cO_\qu{U}$-module and 
\[i:\qu{P}/\hbar\qu{P}\to P|_U\]
is an isomorphism of $\cO_U$-modules. It is not hard to see that for any $\star$-quantization $\gb{M}$ of $M$, the category
of sections $\sX(\gb{M})$ is the category of deformations of $P$ to a $\cO_\gb{M}$-module. Similar examples can be constructed for other stacks
from Section~\ref{ssec:equivariance}.

\subsection{Equivariance and $R$-linearity} \label{ssec:linearity}

Sheaves on quantizations can also be described using representations
of stacks of algebroids (as in for instance \cite{KSarxiv}). 
Let us compare this approach to the framework of inertial actions of
$\cO_{\sQ{M}}^\times$. 

For example, Definition~\ref{def:modules_gerby} defines $\cO$-modules on
a $\star$-stack $\gb{M}$ to be {\bfseries representations} of $\gb{M}$, that is,
$1$-morphisms
$\qu{F} : \gb{M} \to \Sh_{R}$ of stacks of  $R$-linear
categories over $M$. Here $\Sh_{R}$ be the stack of sheaves of $R$-modules on $M$: to an open
$U\subset M$, it assigns the category $\Sh_R(U)$ of sheaves of $R$-modules on $U$.

\

\medskip

\punkt \ 
Let $\qu{F} : \gb{M} \to \Sh_{R}$ be a
representation of $\gb{M}$. For any open $U\subset M$ and any $\alpha\in\gb{M}(U)$,
$\qu{F}(\alpha)$ is a sheaf of $R$-modules on $U$. Since $\qu{F}$ is a morphism of stacks of $R$-linear categories,
the action of $R$ on $\qu{F}(\alpha)$ naturally extends to an action
of the sheaf of $R$-algebras
\[
\mathcal{O}_{\qu{U}_{\alpha}}=\sEnd_{\gb{M}}(\alpha).
\]
We denote
the resulting $\mathcal{O}_{\qu{U}_{\alpha}}$-module by
$\qu{F}_{(\alpha)}$.
In particular, if $\alpha \in \gb{M}(M)$ is a global section, then
we get a functor
\[
\xymatrix@R-2pc{
\left(\text{
representations of $\gb{M}$}
\right)  \ar[r] &
\left(\text{$\mathcal{O}_{\qu{M}_{\alpha}}$-modules}\right) \\
\qu{F}  \ar[r] & \qu{F}_{(\alpha)},
}
\]
which is an equivalence. For general $\gb{M}$, a representation of $\gb{M}$
can be explicitly described as a collection of
$\mathcal{O}_{\qu{U}_{\alpha}}$-modules over neutralized open patches
together with gluing conditions as in \cite{PS}.

\

\medskip

\noindent
One can check that the two approaches to $\cO_\gb{M}$-modules lead to the same result.
\

\smallskip

\punkt {\bfseries Proposition.}\label{prop:linequiv} 
{\it Let $\gb{M}$ be a $\star$-stack over $M$. 
There is an equivalence between the category of representations of $\gb{M}$ and the category of 
$\mathcal{O}_{\sQ{M}_{\op{an}}}^{\times}$-equivariant $1$-morphisms
from $\gb{M}^{\times}$ to the stack $\sX=\Omod{M}$, which we denoted by $\sX(\gb{M})$. The equivalence 
assigns to $\qu{P}\in\sX(\gb{M})$ the composition
\[\xymatrix@-1pc{\gb{M}\ar[r]^{\qu{P}}&\sX\ar@{=}[r]&\Omod{M}\ar[r]&\Sh_R}.\]
Here the last arrow is the forgetful $1$-morphism.} 

\noindent
{\bf Proof.} Let us describe the inverse correspondence.
Suppose that $\qu{F} : \gb{M} \to \Sh_{R}$
is a representation of $\gb{M}$. Then  the $1$-morphism
$\qu{F}|_{\gb{M}^{\times}}$ naturally lifts
to a $1$-morphism
\[
\xymatrix@-1pc{
\gb{M}^{\times} \ar[rr]^-{\qu{F}_{(\bullet)}} \ar[dr]_-{\qu{F}} & & \Omod{M}
\ar[dl] \\
& \Sh_{R} &
}
\]
where
\[
\xymatrix@R-2pc{
\widetilde{F}_{(\bullet)}  : \hspace{-0.1in} & \gb{M}^{\times}
\ar[r] & \Omod{M} \\
& \alpha \ar@{|->}[r] & \left( \widetilde{U}_{\alpha},
\qu{F}_{(\alpha)} \right).
}
\]
Moreover, by construction $\widetilde{F}_{(\bullet)}$ is equivariant
with respect to the natural
$\mathcal{O}_{\sQ{M}_{\op{an}}}^{\times}$-equivariant structures on
$\gb{M}^{\times}$ and $\Omod{M}$.  Thus $\widetilde{F}_{(\bullet)}$ is
naturally an object of $\sX(\gb{M})$, as required.
 \ \hfill $\Box$

\

\medskip

\punkt \ Proposition~\ref{prop:linequiv} can be generalized in the following manner. 
Let ${\sf{S}}$ be a stack of $R$-linear categories over
$M_{\op{an}}$ (above we have ${\sf{S}} = \Sh_{R}$) and let $\sX$
be the stack of $\mathcal{O}_{\sQ{M}_{\op{an}}}$-modules in ${\sf{S}}$. 
Concretely $\sX$ is the stack of categories
on $\sQ{M}_{\op{an}}$ such that, for every $\qu{U} \in
\sQ{M}_{\op{an}}$, the category $\sX(\qu{U})$ is the 
category of
$\cO_{\qu{U}}$-modules in ${\sf{S}}|_U$. In other words,
$\sX(\qu{U})$ is the category of pairs $(\alpha,a)$, where
$\alpha \in {\sf{S}}(U)$ and $a : \cO_{\qu{U}} \to
\sEnd_{{\sf{S}}}(\alpha)$ is a homomorphism of sheaves of $R$-algebras
on $U$.

\

\medskip

\punkt {\bf Proposition.} \label{prop:linequivgen} {\it For any $\star$-quantization $\gb{M}$
of $M$, there is a natural equivalence between the category of
$R$-linear $1$-morphisms $\gb{M} \to {\sf{S}}$ of stacks over
$M\an$, and the category of sections $\mycal{X}(\gb{M})$.}

\noindent
{\bf Proof.} As in Proposition~\ref{prop:linequiv}, the equivalence is given
by the composition with the forgetful $1$-morphism $\sX \to {\sf{S}}$
of stacks over $M\an$, and the proof is similar. \ \hfill $\Box$

\

\medskip

\punkt In particular, Proposition~\ref{prop:linequivgen} provides an equivalence between Definitions \ref{def:ODbimod}
and \ref{def:ODbimodrel}, and the corresponding definitions using $\cO^\times_{\sQ{M}}$-equivariant stacks.

\subsection{The $2$-category of $\star$-quantizations}\label{ssec:2stacks}

If we have an action of a group on a space, equivariant objects can be
viewed as objects on the quotient, which is usually a stack rather
than a space. Similarly, if a sheaf of groups acts on a stack,
equivariant objects can be thought of as objects on a quotient
$2$-stack. In this section, we explore this point of view in our {setting}.

\

\medskip

\punkt \
Let $G$ be a sheaf of groups acting inertially on $\sQ{M}\an$.
Let us define the quotient $2$-stack $[\sQ{M}/G]$. The definition mimics that
of the quotient $1$-stack $[N/H]$ for an action of a group $H$ on a space $N$.
Recall that $[N/H]$ is the stack of groupoids that assigns to every space $S$ 
the category $[N/H](S)$ of pairs $(T,\varphi)$, where $T$ is a $H$-torsor on $S$, 
and $\varphi : T \to N$ is a $H$-equivariant map. In the same spirit, 
$[\sQ{M}/G]$ can be defined as the $2$-stack of $G$-gerbes on $\sQ{M}$. 
This leads to the following definition.

\
\medskip

\punkt {\bf Definition.} \ \label{defn:quotient2} {\it Let $\sX$ be a 
$G$-equivariant stack over $\sQ{M}$. We say that $\sX$ is a {\bf $G$-gerbe}
if $\sX$ is a gerbe over $M$ and the action of $G$ on $\sX$ is simple and transitive
(simple transitive actions are defined in Proposition~\ref{prop:mapsfromgerbe}).

We let $[\sQ{M}/G](M)$ be the $2$-category of $G$-gerbes on $M$. More generally, let
$U\subset M$ be an open subset. Then $G$ restricts to a sheaf of groups $G|_U$ 
acting inertially on $\sQ{U}$, and we let $[\sQ{M}/G](U)$ be the $2$-category of 
$G|_U$-gerbes.}

\

\medskip

\noindent
This definition gives an abstract description of a quotient
$2$-stack. In particular cases, the quotient is identified with the $2$-stacks of
$\star$-stacks and $\star$-quantizations.
\

\smallskip

\punkt {\bf Theorem.} \label{th:2quotient}
{\it 
\begin{itemize}
\item[{\bf (a)}] The correspondence $\gb{M}\mapsto\gb{M}^\times$ gives an equivalence between
the $2$-category of $\star$-stacks on $M$ and the $2$-category of $\cO_{\sQ{M}}^\times$-gerbes 
$[\sQ{M}/\cO_{\sQ{M}}^\times](M)$. (Recall that $\gb{M}^\times$ caries a natural 
simple transitive $\cO_{\sQ{M}}^\times$-equivariant structure, which is described in \ref{ex:quant-equiv}.)

\item[{\bf (b)}] Similarly, the correspondence $\gb{M}\mapsto\gb{M}^{(1)}$ gives an 
equivalence between the $2$-category of $\star$-quantizations of $M$ and the $2$-category of $\cO_{\sQ{M}}^{(1)}$-gerbes $[\sQ{M}/\cO_{\sQ{M}}^{(1)}](M)$.
\end{itemize}}

\noindent
{\bf Proof.} Two parts of the theorem are completely analogous, so we prove only part {\bf (a)}.
Let us construct an inverse correspondence. Fix an $\cO_{\sQ{M}}^\times$-gerbe $\sX$, and
let us construct the corresponding $\star$-stack $\gb{M}$.

On the level of objects, ${\sf{ob}}(\gb{M}) = {\sf{ob}}(\sX)$. It remains to describe morphisms.
Suppose $U_1,U_2\subset M$ are open sets, $\qu{U}_i\in\sQ{M}(U_i)$, and $\alpha_i\in\sX(\qu{U}_i)$
for $i=1,2$. Let us describe $\op{Hom}_{\gb{M}}(\alpha_1,\alpha_2)$. Assume $U_1\subset U_2$, 
otherwise the space of morphisms is empty.

Consider the sheaf $\sHom_{\sX}(\alpha_{1},\alpha_{2})$ on $U_1$. 
The sheaf is a torsor over the sheaf of groups $\sHom_{\sX}(\alpha_1,\alpha_1)$ 
(because $\sX$ is a gerbe).  
For any open subset $V\subset U_1$, the action of $\cO_{\sQ{M}}^\times$ on $\sX$ induces a map
$$\Gamma(V,\cO^\times_{\qu{U}_1})\to\op{Hom}_{\sX}(\alpha_1|_V,\alpha_1|_V),$$
which is bijective because the action is simple transitive. This yields an identification
$$\cO^\times_{\qu{U}_1}=\sHom_{\sX}(\alpha_1,\alpha_1).$$ 
We finally set
\[
\op{Hom}_{\gb{M}}(\alpha_1,\alpha_2) :=  \Gamma\left( U_1,
\sHom_{\sX}(\alpha_1,\alpha_2)\times_{\cO_{\qu{U}_1}^\times}
\cO_{\qu{U}_1}\right). 
\]
Further details (such as construction of the composition of morphisms)
are left to the reader. 
\ \hfill $\Box$

\

\medskip

\punkt \ Theorem~\ref{th:2quotient} clarifies
Definitions~\ref{def:equiv.mor}, \ref{def:equiv.mor1}.  For instance,
Definition~\ref{def:equiv.mor} concerns the inertial action of
$G=\cO^\times_{\sQ{M}}$ on $\sQ{M}$. Naturally, a $G$-equivariant
stack $\sX$ over $\sQ{M}$ can be viewed as a stack over the quotient
$[\sQ{M}/G]$. We can therefore define sections of $\sX$ over objects
of $[\sQ{M}/G]$. By Theorem~\ref{th:2quotient}, the category of global
sections $[\sQ{M}/G](M)$ is equivalent to the category of
$\star$-stacks on $M$. This leads to the definition of sections of
$\sX$ over a $\star$-stack $\gb{M}$ (Definition~\ref{def:equiv.mor}).

\bibliographystyle{alpha} 
\bibliography{analytic}

\end{document}